\documentclass[12pt]{article}

\title{Set optimization - a rather short introduction}

\author{
Andreas H. Hamel\footnote{Free University Bozen, \href{mailto:andreas.hamel@unibz.it}{andreas.hamel@unibz.it}}, Frank Heyde\footnote{TU Bergakademie Freiberg, \href{mailto:frank.heyde@tu-freiberg.de}{frank.heyde@math.tu-freiberg.de}}, Andreas L\"ohne\footnote{Martin-Luther University Halle-Wittenberg, \href{mailto:andreas.loehne@mathematik.uni-halle.de}{andreas.loehne@mathematik.uni-halle.de}}, \\[.25cm] Birgit Rudloff\footnote{Princeton University, \href{mailto:brudloff@princeton.edu}{brudloff@princeton.edu}}, Carola Schrage\footnote{University of the Aosta Valley, \href{mailto:carolaschrage@gmail.com}{carolaschrage@gmail.com}}
}
\date{{\small \today}}

\usepackage{array} % for better arrays (eg matrices) in maths
\usepackage{verbatim} % adds environment for commenting out blocks of text & for better verbatim
\usepackage{amsmath, amssymb, enumerate, color}
\usepackage{stmaryrd}

\newtheorem{theorem}{Theorem}
\newtheorem{corollary}[theorem]{Corollary}
\newtheorem{remark}[theorem]{Remark}
\newtheorem{lemma}[theorem]{Lemma}
\newtheorem{definition}[theorem]{Definition}
\newtheorem{proposition}[theorem]{Proposition}
\newtheorem{example}[theorem]{Example}

\numberwithin{equation}{section}  % Formeln mit f¸hrender Kapitelnummer
\numberwithin{figure}{section}    % Abbildungen mit f¸hrender Kapitelnummer
\numberwithin{table}{section}     % Tabellen mit f¸hrender Kapitelnummer
\numberwithin{theorem}{section}

\newcommand{\of}[1]{\ensuremath{\left( #1 \right)}}

\newcommand{\abs}[1]{\ensuremath{\left| #1 \right|}}
\newcommand{\cb}[1]{\ensuremath{ \left\{ #1 \right\} }}
\newcommand{\sqb}[1]{\ensuremath{ \left[ #1 \right] }}

\newcommand{\bs}{\backslash}

\newcommand{\pend}{ \hfill $\square$ \medskip}

\newcommand{\eps}{\ensuremath{\varepsilon}}

\newcommand{\vp}{\ensuremath{\varphi}}

\renewcommand{\P}{\ensuremath{\mathcal{P}}}

\newcommand{\A}{\ensuremath{\mathcal{A}}}
\newcommand{\F}{\ensuremath{\mathcal{F}}}
\newcommand{\G}{\ensuremath{\mathcal{G}}}

\newcommand{\D}{\ensuremath{\mathcal{D}}}

\newcommand{\lel}{\preccurlyeq}
\newcommand{\leu}{\curlyeqprec}
\newcommand{\Min}{{\rm Min\,}}

\DeclareMathOperator*{\wMin}{wMin}
\DeclareMathOperator*{\Inf}{ Inf}

\newcommand{\R}{\mathrm{I\negthinspace R}}
\newcommand{\OLR}{\overline{\mathrm{I\negthinspace R}}}
\newcommand{\N}{\mathrm{I\negthinspace N}}

\newcommand{\dom}{{\rm dom \,}}
\newcommand{\epi}{{\rm epi \,}}

\newcommand{\gr}{{\rm graph \,}}

\newcommand{\cl}{{\rm cl \,}}

\newcommand{\co}{{\rm co \,}}
\newcommand{\core}{{\rm core \,}}

\newcommand{\isum}{{+^{\negmedspace\centerdot\,}}}
\newcommand{\ssum}{{+_{\negmedspace\centerdot\,}}}
\newcommand{\idif}{{-^{\negmedspace\centerdot\,}}}

\newcommand{\triup}{{\rm \vartriangle}}
\newcommand{\trido}{{\rm \triangledown}}

\newcommand{\diag}{{\rm diag}}

\newcommand{\cone}{{\rm cone\,}}
\newcommand{\Int}{{\rm int\,}}

\newcommand{\C}{\ensuremath{C^+\bs\negthinspace\cb{0}}}

\newcommand{\One}{\mathrm{1\negthickspace I}}

\usepackage[%backref=page,
colorlinks=true, linkcolor=blue, citecolor=blue, urlcolor=blue, filecolor=blue
auskommentieren und unten pdfborder aktivieren]{hyperref}
\definecolor{color0}{gray}{.50}
\definecolor{color1}{rgb}{0,.2,.8}
\definecolor{color2}{rgb}{1,.2,0}
\definecolor{color3}{rgb}{.8,.5,1}
\newcommand{\f}{\color{color1}}

\begin{document}
\maketitle

\begin{abstract}
Recent developments in set optimization are surveyed and extended including various set relations as well as fundamental constructions of a convex analysis for set- and vector-valued functions, and duality for set optimization problems. Extensive sections with bibliographical comments summarize the state of the art. Applications to vector optimization and financial risk measures are discussed along with algorithmic approaches to set optimization problems.
\end{abstract}

\vspace{.5cm}

\vspace{.5cm}

%\vspace{.5cm}

\vfill

\pagebreak

\tableofcontents

\pagebreak

%New section
\section{Introduction}

In his book \cite[p. 378]{Jahn04Book}, J. Jahn states that the set relation approach `opens a new and wide field of research' and the so-called set relations `turn out to be very promising in set optimization.' We share this opinion, and this note aims at a (partial) fulfillment of this promise.

What is ``set optimization?" The answer given in this note concerns minimization problems with set-valued objective functions and is based on a twofold solution concept: Look for a set of arguments each of which has a function value which is minimal in some sense, and all those values generate the infimum of the function. Thus, infimum attainment and minimality are the two, no longer equivalent requirements for a solution of a set optimization problem. It turns out that the set relation infimum is a useful concept in contrast to the vector order infimum which may not exist, and even if it does, it is of no practical use.

What is a motivation to consider set optimization problems? The heart of the problem is the question of how to deal with a non-total order relation, i.e. when there are non-comparable alternatives. The ``complete lattice approach" based on set relations re-gains meaningful and applicable notions of infimum and supremum even if the departing pre-ordered vector space does not have the least upper bound property, is not even a lattice, its positivity cone is not pointed, not normal or has an empty interior. The theory presented in this survey suggests that even vector-valued optimization problems should be treated as set-valued ones. This point of view has already been emphasized in \cite{Loehne11Book} for problems with a solid ordering cone.

According to an old theorem by Szpilrajn \cite{Szpilrajn30FM}, which is well-known in mathematical economics, but less so in vector and set optimization, a partial order (preorder) is the intersection of all linear orders (total preorders) including it. In the same spirit, dual descriptions of objects related to a preorder such as convex functions can be given in terms of half spaces generating total orders, hence dual objects are naturally halfspace- or hyperplane-valued.\footnote{In contrast to many duality results in vector optimization, this can bee seen as a realization of one of the many `duality principles in optimization theory that relate a problem expressed in terms of vectors in a space to a problem expressed in terms of hyperplanes in the space,' see \cite[p. 8]{Luenberger68Book}.} Since the simplest dual object is a linear functional, set-valued replacements for them should be halfspace- or hyperplane-valued as well and ``as linear as possible." This basic idea leads to a new type of duality which is not only strong enough to provide set-valued analogs of the Fenchel-Moreau and the Lagrange duality theorem, but also implies known duality results in vector optimization which are usually stated under much stronger assumptions.

It turns out that convex analysis, in particular duality, does not rely on linearity of functionals or image spaces, but rather on ``conlinearity." The structure of a conlinear space as introduced in \cite{Hamel05Habil} is precisely the part of the structure of a linear space which remains invariant under passing to the power set (with Minkowski addition) or order completion (add a least and greatest element to an ordered vector space). Thus, $\R\cup\cb{-\infty, +\infty}$ is the prototype of a conlinear space. A particular feature is the resulting bifurcation: The extended reals can be provided with inf-addition or sup-addition ({\f see \cite[p. 15]{RockafellarWets98Book}}, but already introduced by J. J. Moreau \cite{Moreau63}), which produces two different conlinear spaces. A preoder on a linear space can be extended in two different ways to the power set of this space, thus producing two different ordered conlinear spaces. It should become clear why this happens and how to deal with this ambiguity.

This survey aims at developing ideas and structures and providing a framework for principal results. Full proofs are only given for new or unpublished results, or if they illustrate an important idea particularly nicely. Sections with bibliographical remarks conclude each part with the goal to put the presented material into perspective with variational analysis and vector optimization theory in view.

Several results are new, mostly complementing those obtained by the authors in several recent publications. For example, Proposition \ref{PropTotalOrderF} discusses the totalness of set relations, Section \ref{SubSecScalarization} relies on an improved general scheme for scalarizations, Theorem \ref{ThmConlinearFunction} characterizes set-valued dual variables parallel to results for continuous linear functions and convex processes, Section \ref{SubSecTransFunc} contains a new general framework for directionally translative functions and Proposition \ref{PropLagrangeSuff} is a new sufficient condition which provides a complementary slackness condition for set optimization.

%%%New section
\section{Set relations and lattices of sets}

%%%New subsection
\subsection{Extending pre-orders from linear spaces to their power sets}

Let $Z$ be a non-trivial real linear space and $C \subseteq Z$ a convex cone with $0 \in C \neq Z$. In particular, $C = \cb{0}$ is allowed. Here, $C$ is said to be a cone if $z \in C$, $t>0$ imply $tz \in C$. By
\[
z_1 \leq_C z_2 \quad \Leftrightarrow \quad z_2 - z_1 \in C
\]
a reflexive and transitive relation $\leq_C$ is defined on $Z$; such a relation is usually called a pre-order. It is compatible with the linear structure of $Z$ in the usual sense, i.e. $z_1, z_2, z \in Z$, $t \geq 0$ and $z_1 \leq_C z_2$ imply $z_1 + z \leq_C z_2 + z$ as well as $tz_1 \leq tz_2$. Obviously,
\[
z_1 \leq_C z_2 \; \Leftrightarrow \; z_2 - z_1 \in C  \; \Leftrightarrow \; z_2 \in z_1 + C  \; \Leftrightarrow \; z_1 \in z_2 - C.
\]
The last two relationships can be used to extend $\leq_C$ from $Z$ to $\mathcal P\of{Z}$, the set of all subsets of $Z$ including the empty set $\emptyset$. Take $A, B \in \mathcal P\of{Z}$ and define
\begin{align*}
A \lel_C B \quad \Leftrightarrow \quad B \subseteq A + C, \\
A \leu_C B \quad \Leftrightarrow \quad  A \subseteq B - C. 
\end{align*}
Here and in the following, we use $+$ to denote the usual Minkowski (element-wise) addition for sets with the conventions $A + \emptyset = \emptyset + A = \emptyset$ for all $A \in \mathcal P\of{Z}$ and $A - C = A + \of{-C}$, $-C = \cb{-c \mid c \in C}$. The following facts are immediate.

\begin{proposition}
\label{PropSetRelProps}
(a) Both $\lel_C$ and $\leu_C$ are reflexive and transitive relations on $\mathcal P\of{Z}$. Moreover, they are not antisymmetric in general, and they do not coincide. 
\\
(b) $A \lel_C B$ $\Leftrightarrow$ $-B \leu_C -A$ $\Leftrightarrow$ $B \lel_{-C} A$. \\
(c) $A \lel_C B$ $\Leftrightarrow$ $A + C \supseteq B + C$;  $A \leu_C B$ $\Leftrightarrow$ $A - C \subseteq B -C$.
\end{proposition}

{\sc Proof.} Left as exercise. \pend

\medskip The property (c) above gives rise to define the set
\[
\mathcal P\of{Z, C} = \cb{A \in \mathcal P\of{Z} \mid A = A + C}
\]
and to observe that it can be identified with the set of equivalence classes with respect to the equivalence relation on $\mathcal P\of{Z}$ defined by
\begin{equation}
\label{EqEquivLel}
A \sim_C B \quad \Leftrightarrow \quad \of{A \lel_C B \; \wedge \; B \lel_C A} \quad \Leftrightarrow \quad A + C = B + C,
\end{equation}
i.e. $\sim_C$ is the symmetric part of $\lel_C$. Likewise,
\[
\mathcal P\of{Z, -C} = \cb{A \in \mathcal P\of{Z} \mid A = A - C}
\]
can be identified with the set of equivalence classes with respect to the equivalence relation
\[
A \sim_{\of{-C}} B \quad \Leftrightarrow \quad \of{A \leu_C B \; \wedge \; B \leu_C A} \quad \Leftrightarrow \quad A - C = B - C.
\]
Below, we will mainly discuss the relation $\lel_C$ which is the appropriate one when it comes to minimization; however, the theory becomes completely symmetric since every statement for the $\lel_C$ relation (and minimization) has a counterpart for $\leu_C$ (and maximization). 

The following proposition relies on (c) of Proposition \ref{PropSetRelProps}. We recall that the infimum of a subset $V \subseteq W$ of a partially ordered set $\of{W, \preceq}$ is an element $\bar w \in W$ (unique if it exists) satisfying $\bar w \preceq v$ for all $v \in V$ and $w \preceq \bar w$ whenever $w \preceq v$ for all $v \in V$. This means that the infimum is the greatest lower bound of $V$ in $W$. The infimum of $V$ is denoted by $\inf V$. Likewise, the supremum $\sup V$ is defined as the least upper bound of $V$. A partially ordered set $\of{W, \preceq}$ is called a lattice if  $\inf\cb{w_1, w_2}$ (and hence $\sup\cb{w_1, w_2}$) exists in $W$ for any two elements $w_1, w_2 \in W$. A lattice $\of{W, \preceq}$ is called (order) complete if each subset has an infimum (and hence a supremum) in $W$.

\begin{proposition}
\label{PropSetRelEC}
The pair $\of{\mathcal P\of{Z, C}, \supseteq}$ is a complete lattice. Moreover,
for a subset $\A \subseteq \mathcal P\of{Z, C}$, the infimum and the supremum of $\A$ are
given by
\begin{equation}
\label{EqInfSup} \inf \A =  \bigcup\limits_{A \in \A} A,\qquad
\sup \A = \bigcap\limits_{A \in \A} A
\end{equation}
where it is understood that $\inf \A = \emptyset$ and $\sup \A = Z$ whenever $\A = \emptyset$. The greatest (top) element of $\mathcal P\of{Z, C}$ with respect to $\supseteq$ is $\emptyset$, the least (bottom) element is $Z$. 
\end{proposition}

In particular, $\supseteq$ is a partial order on $\mathcal P\of{Z, C}$. This is remarkable since this property does not depend on the cone $C$: It can be a trivial cone, i.e. $C = \cb{0}$, or a half space, i.e. $C = \cb{z \in Z \mid \xi\of{z} \geq 0}$ where $\xi$ is a (non-zero) linear function on $Z$, i.e. an element of the algebraic dual of $Z$. Of course, a parallel result holds for $\of{\mathcal P\of{Z, -C}, \subseteq}$.

Note that the convention $\inf \emptyset = \emptyset$ and $\sup \emptyset = Z$ is in accordance with the following monotonicity property: If $\mathcal A_1 \subseteq \mathcal A_2$ then $\inf \mathcal A_1 \subseteq \inf \mathcal A_2$ and $\sup \mathcal A_1 \supseteq \sup \mathcal A_2$ in $\P\of{Z, C}$. 

\medskip {\sc Proof.} To show the first formula in \eqref{EqInfSup} one has to prove two facts: First,
\[
\forall A' \in \mathcal A \colon \bigcup\limits_{A \in \A} A \supseteq A',
\]
and secondly, for $B \in \mathcal P\of{Z, C}$
\[
\of{\forall A \in \mathcal A \colon B \supseteq A} \; \Rightarrow \; B \supseteq \bigcup\limits_{A \in \A} A.
\]
Both claims are obvious. The second formula of \eqref{EqInfSup} also follows from the definition of the supremum with respect to $\supseteq$. The lattice property is a consequence of \eqref{EqInfSup}.
\pend

\begin{remark}
One could also use other representatives of the equivalence classes defined by \eqref{EqEquivLel}
\[
\cb{A' \in \mathcal P\of{Z} \mid  A \lel_C A' \; \wedge \; A' \lel_C A}.
\]
As as rule, one has to impose additional assumptions, for example a non-empty interior of the cone $C$. An example is the infimal set approach of J. Nieuwenhuis \cite{Nieuwenhuis80Opt} which has been extended in \cite{Loehne05Diss} and \cite{LoehneTammer07Opt} (compare also \cite{Tanino88JMAA}, \cite{Loehne11Book}). This approach is summarized in Example \ref{ExSelfInfimalSets} below. 
\end{remark}

%%%New subsection
\subsection{Comments on set relations}

In the (vector and set) optimization community, D. Kuroiwa is credited for the introduction of the two "set relations" $\lel_C$ and $\leu_C$ above and, indeed, he was the first who used them for defining optimality notions for optimization problems with a set-valued objective function, compare \cite{Kuroiwa96}, \cite{KuroiwaTanakaTruong97} and several reports \cite{Kuroiwa96R, Kuroiwa97R, Kuroiwa98-1, Kuroiwa98-2, Kuroiwa98-3, Kuroiwa99RIMS} published by RIMS Kokyuroku 1996-1999. However, it should be noted that these ``set relations" were in use much earlier in different contexts. 

In the 1993 paper \cite{Brink93AU}, C. Brink describes an algebraically
motivated approach to power structures where the two relations $\leu$ and $\lel$ (analog extensions of a preorder on a general set, not necessarily a vector space) are denoted by $R_0^+$ and $R_1^+$, respectively. These and similar structures mostly defined on finite or countable sets are widely used in theoretical computer science as becomes clear from the reference list of \cite{Brink93AU}. For example, in \cite[Definition 1]{Pollak88LNM} the following definition is used: A set $A$ `can be reduced to' another set $B$ if for all $a \in A$ there is $b \in B$ such that $a \leq b$ for some partial order $\leq$, thus $A \leu B$, which is parallel to the definition of $\leu_C$ above.

Z. G. Nishianidze \cite{Nishianidze84} also used the relations $\lel$ and $\leu$ in the context of fixed point theory. This reference was brought to our attention by J. Jahn. Constructions mainly motivated by applications in economics and social choice theory can be found e.g. in \cite{Bossert89JET}, \cite{NehringPuppe96}. Compare also the references therein, especially \cite{KannaiPeleg84}. In \cite{Fishburn92JET}, set relations (on finite sets) and corresponding best choice problems are motivated by committee selection, governing coalition formation, product line formation and similar problems.

As pointed out in \cite{Hamel05Habil}, the earliest reference known to us is the paper \cite{Young31} by R. C. Young. It already contains the
definitions of $\lel_C$ and $\leu_C$ implicitly and presents applications to the analysis of upper and lower limits of sequences of real numbers.  

Another field of application for set relations is interval mathematics. In the survey \cite[Section 2.2]{Nickel75}  from 1975, an order relation is investigated which is defined on the set of order intervals of a partially ordered set $M$. This relation coincides with $\lel_C \cap \leu_C$  if $M = Z$ and $\leq_C$ is a partial order on $Z$. It has also been discussed, for example, in \cite{JahnTruong11JOTA} and \cite{Jahn13JOTA-OF}. K.-U. Jahn \cite{Jahn84} applies it in fixed point theory for interval functions, and  K. D. Schmidt \cite{Schmidt86} relates it to general ordered convex cones. Later, the ``set-less-or-equal relation" became a standard part of  the FORTRAN 95 specification for interval arithmetic, see \cite{ChiriaevWalster98}. We point out that the ``set less" relation in \cite{Jahn13JOTA-OF} actually is the ``set-less-or-equal" relation in \cite[Section 12.8]{ChiriaevWalster98} and also coincides with $\lel_C \cap \leu_C$.

In \cite{KuroiwaTanakaTruong97}, one can find a systematic investigation of six extensions of a pre-order $\leq_C$ on a topological linear
space  generated by a convex ordering cone $C$ with nonempty interior to its power set; the relations $\lel_C$ and $\leu_C$ are proven to be the only such relations which are reflexive and transitive; definitions for the convexity of set-valued functions with respect to several order relations are given. Subsequent papers of the three authors of \cite{KuroiwaTanakaTruong97} contain applications to optimization problems with  a set-valued objective, see for example \cite{Kuroiwa03JIOS}, \cite{Kuroiwa01}, \cite{Truong03JMAA}. For this topic, compare also the book \cite{Jahn04Book}, especially Chapter V. The recent paper \cite{JahnTruong11JOTA} contains even more set relations.

After 2005, many authors adopted the concepts related to $\lel_C$ and $\leu_C$, see, among an increasing number of others, \cite{AlonsoRodriguezMarin05NA}, \cite{HamelLoehne06JNCA}, \cite{HernandezRodriguezMarin07JOTA, HernandezRodriguezMarin07NA, HernandezRodriguezMarin07JMAA}, \cite{Maeda10AMC, Maeda12JOTA}, \cite{KuwanoTanakaYamada10JNCA}. Quite recently, robustness for vector optimization problems has been linked to the two (and other) set relations, see \cite{IdeKoebis13R24}, \cite{IdeSchoebel13R27}, \cite{IdeKoebisKuroiwaSchoebelTammer14FPTA}.

In \cite{Hamel05Habil}, \cite{Loehne05Diss}, \cite{LoehneTammer07Opt}, it has been realized that the set relations unfold their full potential in the framework of complete lattices; Propositions \ref{PropSetRelEC} above and infimal set versions of it such as \cite[Proposition 1.52]{Loehne11Book} may serve as a blueprint for this idea. Because of Proposition \ref{PropSetRelEC} (which can be found, even in a more general set-up, in \cite[Theorem 6]{Hamel05Habil} and, for a different image space, \cite[Proposition 1.2.3]{Loehne05Diss}) we call this approach the ``complete lattice approach" to set optimization.\footnote{For apparent reasons, we would like to call this just ``set optimization," but this term is currently used for just too many other purposes.}

%%%New subsection
\subsection{Inf-residuated conlinear spaces of sets}

We start with a definition which provides the algebraic framework for the image space analysis. It is taken from \cite{Hamel05Habil} where references and more
material about structural properties of conlinear spaces can be found.

\begin{definition}
\label{DefConlinearSpace} A nonempty set $W$ together with two algebraic operations $+
\colon W \times W \to W$ and $\cdot \colon \R_+ \times W \to W$ is called a conlinear
space provided that
\\
(C1) $\of{W, +}$ is a commutative semigroup with neutral element $\theta$,
\\
(C2) (i) $\forall w_1, w_2 \in W$, $\forall r \in \R_+$: $r \cdot \of{w_1 + w_2} = r
\cdot w_1 + r \cdot w_2$, (ii) $\forall w \in W$, $\forall r, s \in \R_+$: $s \cdot \of{r
\cdot w} = \of{sr} \cdot w$, (iii) $\forall w \in W$: $1 \cdot w = w$, (iv) $0 \cdot
\theta = \theta$.

An element $w \in W$ is called a convex element of the conlinear space $W$ if
\[
\forall s, t \geq 0 \colon \of{s+t} \cdot w = s \cdot w + t \cdot w.
\]

A conlinear space $\of{W, +, \cdot}$ together with a partial order $\preceq$ on $W$ (a
reflexive, antisymmetric, transitive relation) is called ordered conlinear space provided
that (v) $w, w_1, w_2 \in W$, $w_1 \preceq w_2$ imply $w_1 + w \preceq w_2 + w$, (vi)
$w_1, w_2 \in W$, $w_1 \preceq w_2$, $r \in \R_+$ imply $r \cdot w_1 \preceq r\cdot w_2$.

A non-empty subset $V \subseteq W$ of the conlinear space $\of{W, +, \cdot}$ is called a conlinear subspace of $W$ if (vii) $v_1, v_2 \in V$ implies $v_1 + v_2 \in V$ and (viii) $v \in V$ and $t \geq 0$ imply $t \cdot v \in V$.
\end{definition}

It can easily be checked that a conlinear subspace of a conlinear space again is a conlinear space. Note that an important feature of the above definition is the missing second distributivity law which is used to define convex elements.

\begin{example}
(a) The Minkowski addition $+$ has already been extended to $\mathcal P\of{Z, C}$ and $\mathcal P\of{Z, -C}$ (see the paragraph before Proposition \ref{PropSetRelProps}). The multiplication with non-negative numbers is extended to $\mathcal P\of{Z, C}$ by defining $t \cdot A = \cb{ta \mid a \in A}$ for $A \in \mathcal P\of{Z, C}\bs\cb{\emptyset}$, $t > 0$ and
\[
0 \cdot A = C, \quad t \cdot \emptyset = \emptyset
\]
for all $A \in \mathcal P\of{Z, C}$ and $t > 0$. In particular, $0 \cdot \emptyset = C$ by definition, and we will drop the $\cdot$ in most cases. Since the same can be done for  $\mathcal P\of{Z, -C}$, the triples $\of{\mathcal P\of{Z, C}, +, \cdot}$ and  $\of{\mathcal P\of{Z, -C}, +, \cdot}$ are conlinear spaces in the sense of Definition \ref{DefConlinearSpace}.

Note that it does not hold:
\[
\forall s, t \geq 0, \; \forall A \in \mathcal P\of{Z, C} \colon \of{s + t} \cdot A = s\cdot A + t \cdot A.
\]
Counterexamples are provided by non-convex sets $A \subseteq Z$. Therefore, $\of{\mathcal P\of{Z, C}, +, \cdot, \supseteq}$ is neither an ordered semilinear space \cite{Rubinov77RMS} nor a semi-module over the semi-ring $\R_+$ \cite{Zimmermann81Book}.

(b) The extended real numbers $\OLR = \R\cup\cb{-\infty, +\infty}$ provide two more examples. Supplied with the inf-addition $\isum$ and the sup-addition $\ssum$, respectively, one obtains two (different!) conlinear spaces. For terminology and references, see \cite{RockafellarWets98Book} and \cite{HamelSchrage12JCA}.
\end{example}

The next result connects the conlinear structure on $\of{\mathcal P\of{Z, C}, +, \cdot}$ with the order structure of $\of{\mathcal P\of{Z, C}, \supseteq}$.

\begin{proposition}
\label{PropSetRelAlgebra}
(a) $A, B, D, E \in \mathcal P\of{Z, C}$, $A \supseteq B$, $D \supseteq E$ $\Rightarrow$ $A + D \supseteq B + E$, \\
(b) $A, B \in \mathcal P\of{Z, C}$, $A \supseteq B$, $s \geq 0$  $\Rightarrow$ $sA \supseteq sB$, \\
(c) $\A \subseteq \mathcal P\of{Z, C}$, $B \in \mathcal P\of{Z, C}$ $\Rightarrow$ 
\begin{align}
\label{EqInfAdd}
\inf\of{\mathcal A + B} & = \of{\inf \A} + B \\
\label{EqSupAdd}
\sup\of{\mathcal A + B} & \supseteq \of{\sup \A} + B.
\end{align}
where $\mathcal A + B = \cb{A + B \mid A \in \A}$.
\end{proposition}

{\sc Proof.} Exercise. \pend

\medskip The following example shows that \eqref{EqSupAdd} does not hold with equality in general.

\begin{example} Let $Z = \R$, $C = \R_+$, $\A = \cb{[t, \infty) \mid t \geq 0}$, $B = \R$. Then,
\[
\forall t \geq 0 \colon [t, \infty) + \R = \R \quad \text{and} \quad \sup \A = \bigcap_{t \geq 0}[t, \infty) = \emptyset,
\]
so $\sup\of{\mathcal A + B}  = \R \neq \emptyset = \of{\sup \A} + B$.
\end{example}

\medskip Items (a) and (b) of the previous proposition show that $\of{\mathcal P\of{Z, C}, +, \cdot, \supseteq}$  (as well as  $\of{\mathcal P\of{Z, -C}, +, \cdot, \subseteq}$) carries the structure of an ordered conlinear space. Moreover, Proposition \ref{PropSetRelEC} shows that they are also complete lattices. The innocent looking equation \eqref{EqInfAdd} has far reaching consequences. In lattice theoretical terms, it means that $\of{\mathcal P\of{Z, C},+, \cdot, \supseteq}$ is {\em inf-residuated} (but not {\em sup-residuated} in general). The opposite is true for $\of{\mathcal P\of{Z, -C}, +, \cdot, \subseteq}$: it is sup-, but not inf-residuated. The following proposition is an explanation for the word ``inf-residuated."

\begin{proposition}
\label{PropInfRes}
The relationship \eqref{EqInfAdd} given in (c) of Proposition \ref{PropSetRelAlgebra} is equivalent to: For each $A, B \in \mathcal P\of{Z, C}$, the set
\[
\cb{D \in \mathcal P\of{Z, C} \mid A \supseteq B + D}
\]
has a least element (with respect to $\supseteq$).
\end{proposition}

{\sc Proof.} Assume \eqref{EqInfAdd} and fix $A, B \in \mathcal P\of{Z, C}$. Define
\[
\hat D = \inf\cb{D \in \mathcal P\of{Z, C} \mid A \supseteq B + D}.
\]
By \eqref{EqInfAdd} and \eqref{EqInfSup}, 
\begin{align*}
B + \hat D & = B +  \inf\cb{D \in \mathcal P\of{Z, C} \mid A \supseteq B + D} \\
	& = \inf\cb{B + D \in \mathcal P\of{Z, C} \mid A \supseteq B + D}\\
	& = \bigcup\cb{B + D \in \mathcal P\of{Z, C} \mid A \supseteq B + D} \subseteq A
\end{align*}
which means $\hat D \in \cb{D \in \mathcal P\of{Z, C} \mid A \supseteq B + D}$, so $\hat D$ is the desired least element. The converse direction is left as an exercise. \pend

\medskip The inf-residuation of $A, B \in \mathcal P\of{Z, C}$ is denoted
\begin{equation}
\label{EqInfRes}
A \idif B = \inf\cb{D \in  \mathcal P\of{Z, C} \mid A \supseteq B + D}.
\end{equation}
This operation will serve as a substitute for the difference in linear spaces. Indeed, for $Z = \R$, $C = \R_+$, $A = a + \R_+$, $B = b + \R_+$, $a, b \in \R$ one obtains
\[
A \idif B = \cb{r \in \R \mid b + r + \R_+ \subseteq a + \R_+} = \cb{r \in \R \mid b - a + r + \R_+ \subseteq \R_+} = a-b + \R_+.
\]
Compare Example \ref{ExRSpace} below for more about the extended reals. The following proposition states two elementary properties of $\idif$. A full calculus exists for $\idif$ which to a large extend can be derived from known results in lattice/residuation theory. One example is Proposition \ref{PropInfResCalc} below which can be understood as a residuation version of ``the negative of the infimum is the supremum of the negative."

\begin{proposition}
\label{PropBasicIdif} Let $A, B \in \mathcal P\of{Z, C}$. Then
\begin{equation}
\label{EqPointIdif}
A \idif B = \cb{z \in Z \mid B + z \subseteq A},
\end{equation}
and if $A$ is closed (convex) then $A \idif B$ is closed (convex) where $Z$ is required to be a topological linear space if closedness is involved.
\end{proposition}

{\sc Proof.} The proof of the equation is immediate from \eqref{EqInfSup} and \eqref{EqInfRes}, and the second claim follows from the first and
\[
\cb{z \in Z \mid B + z \subseteq A} = \bigcap_{b \in B}\cb{z \in Z \mid z \in A + \cb{-b}}.
\]
Of course, $A + \cb{-b}$ is closed (convex) if $A$ is closed (convex), and these properties are stable under intersection. \pend

\begin{remark}
We would like to draw the reader's attention to the fact that the structure of an ordered conlinear space which also is an inf-residuated complete lattice is ``rich enough" to serve as an image space in convex analysis. In fact, this structure is shared by $\OLR$ with inf-addition and $(\mathcal P\of{Z, C}, +, \cdot, \supseteq)$ (as well as others, see below). Completely symmetric counterparts are provided by $\OLR$ with sup-addition and $(\mathcal P\of{Z, -C}, +, \cdot, \subseteq)$ which are sup-residuated complete lattices. The transition from one to the other, provided by multiplication with $-1$, is a `duality' in the sense of \cite{Singer97Book}. Although elements of this structure were well-known and have been used before (see the comments section below), the development of this framework for a ``set-valued" convex/variational analysis and optimization theory is one contribution of the authors of this survey. One nice feature is that this structure admits to establish many results for vector/set-valued functions in the same way as for extended real-valued functions -- not surprising after one realizes the similarities between the extended reals and conlinear spaces of sets.
\end{remark}

We conclude this section by providing more examples of inf-residuated complete lattices of sets which will be used later on.

\begin{example}
\label{ExFSpace}
Let $Z$ be a topological linear space and $C \subseteq Z$ a convex cone with $0 \in C$. The set
\[
\F\of{Z, C} = \cb{A \subseteq Z \mid A = \cl\of{A + C}}
\]
clearly is a subset of $\mathcal P\of{Z, C}$, but not closed under (Minkowski) addition. Therefore, we define an associative and commutative binary operation $\oplus \colon
\F\of{Z, C} \times \F\of{Z, C} \to \F\of{Z, C}$ by
\begin{eqnarray}
\label{EqSumInP} A \oplus  B = \cl\of{A+B}
\end{eqnarray}
for $A, B \in \F\of{Z, C}$. The element-wise multiplication with 
non-negative real numbers is extended by
\[
0 \odot A = \cl C, \quad t \odot \emptyset = \emptyset
\]
for all $A \in \F\of{Z, C}$ and $t > 0$. In particular, $0 \odot \emptyset = \cl C$ by
definition, and we will drop $\odot$ in most cases. The triple $\of{\F\of{C}, \oplus, \odot}$ is a conlinear space with neutral element $\cl C$. 

On $\F\of{Z, C}$, $\supseteq$ is a partial order which is compatible with the
algebraic operations just introduced. Thus, $\of{\F\of{Z, C}, \oplus, \odot, \supseteq}$ is a partially ordered, conlinear space. 

Moreover, the pair $\of{\F\of{Z, C}, \supseteq}$ is a complete lattice, and if $\mathcal A \subseteq \F\of{Z, C}$ then
\[
 \inf_{\of{\F\of{Z, C}, \supseteq}} \mathcal A = \cl\bigcup\limits_{A \in \mathcal A} A, \quad
 \sup_{\of{\F\of{Z, C}, \supseteq}} \mathcal A = \bigcap\limits_{A \in \mathcal A} A
\]
where again $\inf_{\of{\F\of{Z, C}, \supseteq}} \mathcal A =
\emptyset$ and $\sup_{\of{\F\of{Z, C}, \supseteq}} \mathcal A = Z$ whenever $\mathcal A = \emptyset$. 

The inf-residuation in $\F\of{Z, C}$ is defined as follows: For $A, B \in
\F\of{Z, C}$, set
\begin{equation}
\label{EqInfResiduation} A \idif B = \inf_{\of{\F\of{Z, C}, \supseteq}}\cb{D \in \F\of{Z, C} \mid B + D \subseteq A} =
\cb{z\in Z \mid B + z\subseteq A}.
\end{equation}
Note that, for $A \in \mathcal \F\of{Z, C}$, the set on the right hand side of \eqref{EqInfResiduation} is indeed closed by Proposition \ref{PropBasicIdif}. 
\end{example}

\begin{example}
\label{ExSelfInfimalSets}
Let $Z$ be a topological linear space and $C \subsetneq Z$ be a closed convex cone with $\emptyset \neq \Int C \neq Z$.
The set of {\em weakly minimal} points of a subset $A \subseteq Z$ (with respect to $C$) is defined by
\[
\wMin A = \cb{y \in A \mid (\cb{y} -\Int C) \cap A = \emptyset}. 
\]
For $A \in \F(Z, C)$, it can be shown (\cite[Proposition 1.40 and Corollary 1.44]{Loehne11Book}) that $\wMin A \neq \emptyset$ if and only if $A \not\in \cb{Z, \emptyset}$. This justifies the following construction. For $A \in \F(Z, C)$, set
\[
\Inf A = \left\{
	\begin{array}{ccc}
	\wMin A & \colon & A \not\in \cb{Z, \emptyset} \\
	\cb{-\infty} & \colon & A = Z \\
	\cb{+\infty} & \colon & A = \emptyset.
	\end{array}
	\right.
\]
Then $\Inf A \subseteq Z \cup\cb{\pm\infty}$, and $\Inf A$ is never empty. The set $A$ can be reconstructed from $\Inf A$ by
\[
A = \left\{
	\begin{array}{ccc}
	\Inf A \oplus C & \colon & \Inf A \not\in \cb{\cb{-\infty}, \cb{+\infty}} \\
	  Z & \colon & \Inf A = \cb{-\infty} \\
	 \emptyset & \colon & \Inf A = \cb{+\infty}.
	\end{array}
	\right.
\]
Defining the set $\mathcal I(Z, C) = \cb{\Inf A \mid A \in \F(Z, C)}$ (the collection of `self-infimal' subsets of $Z \cup\cb{\pm\infty}$, \cite[Definition 1.50]{Loehne11Book}) and appropriate algebraic operations as well as an order one obtains $\F(Z, C) = \cb{B \oplus C \mid B \in \mathcal I(Z, C)}$. Moreover, $\mathcal I(Z, C)$ and $\F(Z, C)$ are algebraically and order isomorphic ordered conlinear spaces (compare Proposition 1.52 of \cite{Loehne11Book}). The reader is referred to \cite{Nieuwenhuis80Opt,Tanino88JMAA,LoehneTammer07Opt} and \cite{Loehne11Book} for more details concerning infimal (and supremal) sets.
\end{example}

\begin{example}
\label{ExGSpace} Let $Z$, $C$ be as in Example \ref{ExFSpace}. The set
\[
\G\of{Z, C}  = \cb{A \subseteq Z \mid A = \cl\co\of{A + C}} \subseteq \F\of{Z, C}
\]
together with the operations $\oplus$ and $\odot$ introduced in Example \ref{ExFSpace} is a conlinear subspace of $\of{\F\of{C}, \oplus, \odot}$. In fact, $\G\of{Z, C}$ precisely contains the convex elements of $\F\of{Z, C}$.  Moreover, the pair $\of{\G\of{Z, C}, \supseteq}$ is a complete lattice, and for $\emptyset \neq \mathcal A \subseteq \G\of{Z, C}$
\[
\inf_{\of{\G\of{Z, C}, \supseteq}} \mathcal A = \cl\co\bigcup\limits_{A \in \mathcal A} A
\]
whereas the formula for the supremum is the same as in $\F\of{Z,C}$. The inf-residuation in $\of{\G\of{Z, C}, \oplus, \odot, \supseteq}$ is the same as in $\of{\F\of{Z, C}, \oplus, \odot, \supseteq}$ which is a consequence of \eqref{EqInfResiduation} and Proposition \ref{PropBasicIdif}.
\end{example}

\begin{example}\label{ExGSpaceInfimal}
If in Example \ref{ExSelfInfimalSets} and under the assumptions therein, $\mathcal{F}(Z,C)$ is replaced by $\mathcal{G}(Z,C)$, we obtain a conlinear space $\mathcal{I}_{co}(Z, C)$, which is a subspace of $\mathcal{I}(Z, C)$ that is algebraically and order isomorphic to $\mathcal{G}(Z,C)$. For further details, the reader is referred to \cite[Section 1.6]{Loehne11Book}.
\end{example}

Note that parallel results are obtained for $\F\of{Z,-C}$, $\G\of{Z,-C}$ with the same algebraic operations as in $\F\of{Z, C}$, $\G\of{Z, C}$ and the order relation $\subseteq$.

\begin{example}
\label{ExRSpace}
Let us consider $Z = \R$, $C = \R_+$. Then
\[
\F\of{\R, \R_+} = \G\of{\R, \R_+} = \cb{[r, +\infty) \mid r \in \R}\cup\cb{\R}\cup\cb{\emptyset}.
\]
Moreover, by
\[
a  = \inf_{\of{\R, \leq}}A \quad \text{for} \quad A \in \G\of{\R, \R_+} \\
\]
and
\[
A  = \left\{
	\begin{array}{ccc}
	\R & : & a = -\infty \\
	\left[a, +\infty\right) & : & a \in \R \\
	\emptyset & : & a = +\infty
	\end{array}
	\right.
\]
we obtain an algebraic and order isomorphism between $\of{\G\of{\R, \R_+}, \oplus, \odot, \supseteq}$ and $\of{\OLR, \isum, \cdot, \leq}$ where $\isum$ is the inf-addition (see \cite{RockafellarWets98Book}) on $\OLR = \R\cup\cb{\pm\infty}$ with $\of{+\infty} \isum r = r \isum \of{+\infty} = +\infty$ for all $r \in \OLR$, and $\cdot$ is an extension of the multiplication of non-negative real numbers by elements of $\OLR$. Note that $0 \cdot(-\infty) = 0 \cdot (+\infty) = 0$ since otherwise $\of{\OLR, \isum, \cdot}$ is not a conlinear space. Of course, $A \supseteq B$ if, and only if, $\inf_{\of{\R, \leq}}A \leq \inf_{\of{\R, \leq}}B$. Thus, $\of{\OLR, \isum, \cdot, \leq}$ is an ordered conlinear space which is a complete lattice with respect to $\leq$. Moreover,
\[
\forall M \subseteq \OLR, \; \forall r\in \OLR \colon \inf_{(\OLR, \leq)}\of{M \isum \cb{r}} = r \isum \inf_{(\OLR, \leq)}M,
\]
which admits the introduction of the inf-residuation in $\OLR$ \cite{HamelSchrage12JCA}. Here, $M \isum \cb{r} = \cb{m \isum r \mid m \in M}$. We have
\[
r \idif s  = \inf\cb{t \in \R \mid r \leq s \isum t}
\]
for all $r, s\in\OLR$ with some strange properties. For examples, expressions like $\of{+\infty} \idif \of{-\infty}$ are well-defined and even useful as shown in \cite{HamelSchrage12JCA}, \cite{HamelSchrage13PJO}.
\end{example}

\begin{remark}
As a simple, but instructive exercise the reader should try to establish the isomorphism between $\of{\G\of{\R, -\R_+}, \oplus, \cdot, \subseteq}$ and $\of{\OLR, \ssum, \cdot, \leq}$ where $\ssum$ denotes the ``sup-addition" \cite{RockafellarWets98Book}. This shows that the reason why `there's no single symmetric way of handling $\of{+\infty} +\of{-\infty}$' is basically the same as the one for having two ``canonical" extensions of a vector pre-order to the power set of the vector space.\footnote{R. T. Rockafellar and R.-B. Wets also remark on  p. 15 of \cite{RockafellarWets98Book} that the second distributivity law does not extend to all of $\OLR$ which is another motivation for the concept of ``conlinear" spaces. Finally, it is interesting to note that the authors of \cite{RockafellarWets98Book} consider it a matter of cause to associate minimization with inf-addition (see p. 15). In the set optimization community, there is no clear consensus yet about which relation to use in what context and for what purpose. However, this note makes a clear point towards \cite{RockafellarWets98Book}: associate $\lel_C$ with minimization and $\leu_C$ with maximization because the theory works for these cases. One should have a very strong reason for doing otherwise and be advised that in this case many standard mathematical tools just don't work.}
\end{remark}

The image space $\G\of{Z, C}$ will feature prominently in duality theories for set-valued functions/optimization problems. The last example shows that it shares almost all properties with the extended reals provided with the inf-addition. The notable exception is that the order $\supseteq$ is not total in general. The following result clarifies this question.

\begin{proposition}
\label{PropTotalOrderF} Let $Z$ be a locally convex space and $C \subseteq Z$ a convex cone with $0 \in C$ and $Z \neq \cl C$. Then $\supseteq$ is total on $\F(Z, C)$ if, and only if, $\cl C$ coincides with a half-space $H^+(z^*) := \cb{z \in Z \mid z^*(z) \geq 0}$ for some $z^* \in \C$.
\end{proposition}

{\sc Proof.} The ``if" part is immediate. For the ``only if" part, assume $\supseteq$ is total on  $\F(Z, C)$ and $\cl C$ is not a half space. Then, there are $z^* \in \C$ and $\bar z \in Z$ such that 
\[
\cl C \subseteq H^+(z^*) \quad \text{and} \quad \bar z \in H^+(z^*)\bs \cl C.
\]
Indeed, the existence of $z^* \in \C$ and the first inclusion follow from a separation argument, the second from the assumption. We claim that 
\[
\forall s \in \R \colon H(z^*, s) := \cb{z \in Z \mid z^*(z) \geq s} \not\subseteq \cl C.
\]
In order to verify the claim, assume $H(z^*, s) \subseteq \cl C$ for some $s \in \R$. Then, there is $z_s \in Z$ such that $H(z^*, s) = z_s + H^+(z^*)$ and $z^*(z_s) = s$. This implies
\[
\forall t > 0 \colon z_s + t\bar z \in H(z^*, s)
\] 
since $z^*(z_s + t\bar z) = s + tz^*(\bar z) \geq s$. By assumption, $ z_s + t\bar z \in \cl C$ for all $t>0$, hence
\[
\forall t > 0 \colon \frac{1}{t}z_s + \bar z \in \cl C
\]
which in turn gives $\bar z \in \cl C$, a contradiction. This proves the claim, i.e. $\cl C \not\supseteq H(z^*, s)$ for all $s \in \R$. Since $\supseteq$ is total,
\[
\cl C \subseteq \bigcap_{s \in \R}H(z^*, s) = \emptyset,
\]
a contradiction. \pend

%%%New subsection
\subsection{Comments on conlinear spaces and residuation}

The term 'conlinear space' was coined in \cite{Hamel05Habil} because of the analogies to 'convex cones' on the one hand and to linear spaces on the other hand.

A related concept is the one of quasilinear spaces or almost linear spaces as defined in, for example, \cite{Mayer70Comp}  and \cite{Godini85JAT}, respectively. A quasilinear (or almost linear) space satisfies all the axioms of a linear space, but the second distributivity law which is required only for non-negative reals. Hence $\of{\mathcal P(Z), +, \cdot}$, $\of{\mathcal P(Z, C), +, \cdot}$ and $\of{\F(Z, C), \oplus, \cdot}$ are not quasilinear spaces in general. With respect to interval mathematics, K. D. Schmidt \cite[Section 4]{Schmidt86} observed `$[\ldots]$ it seems to be convenient to generalize one step further and to restrict the multiplication by scalars to positive scalars alone.'  Keeping all the other requirements for a quasilinear space we obtain an abstract convex cone in the sense of B. Fuchssteiner, W. Lusky \cite{FuchssteinerLusky81Book}. In \cite{KeimelRoth92LNM}, the same concept is the basic one, sometimes a convex cone even without a zero element. Abstract convex cones also coincide with semilinear spaces as probably introduced by A. M. Rubinov \cite{Rubinov77RMS} (he refers to a 1975 joint paper with S. S. Kutateladze) and recalled, for example, in \cite[Definition 2.6]{GloverJeyakumarRubinov02JCA}. 

We remark that convex cones in the sense of \cite{FuchssteinerLusky81Book} and semilinear spaces (with a ``zero") in the sense of \cite[Definition 2.6]{GloverJeyakumarRubinov02JCA} are also semi-modules over $\R_+$ (and even semivectorspaces) as defined by U. Zimmermann in  \cite[Section 5]{Zimmermann81Book}. Finally,  another close relative of a conlinear space is a semivector space in the sense of  \cite{PrakashSertel74SF}. P. Prakash and M. R. Sertel (see also \cite{PrakashSertel76PAMS}) defined this structure in the early Seventies and observed that the collections of non-empty and non-empty convex sets of a vector space form a semivector spaces. In a semivector space, the existence of a neutral element with respect to the addition is not required. Therefore, it might be considered as the ``weakest" algebraic concept discussed here.

R. Dedekind \cite{Dedekind1872Book} already introduced the residuation concept and used it in order to construct the real numbers as `Dedekind sections' of rational numbers. Among others, R. P. Dilworth and M. Ward turned it into a standard tool in abstract lattice theory, see \cite{Ward37AJM}, \cite{Dilworth38BAMS}, \cite{DilworthWard38PNAS, DilworthWard39TAMS} and many followers.

Sometimes, the right hand side of \eqref{EqInfResiduation} is called the geometric difference \cite{Pontrjagin80MS}, star-difference (for example in \cite{Volle88MP}), or Minkowski difference \cite{Hadwiger50MZ} of the two sets $A$ and $B$, and H. Hadwiger should probably be credited for its introduction. The relationship to residuation theory  (see, for instance, \cite{BlythJanowitz72Book}, \cite{GalatosJipsenKowalskiOno07Book}) has been established in \cite{HamelSchrage12JCA}. At least, we do not know an earlier reference. In the context of abstract duality, residuation has been used, for example, in \cite{MartinezLegazSinger95JCA}, \cite{GetanMartinezLegazSinger03JMS} and also in idempotent analysis (see \cite[Section 3.3]{GaubertMaxPlus97}, for example). Note that in \cite{GetanMartinezLegazSinger03JMS}, the set $\OLR$ is supplied both with $\isum$ and $\ssum$ at the same time, and this idea is extended to `the canonical enlargement' of a `boundedly complete lattice ordered group' (see \cite[Section 3]{GetanMartinezLegazSinger03JMS}) which is different from the point of view of this survey. On the other hand, $\F(Z, C)$ and $\G(Z, C)$ are special cases of $(A, \preceq)$ in \cite[Section 2]{GetanMartinezLegazSinger03JMS}, but therein the conlinear structure is not used.

%%%New section
\section{Minimality notions}

%%%New subsection
\subsection{Basic concepts}

\label{SubsecMinConcepts}

This section is concerned with the question of how to define ``infimum attainment" and ``minimality." We shall focus on the relation $\supseteq$ on $\F\of{Z, C}$ and $\mathcal G\of{Z, C}$ noting that there are parallel concepts and results for $\subseteq$ on $\mathcal F\of{Z, -C}$, $\mathcal G\of{Z, -C}$. In the remainder of the paper, the infimum or supremum is always taken in the corresponding space of elements, for example
\[
\inf \mathcal A = \cl\bigcup_{A \in \mathcal A} A
\]
whenever $\mathcal A \subseteq \F(Z, C)$ whereas for $\mathcal A \subseteq \F(Z, C)$
\[
\inf \mathcal A = \cl\co\bigcup_{A \in \mathcal A} A.
\]

With the constructions from the previous section in view, we have  at least two possibilities for a minimality notion. Given a set $\mathcal A \subseteq \mathcal F\of{Z, C}$ or $\mathcal A \subseteq \mathcal G\of{Z, C}$, look for
\begin{description}
\item[(I)] $\inf \mathcal A = \cl\bigcup_{A \in \mathcal A}A$ or  $\inf \mathcal A = \cl\co\bigcup_{A \in \mathcal A}A$, respectively, or
\item[(II)] minimal elements with respect to $\supseteq$, i.e. $B \in \mathcal A$ satisfying
\[
A \in \mathcal A, \; A \supseteq B \quad \Rightarrow \quad A = B.
\]
\end{description}
Note that the second possibility corresponds to the so-called set criterion which became popular due to the work of D. Kuroiwa and collaborators: One looks for minimal elements of $\mathcal A \subseteq \mathcal P\of{Z}$ with respect to $\lel_C$. Since $\lel_C$ is not antisymmetric in general one has to look for 
$B \in \mathcal A$ satisfying
\begin{description}
\item[(IIa)]
\[
A \in \mathcal A, \; A \lel_C B \quad \Rightarrow \quad B \lel_C A.
\]
\end{description}
Neither of the two possibilities above has been considered first. Rather, the following problem has been studied since the 1980ies by H. W. Corley \cite{Corley87JOTA}, Dinh The Luc \cite{Luc89Book} and others, and it is still popular.
\begin{description}
\item[(III)] Find minimal elements of $\bigcup_{A \in \mathcal A}A$ with respect to $\leq_C$, i.e. find $b \in \bigcup_{A \in \mathcal A}A$ satisfying
\[
a \in \bigcup_{A \in \mathcal A}A, \; a \leq_C b \quad \Rightarrow \quad b \leq_C a.
\]
\end{description}
In this way, a set optimization problem is reduced to a vector optimization problem, and sometimes this problem is referred to as the vector criterion in set optimization. Note that, in some way, it involves the infimum of $\mathcal A$ in $\mathcal P(Z, C)$.

\begin{example}\label{ExMinInf}
Let $Z = \R^2$, $C = \cb{0} \times \R_+$ and $\mathcal A = \cb{A_t \mid t \in \sqb{0,1}}$ where
\[
A_t = \sqb{-1+t, t} \times \R_+.
\]
Then each $A_t$ is minimal with respect to $\supseteq$ and 
\[
\inf \mathcal A = A_0 \cup A_1 = \sqb{-1, 1} \times \R_+.
\]
This shows that not all minimal elements are required to generate the infimum which prepares the following definition.
\end{example}

\begin{definition}
\label{DefMinSet}
Let $\mathcal A \subseteq \F(Z, C)$ or $\mathcal A \subseteq \G(Z, C)$. 

(a) A set $\mathcal B \subseteq \mathcal A$ is said to generate the infimum of $\mathcal A$ if
\[
\inf \mathcal B = \inf \mathcal A.
\]

(b) An element $\bar A \in \mathcal A$ is called minimal for $\mathcal A$ if it satisfies
\[
A \in \mathcal A, \; A \supseteq \bar A \quad \Rightarrow \quad A = \bar A.
\]
The set of all minimal elements of $\mathcal A$ is denoted by $\Min \mathcal A$. 
\end{definition}

Parallel definitions apply to generators of the supremum and maximal elements. Of course, $\mathcal A$ always generates the infimum of $\mathcal A$. On the other hand, a set of minimal elements of $\mathcal A$ does not necessarily generate the infimum of $\mathcal A$. In Example \ref{ExMinInf}, every subset of $\mathcal A$ including $A_0$ and $A_1$ consists of minimal elements and generates the infimum, i.e., in general, sets of minimal elements generating the infimum are not singletons. If $\mathcal A \subseteq \mathcal P\of{\R, \R_+}$, then a single element $A \in \mathcal A$ generates the infimum of $\mathcal A$ if, and only if, it is a minimal one. Definition \ref{DefMinSet} leads to the following ``complete lattice approach." Given a set $\mathcal A \subseteq \F(Z, C)$ or $\mathcal A \subseteq \G(Z, C)$ look for
\begin{description}
\item[(IV)] a set $\mathcal B \subseteq \mathcal A$ such that
\[
\inf \mathcal B = \inf \mathcal A \quad \text{and} \quad \mathcal B \subseteq \Min \mathcal A.
\]
\end{description}

Hence, the minimality notion of the ``complete lattice approach" consists of looking for sets of minimal elements which generate the infimum. We turn these notions into a solution concept for set optimization problems. The following definition is a special case of the general one given in \cite{HeydeLoehne11Opt}.

\begin{definition}
\label{DefSetSolution}
Let $X$ be a non-empty set, $f \colon X \to \F\of{Z, C}$ (or $f \colon X \to \G\of{Z, C}$) a function and $f[X] = \cb{f(x) \mid x \in X}$. 

(a) A set $M \subseteq X$ is called an infimizer for $f$ if
\[
\inf f[M] =  \inf f[X].
\]

(b) An element $\bar x \in X$ is called a minimizer of $f$ if $f(\bar x)$ is minimal for $f[X]$.

(c) A set $M \subseteq X$ is called a solution of the problem
\[
\tag{P} \text{minimize} \quad f(x) \quad \text{subject to} \quad x \in X
\]
if $M$ is an infimizer for $f$, and each $\bar x \in M$ is a minimizer of $f$. It is called a full solution if the set $f[M]$ includes all minimal elements of $f[X]$.
\end{definition}

Thus, solutions of set minimization problems in the ``complete lattice" sense are infimizers consisting only of minimizers. Again, parallel definitions apply to solutions of maximization problems which will later appear in duality results. One more concept is needed for a Weierstra{\ss} type theorem.

\begin{definition}
\label{DefDominationProp}
A set $\mathcal A \subseteq \F(Z, C)$ (or $\mathcal A \subseteq \G(Z, C)$) is said to satisfy the domination property if
\[
\forall A \in \mathcal A, \; \exists \bar A \in \Min\mathcal A \colon \bar A \supseteq A.
\]
\end{definition}

\begin{proposition}
\label{PropDominationProp}
Let $f \colon X \to \F\of{Z, C}$ (or $f \colon X \to \G\of{Z, C}$) be a function and $f[X]$ satisfy the domination property. Then
\[
M = \cb{x \in X \mid f(x) \in \Min f[X]}
\]
is a full solution of $\of{P}$.
\end{proposition}

{\sc Proof.} The domination property yields the first of the following inclusions while the second one follows from $M \subseteq X$:
\[
\inf_{x \in M}f(x) \supseteq \inf_{x \in X}f(x) \supseteq \inf_{x \in M}f(x).
\]
This already completes the proof since $M$ comprises all minimizers of $f[X]$. \pend

%%%New subsection
\subsection{Comments on solution concepts in set optimization}

The appearance of set-valued functions in optimization theory was mainly motivated by unifying different forms of constraints, see \cite{Borwein77MP} and also \cite{Oettli81}, \cite{Oettli82}. Problem (P) in \cite[p. 196]{Borwein81MS} seems to be the first explicit set-valued optimization problem. J. M. Borwein defines its optimal value as the infimum with respect to the underlying vector order and assumes that the image space is conditional order complete, i.e. every subset which is bounded from below (above) has an infimum (supremum) in the space. Clearly, a necessary condition for this is that the image space is a vector lattice. This restricts the applicability of such results considerably and besides, the vector infimum/supremum does not produce solution concepts which are useful in applications.

In \cite{Postolica86-1}, \cite{Postolica86-2}, V. Postolica formulates an optimization problem with a set-valued objective and uses the minimality concept (III) above.

H. W. Corley \cite{Corley87JOTA}, \cite{Corley88JOTA} defined `the maximization of a set-valued function with respect to a cone in possibly infinite dimensions'  mainly motivated by the fact that, in duality theories for multiobjective problems as established by T. Tanino and Y. Sawaragi \cite{SawaragiTanino79JOTA}, `dual problems took this form' (quotes from \cite[p.  489]{Corley87JOTA}). The same motivation can be found in Dinh The Luc's book \cite{Luc89Book} in which vector optimization problems with a set-valued objective are  investigated using the approach (III).

Both authors considered optimality in the sense of (III) above: Take the union of all objective values and then look for (weakly, properly etc.) minimal points in this union with respect to the vector order. This approach has been the leading idea ever since, among the many followers are \cite{Lin94JMAA}, \cite{Ferro96JOTA}, \cite{Ferro97SVAN}, \cite{LiChen97JMAA}, \cite{ChenJahn98MMOR}, \cite{Li99JMAA}, \cite{FloresBazan01MMOR}, \cite{LalithaDuttaGovil03JAMS} (just to mention a few), and even the more recent \cite{Truong05JMAA}, \cite{Sach05JOTA}, \cite{CrespiGinchevRocca06MMOR}, \cite{CrespiGinchevRocca06PJO}, \cite{MordukhovichTruong07CC}, \cite{MordukhovichTruong07AM}, \cite[Sections 7.1.3., 7.4.2.]{BotGradWanka09Book}, \cite{FloresBazanJimenez09SICON}, \cite{RodriguezMarinSama13JOTA} and many more. We call this approach the vector approach to set optimization.

The picture changed when the set relations were popularized by Kuroiwa and his co-authors \cite{KuroiwaTanakaTruong97} and \cite{Kuroiwa97R, Kuroiwa98-1, Kuroiwa98-2}, \cite{Kuroiwa03JIOS}. Still, it took several years until the idea to use (II) above as a solution concept for set-valued optimization problems became more popular, see \cite{HamelLoehne02R, HamelLoehne06JNCA}, \cite{AlonsoRodriguezMarin05NA}, \cite{Truong05JOTA}, \cite{HernandezRodriguezMarin07NA},  \cite{HernandezRodriguezMarin07JOTA}, \cite{HernandezRodriguezMarin07JMAA}, \cite{HernandezRodriguezMarin07PJO}, \cite{LiTeoZhang09NA} and also Chapter 5 of Jahn's book \cite{Jahn04Book}. The basic idea is, of course, to ``lift" the concept of minimal (= non-dominated) image points from elements of a vector space to elements of the power set of the vector space. Therefore, we call this approach the set relation approach to set optimization.

Roughly another ten years later, it has been realized that the so-called set relations can be utilized in a more subtle manner which is described in the previous section: Via equivalence classes with respect to the two pre-orders and hull operations one defines (conlinear) spaces of sets which enjoy rich algebraic and order structures. The set relations somehow disappear from the final picture since they serve as a tool to construct the image spaces in which the subset or superset inclusion appears as a partial order. This approach, which we call the ``complete lattice approach" to set optimization has been developed in the two theses \cite{Loehne05Diss}, \cite{Hamel05Habil} and led to the solution concept in \cite{HeydeLoehne11Opt}, \cite{Loehne11Book} which is the basis for the Definitions \ref{DefMinSet} and \ref{DefSetSolution} above. One may realize that the complete lattice approach (IV) absorbs both of (I) and (II) as well as (IIa).

%New section
\section{Set-valued functions}

%%%New subsection
\subsection{Basic concepts}

Let $X$ be another linear space and $f \colon X \to \mathcal P\of{Z, C}$ a function. The goal is to develop a convex analysis for such functions $f$. We start by recalling a popular definition. A function $\hat f \colon X \to \mathcal P\of{Z}$ is called $C$-convex (see e. g. \cite[Definition 1.1]{Borwein77MP}) if
\begin{equation}
\label{EqCConvex}
t \in \of{0,1}, \; x_1, x_2 \in X \; \Rightarrow \;
	\hat f\of{tx_1 + \of{1-t}x_2} + C \supseteq t\hat f\of{x_1} +  \of{1-t}\hat f\of{x_2},
\end{equation}
and it is called $C$-concave (\cite[p. 117]{Luc89Book}) if
\begin{equation}
\label{EqCConcave}
t \in \of{0,1}, \; x_1, x_2 \in X \; \Rightarrow \;
	t\hat f\of{x_1} +  \of{1-t}\hat f\of{x_2} \subseteq \hat f\of{tx_1 + \of{1-t}x_2} - C.
\end{equation}
Of course, the $C$-convexity inequality is just
\[
\hat f\of{t x_1 + \of{1-t}x_2} \lel_C t\hat f\of{x_1} +  \of{1-t}\hat f\of{x_2},
\]
and the $C$-concavity inequality
\[
 t\hat f\of{x_1} +  \of{1-t}f\of{x_2} \leu_C \hat f\of{tx_1 + \of{1-t}x_2}.
\]
Here is another interesting feature of the set-valued framework. If $f$ maps into $\mathcal P\of{Z, C}$, then the cone $C$ can be dropped from \eqref{EqCConvex} whereas \eqref{EqCConcave} becomes meaningless for many interesting cones $C$ (for example, for generating cones, i.e. $C - C = Z$). The opposite is true for  $\mathcal P\of{Z, -C}$-valued functions. This gives a hint why convexity (and minimization) is related to  $\mathcal P\of{Z, C}$-valued functions and concavity (and maximization) to  $\mathcal P\of{Z, -C}$-valued ones.

The graph of a function $\hat f \colon X \to \mathcal P\of{Z}$ is the set
\[
\gr \hat f = \cb{\of{x ,z} \in X \times Z \mid z \in \hat f\of{x}},
\]
and the domain is the set
\[
\dom \hat f = \cb{x \in X \mid f\of{x} \neq \emptyset}.
\]

\begin{definition}
A function $f \colon X \to \mathcal P\of{Z, C}$ is called

(a) convex if $\gr f$ is a convex subset of $X \times Z$,

(b) positively homogeneous if $\gr f$ is a cone in $X \times Z$,

(c) sublinear if $\gr f$ is a convex cone $X \times Z$,

(d) proper if $\dom f \neq \emptyset$ and $f\of{x} \neq Z$ for all $x \in X$.
\end{definition}

\begin{proposition}
\label{PropConvex1}
A function $f \colon X \to \mathcal P\of{Z, C}$ is convex if, and only if,
\begin{equation}
\label{EqConvex}
t \in \of{0,1}, \; x_1, x_2 \in X \; \Rightarrow \;
	f\of{tx_1 + \of{1-t}x_2}  \supseteq t f\of{x_1} +  \of{1-t} f\of{x_2}.
\end{equation}
It is positively homogeneous if, and only if,
\begin{equation}
\label{EqPosHom}
t>0, \; x \in X \; \Rightarrow \;
	f\of{tx}  \supseteq tf\of{x},
\end{equation}
and it is sublinear if, and only if,
\begin{equation}
\label{EqSublinear}
s, t > 0, \; x_1, x_2 \in X \; \Rightarrow \;
	 f\of{sx_1 + tx_2} \supseteq s f\of{x_1} +  t f\of{x_2}.
\end{equation}
\end{proposition}

{\sc Proof.} Exercise. \pend

\medskip A parallel result for concave $\mathcal P\of{Z, -C}$-valued functions can be established. As a straightforward consequence of Proposition \ref{PropConvex1} we obtain the following facts.

\begin{proposition}
\label{PropConvex2}
Let $f \colon X \to \mathcal P\of{Z, C}$ be a convex function. Then

(a) $f\of{x}$ is convex for all $x \in X$, i.e. $f$ is convex-valued,

(b) $\cb{x \in X \mid z \in f\of{x}}$ is convex for all $z \in Z$,

(c) $\dom f$ is convex.
\end{proposition}

{\sc Proof.} Another exercise. \pend

\medskip In the remainder of this subsection, let $X$ and $Z$ be topological linear spaces. We shall denote by $\mathcal N_X$ and $\mathcal N_Z$ a neighborhood base of $0 \in X$ and $0 \in Z$, respectively.

\begin{definition}
\label{DefClosedness}
A function $f \colon X \to \mathcal P\of{Z, C}$ is called

(a) closed-valued if $f(x)$ is a closed set for all $x \in X$,

(b) level-closed if $\cb{x \in X \mid z \in f(x)}$ is closed for all $z\in Z$,

(c) closed if $\gr f$ is a closed subset of $X \times Z$ with respect to the product topology.
\end{definition}

\begin{remark}
\label{RemLevelClosed} A function $f \colon X \to \F(Z, C)$ is level-closed if, and only if, $\cb{x \in X \mid f(x) \supseteq A}$ is closed for all $A \in \mathcal F\of{Z, C}$ which may justify the term ``level-closed." Indeed, this follows from $\cb{z} \oplus C \in \F(Z, C)$ and
\[
\forall A \in \F(Z, C) \colon \cb{x \in X \mid f(x) \supseteq A} = \bigcap_{a \in A}\cb{x \in X \mid a \in f(x)}.
\]
Level-closedness is even equivalent to closedness if $\Int C \neq \emptyset$, see \cite[Proposition 2.38]{Loehne11Book}, even for functions mapping into a completely distributive lattice as in \cite{LiuLuo91FSS}, but not in general.
\end{remark}

\begin{example}
\label{ExMusellisFunction} This example is taken from \cite[Example 3.1]{Muselli00JOTA}.
Let $X = \R$, $Z = \R^2$, $C = \cb{\of{0, t}^T \mid t \geq 0}$ and consider the function
\[
f(x) = \left\{
	\begin{array}{ccc}
	\of{\begin{array}{c}
	x \\ x+1
	\end{array}} + C & : & 0 \leq x < 1 \\[.2cm]
	\of{\begin{array}{c}
	1 \\ 4
	\end{array}} + C & : & x = 1 \\[.2cm]
	\emptyset & : & \text{otherwise}
	\end{array}
	\right.
\]
Defining sequences by 
\[
x^k = 1 - \frac{1}{k} \quad \text{and} \quad 
	z^k = \of{\begin{array}{c}
	1 - \frac{1}{k} \\ 2 - \frac{1}{k}
	\end{array}}
\]
we obtain $z^k \in f(x^k)$ for all $k = 1, 2, \ldots$, $x^k \to 1$, $z^k \to \of{1,2}^T$ and $\of{1,2}^T \not\in f(1)$, thus $\gr f$ is not closed. On the other hand,
\[
\cb{x \in X \mid z \in f(x)} = 
	\left\{
	\begin{array}{ccc}
	\cb{z_1} & : & 0 \leq z_1 < 1 \; \text{and} \; z_2 \geq z_1 + 1 \\[.15cm]
	\cb{1} & : & z_1 = 1 \; \text{and} \; z_2 \geq 4 \\[.15cm]
	\emptyset & : & \text{otherwise}
	\end{array}
	\right.
\]
thus $f$ is level-closed.
\end{example}

The following result is immediate.

\begin{proposition}
\label{PropClosedFunction}
Let $f \colon X \to \P(Z,C)$ be a closed function. Then $f$ is closed-valued and level-closed.
\end{proposition}

{\sc Proof.} Yet another exercise.
\pend

\medskip Proposition \ref{PropClosedFunction} shows that a closed $\P(Z,C)$-valued function actually maps into $\F(Z,C)$. Therefore, we can restrict the discussion of lower semicontinuity and closedness to $\mathcal F(Z,C)$-valued functions.  The following definition introduces two more related notions.

\begin{definition}
\label{DefLscLevelClosed} A function $f \colon X \to \F(Z, C)$ is called lattice-lower semicontinuous (lattice-l.s.c.) at $\bar x \in X$ iff
\begin{equation}
\label{EqLsc}
f(x) \supseteq \liminf\limits_{x \to \bar x} f(x)  = \sup\limits_{U \in \mathcal N_X}\inf\limits_{x \in \bar x + U}f(x) = 
	\bigcap\limits_{U \in \mathcal N_X}\cl\bigcup\limits_{x \in \bar x + U}f(x).
\end{equation}
It is called lattice-lower semicontinuous iff it is lattice-l.s.c. at every $\bar x \in X$.
\end{definition}

Parallel definitions apply for $\G(Z, C)$-valued functions. The next result shows the equivalence of lattice-lower semicontinuity and closedness for $\F(Z, C)$-valued functions.

\begin{proposition}
\label{PropLatticeLscClosedness}
A function $f \colon X \to \F(Z, C)$ is lattice-l.s.c. if, and only if, it is closed.
\end{proposition}

{\sc Proof.} The proof of Proposition 2.34 in \cite{Loehne11Book} also applies to this case as already discussed in \cite[p. 59]{Loehne11Book}.
\pend

\medskip The following result contains the heart of the argument for the Weierstra{\ss} type theorem.

\begin{proposition}
\label{PropCompactDom}
Let $f \colon X \to \F\of{Z, C}$ be a level-closed function such that $\dom f$ is compact. Then $f[X]$ satisfies the domination property.
\end{proposition}

{\sc Proof.} This is a special case of Proposition 2.38 in \cite{Loehne11Book}. \pend

\begin{theorem}
\label{ThmWeierstrass}
Let $f \colon X \to \F\of{Z, C}$ be a level-closed function such that $\dom f$ is compact. Then $\of{P}$ has a full solution.
\end{theorem}

{\sc Proof.}  This directly follows from Proposition \ref{PropDominationProp} and \ref{PropCompactDom}. \pend

\medskip Because of Proposition \ref{PropClosedFunction} and \ref{PropLatticeLscClosedness}, lattice-lower semicontinuity or closedness are sufficient conditions for level-closedness.

We turn to upper semi-continuity type properties which will mainly be used to establish sufficient conditions for convex duality results.

\begin{definition}
\label{DefLatticeUsc}
A function $f \colon X \to \F(Z, C)$ is called  lattice-upper semicontinuous (lattice-u.s.c.)  at $\bar x \in X$ if
\[
 \limsup\limits_{x \to \bar x} f(x) = \inf\limits_{U \in \mathcal N_X}\sup\limits_{x \in \bar x + U}f(x) = 
	\cl\bigcup\limits_{U \in \mathcal N_X}\bigcap\limits_{x \in \bar x + U}f(x) \supseteq f(\bar x).
\]
It is called lattice-upper semicontinuous (lattice-u.s.c.) if it is lattice-u.s.c. at every $x \in X$.
\end{definition}

Because of Proposition \ref{PropConvex2}, we only need to consider $\G(Z, C)$-valued functions in the following result.

\begin{proposition}
\label{PropFGSemiCont}
Let $X$ be a locally convex topological linear space and $\mathcal N_X$ a neighborhood base of $0 \in X$ consisting of convex sets. Let $f \colon X \to \of{\F(Z,C),\supseteq}$ be convex. Then, $f$ is lattice-l.s.c. (lattics-u.s.c.) at $\bar x \in X$ if, and only if, it is lattice-l.s.c. (lattice-u.s.c.) as a function into $\of{\G(Z,C),\supseteq}$ at $\bar x$.
\end{proposition}

{\sc Proof.} It is easy to prove that if $f$ is convex, then for all $x\in X$ and all $U\in \mathcal N_X$ the set $\bigcup\limits_{x \in \bar x + U} f(x)$ is convex, hence
\[
\bigcap\limits_{U\in \mathcal N_X}\cl\bigcup\limits_{x \in \bar x + U}f(x) = \bigcap\limits_{U \in \mathcal N_X}\cl\co\bigcup\limits_{x \in \bar x + U}f(x).
\]
With the definition of $\liminf$ in view, the case of lattice-lower semi-continuity follows.

Concerning lattice upper semi-continuity, take 
\[
z_1, \; z_2 \in \bigcup\limits_{U\in \mathcal N_X}\bigcap\limits_{x \in \bar x + U}f(x).
\]
Then, there are $U_1, U_2 \in \mathcal N_X$ such that $z_i \in \bigcap\limits_{x \in \bar x + U_i}f(x)$ for $i = 1,2$. Since $\mathcal N_X$ is a neighborhood base of $0 \in X$ there is $V \in \mathcal N_X$ such that $V \subseteq U_1 \cap U_2$. Hence
\[
\forall x \in \bar x + V \colon z_1, z_2 \in f\of{x}.
\]
Since $f\of{x}$ is a convex set, this implies
\[
\forall t \in \of{0,1}, \forall x\in \bar x + V \colon tz_1 + \of{1-t}z_2 \in f\of{x},
\]
hence $tz_1 + \of{1-t}z_2 \in \bigcup\limits_{U\in\mathcal U}\bigcap\limits_{x \in \bar x + U}f(x)$. This shows that the latter is a convex set. Consequently,
\[
\cl\co\bigcup\limits_{U\in\mathcal U}\bigcap\limits_{x \in \bar x + U}f(x)
= \cl\bigcup\limits_{U\in\mathcal U}\bigcap\limits_{x \in \bar x + U}f(x).
\]
The claim for lattice-upper semi-continuity follows from the definition of $\limsup$. \pend

%{\ff What should we do with this? Proposition: $\of{\bar x, \bar z} \in \Int \gr f$ $\Rightarrow$ $f$ is lattice-u.s.c. at $\bar x$. The converse is not %true in general, i.e. $f$ may be lattice-u.s.c., but $\Int \gr f = \emptyset$.}

%%%New subsection
\subsection{Scalarization of $\mathcal G\of{Z, C}$-valued functions}

\label{SubSecScalarization}

In the following, we assume that $Z$ is a non-trivial locally convex linear space with topological dual $Z^*$. For $A \subseteq Z$, define the extended real-valued functions $\sigma^\triup_A \colon Z^* \to \OLR$ and
$\sigma^\trido_A \colon Z^* \to \OLR$ by
\[
\sigma^\triup_A\of{z^*} = \inf_{a \in A}z^*\of{a} \quad \text{and} \quad
    \sigma^\trido_A\of{z^*} = \sup_{a \in A} z^*\of{a},
\]
respectively. Of course, $\sigma^\trido_A$ is the classical support function of $A$ and $\sigma^\triup_A\of{z^*} = -\sigma^\trido_A\of{-z^*}$ a version of it.
It is well-known (and a consequence of a separation argument) that $A \in \mathcal G\of{Z, C}$ if, and only if,
\begin{equation}
\label{EqDualSet}
A = \bigcap\limits_{z^*\in C^+\bs\{0\}}
    \cb{z \in  Z \mid  \sigma^\triup_A\of{z^*} \leq z^*\of{z}}.
\end{equation}
Moreover, one easily checks for $A, B  \in \mathcal G\of{Z, C}$,
\begin{equation}
\label{EqAdditiveSupport}
\forall z^*\in C^+\bs\{0\} \colon \sigma^\triup_{A \oplus B}\of{z^*} = \sigma^\triup_A\of{z^*} \isum \sigma^\triup_B\of{z^*}.
\end{equation}

\begin{lemma}
\label{LemInfSupScalarized}
If $\mathcal A \subseteq \mathcal G\of{Z, C}$ then
\begin{align}
\label{EqLemSca1}
\forall z^*\in C^+\bs\{0\} \colon & \sigma^\triup_{\inf \mathcal A}\of{z^*} = \inf\cb{\sigma^\triup_A\of{z^*} \mid A \in \mathcal A} ,\\
\label{EqLemSca2}
\forall z^*\in C^+\bs\{0\} \colon & \sigma^\triup_{\sup \mathcal A}\of{z^*} \geq \sup\cb{\sigma^\triup_A\of{z^*} \mid A \in \mathcal A}.
\end{align}
Moreover,
\begin{align}
\label{EqLemSca3}
\inf \mathcal A = \bigcap_{z^*\in C^+\bs\{0\}}\cb{z \in Z \mid \inf\cb{\sigma^\triup_A\of{z^*} \mid A \in \mathcal A} \leq z^*\of{z}}, \\
\label{EqLemSca4}
\sup \mathcal A = \bigcap_{z^*\in C^+\bs\{0\}}\cb{z \in Z \mid \sup\cb{\sigma^\triup_A\of{z^*} \mid A \in \mathcal A} \leq z^*\of{z}}
\end{align}
\end{lemma}

{\sc Proof.} If $\mathcal A \subseteq \cb{\emptyset}$ then there is nothing to prove. Otherwise, 
\[
\forall A \in\mathcal A \colon \sigma^\triup_{\inf \mathcal A}\of{z^*} = \inf_{z \in \inf\mathcal A} z^*\of{z} \leq \sigma^\triup_A\of{z^*},
\]
hence $\sigma^\triup_{\inf \mathcal A}\of{z^*} \leq \inf\cb{\sigma^\triup_A\of{z^*} \mid A \in \mathcal A}$. Conversely, 
\[
\forall z \in \bigcup_{A \in \mathcal A}A \colon  z^*\of{z} \geq \inf\cb{\sigma^\triup_A\of{z^*} \mid A \in \mathcal A},
\]
hence $\sigma^\triup_{\inf \mathcal A}\of{z^*} = \inf\cb{z^*\of{z} \mid z \in \bigcup_{A \in \mathcal A}A} \geq \inf\cb{\sigma^\triup_A\of{z^*} \mid A \in \mathcal A}$ since the support function of a set coincides with the support function of its closed convex hull. This proves \eqref{EqLemSca1} which in turn immediately implies \eqref{EqLemSca3}. 

Moreover, if $z \in \sup \mathcal A = \bigcap_{A \in \mathcal A}A$ then
\[
\forall A \in \mathcal A \colon z^*\of{z} \geq \inf_{a \in A}z^*\of{a} = \sigma^\triup_A\of{z^*} 
\]
which already proves \eqref{EqLemSca2}.
% and also
%\[
%\bigcap_{A \in \mathcal A}A \subseteq \cb{z \in Z \mid z^*\of{z} \geq \sup\cb{\sigma^\triup_A\of{z^*} \mid A \in \mathcal A}}.
%\]
Finally, for all $z^*\in C^+\bs\{0\}$
\[
\cb{z \in Z \mid z^*\of{z} \geq \sup\cb{\sigma^\triup_A\of{z^*} \mid A \in \mathcal A}} = \bigcap_{A \in \mathcal A}\cb{z \in Z \mid z^*\of{z} \geq \sigma^\triup_A\of{z^*}},
\]
hence
\begin{align*}
& \bigcap_{z^*\in C^+\bs\{0\}}\cb{z \in Z \mid z^*\of{z} \geq \sup\cb{\sigma^\triup_A\of{z^*} \mid A \in \mathcal A}}  
	 = \bigcap_{z^*\in C^+\bs\{0\}}\bigcap_{A \in \mathcal A}\cb{z \in Z \mid z^*\of{z} \geq \sigma^\triup_A\of{z^*}}  \\
	& = \bigcap_{A \in \mathcal A}\bigcap_{z^*\in C^+\bs\{0\}}\cb{z \in Z \mid z^*\of{z} \geq \sigma^\triup_A\of{z^*}} =
	\bigcap_{A \in \mathcal A} A = \sup\mathcal A
\end{align*}
according to \eqref{EqDualSet}, and this is just \eqref{EqLemSca4}. \pend

\medskip The following example shows that the inequality in \eqref{EqLemSca2} can be strict. Consider $\mathcal A = \cb{\cb{a} + \R^2_+ \mid a =(a_1, a_2)^T \in \R^2, \; a_1 \geq 0, \; a_2 \geq 0, \; a_1+a_2 = 1} \subseteq \mathcal G\of{\R^2, \R^2_+}$ and $z^* = \of{1, 1}^T$. Then $\sigma^\triup_A\of{z^*} = 1$ for all $A \in \mathcal A$ and $\sigma^\triup_{\sup \mathcal A}\of{z^*} = 2$.

The inf-residuation in $\mathcal G\of{Z, C}$ can also be represented via scalarization.

\begin{proposition}
\label{PropScalarSetDiff} For all $A, B\in \mathcal G\of{Z, C}$,
\[
A \idif B = \bigcap\limits_{z^*\in C^+\bs\{0\}}
    \cb{z\in Z \mid \sigma^\triup_A\of{z^*} \idif \sigma^\triup_B\of{z^*} \leq z^*\of{z}}.
\]
In particular, if $A = \cb{z \in Z \mid \sigma^\triup_A\of{z^*} \leq z^*(z)} \of{= A \oplus H^+(z^*)}$ for
$z^*\in C^+\bs\{0\}$, then
\begin{align*}
A \idif B = \cb{z \in Z \mid \sigma^\triup_A\of{z^*} \idif \sigma^\triup_B\of{z^*} \leq
z^*\of{z}}.
\end{align*}

Moreover,
\[
\forall z^* \in C^+\bs\{0\} \colon \sigma^\triup_{A \idif B}\of{z^*} \geq \sigma^\triup_A\of{z^*} \idif \sigma^\triup_B\of{z^*}
\]
with equality if $A = \cb{z \in Z \mid \sigma^\triup_A\of{z^*} \leq z^*(z)} \of{= A \oplus H^+(z^*)}$. 
\end{proposition}

{\sc Proof.} See \cite[Proposition 5.20]{HamelSchrage12JCA} while recalling $H^+(z^*) = \cb{z \in Z \mid z^*\of{z} \geq
0}$ for $z^* \in Z^*$.  \pend

\medskip The following result can be seen as a ``$-inf = sup -$" rule for the inf-residuation in $\G(Z, C)$. It turns out to be useful later on.

\begin{proposition}
\label{PropInfResCalc} Let $\mathcal A \subseteq \G\of{Z, C}$, $z^* \in \C$  and $H^+(z^*) = \cb{z \in Z \mid z^*(z) \geq 0}$. Then
\begin{align}
\label{EqSubzeroMinusInf}
H^+(z^*) \idif \inf \mathcal A & = \sup_{A \in \mathcal A}\sqb{H^+(z^*) \idif A}, \\
\label{EqSubzeroMinusSup}
H^+(z^*) \idif \sup \mathcal A &\supseteq \inf_{A \in \mathcal A}\sqb{H^+(z^*) \idif A}.
\end{align}
If $A \oplus H^+(z^*) = A$ for all $A \in \mathcal A$ then \eqref{EqSubzeroMinusSup} is satisfied as an equation.
\end{proposition}

{\sc Proof.} Formula \eqref{EqSubzeroMinusInf} directly follows from
\begin{multline*}
H^+(z^*) \idif \inf \mathcal A  = \cb{z \in Z \mid \cl\co\bigcup_{A \in \mathcal A}A + z \subseteq H^+(z^*)} \\
	 = \cb{z \in Z \mid \forall A \in \mathcal A \colon A + z \subseteq H^+(z^*)} 
	 = \bigcap_{A \in \mathcal A}\cb{z \in Z \mid  A + z \subseteq H^+(z^*)}.
\end{multline*} 
	The proof of \eqref{EqSubzeroMinusSup} makes use of the fact $B_1 \subseteq B_2$ $\Leftrightarrow$ $H^+(z^*) \idif B_2 \subseteq H^+(z^*) \idif B_1$. Applying it to $B_1 = \bigcap_{A \in \mathcal A}A$ and $B_1 = A$ we obtain \eqref{EqSubzeroMinusSup}. The equality case can be proven with the help of Lemma \ref{LemInfSupScalarized} and Proposition \ref{PropScalarSetDiff}. \pend

\medskip A simple counterexample for equality in \eqref{EqSubzeroMinusSup} is as follows: $Z = \R^2$, $C = \R^2_+$, $\mathcal A = \cb{A_1, A_2}$ with $A_1 = \of{1, 0}^T + \R^2_+$, $A_2 = \of{0, 1}^T + \R^2_+$ and $z^* = \of{1, 1}^T$. Both  \eqref{EqSubzeroMinusInf} and  \eqref{EqSubzeroMinusSup} are valid for more general sets than $H^+(z^*)$, but this is not needed in the following.

The previous results establish a one-to-one relationship between $\G\of{Z, C}$ and the set
\[
\Gamma\of{Z^*, C^+} = \cb{\sigma \colon C^+ \to \OLR \mid \; \sigma \; \text{is superlinear and has a closed hypograph}}.
\]
On $\Gamma\of{Z^*, C^+}$, we consider the pointwise addition $\isum$ and the pointwise multiplication with non-negative numbers $\cdot$. Finally, two elements of $\Gamma\of{Z^*, C^+}$ are compared pointwise, and we write $\sigma \leq \gamma$ whenever
\[
\forall z^* \in C^+ \colon \sigma\of{z^*} \leq \gamma\of{z^*}.
\] 
The one-to-one relationship includes the algebraic structure as well as the order structure.

\begin{proposition}
The quadrupel $\of{\Gamma\of{Z^*, C^+}, \leq, \isum, \cdot}$ is an inf-residuated conlinear space which is algebraically and order isomorphic to $\of{\G\of{Z, C}, \supseteq, \oplus, \odot}$.
\end{proposition}

{\sc Proof.} The formulas
\[
\sigma^\triup_A\of{z^*} = \inf_{a \in A}z^*\of{a}, \quad A^\triup_\sigma = \bigcap_{z^* \in C^+\bs\{0\}}\cb{z \in Z \mid \sigma\of{z^*} \leq z^*\of{z}}
\]
and
\begin{equation}
\label{EqASigma}
\sigma^\triup_{A^\triup_\sigma} = \sigma, \quad A^\triup_{\sigma^\triup_A} = A
\end{equation}
provide the relationship; the algebraic isomorphism is provided by
\[
\sigma^\triup_{A \oplus B} = \sigma^\triup_A \isum \sigma^\triup_B, \quad A^\triup_\sigma \oplus A^\triup_\gamma = A^\triup_{\sigma \isum \gamma}
\] 
and for $t \geq 0$
\[
\sigma^\triup_{tA} = t \sigma^\triup_{A}, \quad tA^\triup_\sigma = A^\triup_{t\sigma};
\]
the order isomorphism is provided by
\[
A \supseteq B \quad \Leftrightarrow \quad \sigma^\triup_A \leq \sigma^\triup_B
\]
and \eqref{EqASigma}. \pend

\begin{corollary}
\label{CorMinimizingScalarized} Let $\mathcal A \subseteq \G(Z, C)$. Then:

(a) A set $\mathcal B \subseteq \mathcal A$ generates the infimum of $\mathcal A$ if, and only if, 
\[
\sigma^\triup_{\inf \mathcal B} = \sigma^\triup_{\inf \mathcal A}. 
\]

(b) $\bar A \in \mathcal A$ is minimal for $\mathcal A$ if, and only if, $\sigma^\triup_{\bar A}$ is a minimal element of
\[
\cb{\sigma^\triup_{A} \mid A \in \mathcal A}
\]
with respect to the point-wise order in $\Gamma\of{Z^*, C^+}$,
\end{corollary}

{\sc Proof.} This is an obvious consequence of the previous results. \pend

\medskip One may think that this straightforward result reduces $\G(Z, C)$-valued (= set-valued) problems to vector optimization problem since the functions $\sigma^\triup_{A}$ could be considered as elements of some function space with point-wise order. Such an approach can be found in \cite{Jahn13JOTA-OF}. The problem with this point of view is that the functions $\sigma^\triup_{A}$ may attain (and frequently do) the values $-\infty$ and/or $+\infty$. Therefore, the difficulty is conserved by passing from $\G(Z, C)$ to $\Gamma\of{Z^*, C^+}$ since the latter is an ordered conlinear space which, in general, cannot be embedded  into a linear space of functions.

%\begin{remark}
%{\ff What should we do with this?} It may appear weird that the conlinear space $\G(Z, C)$ is one-to-one with a conlinear space of concave and %positively homogeneous functions which is supplied with inf-addition. Doesn't it? 
%\end{remark} 

We turn the above ideas into a scalarization concept for set-valued functions. Let $X$ be a topological linear space and 
$f \colon X \to \mathcal P\of{Z}$, $z^* \in C^+$ be given. Define an extended real-valued
function $\vp_{f, z^*} \colon X \to \OLR = \R\cup\cb{\pm\infty}$ by
\begin{equation}
\label{EqReScalarize} \vp_{f, z^*}\of{x} = \sigma^\triup_{f\of{x}}\of{z^*} = \inf_{z \in f\of{x}}z^*\of{z}.
\end{equation}
The new symbol $\vp_{f, z^*}$ is justified by the fact that we want to emphasize the dependence on $x$ rather than on $z^*$. From \eqref{EqDualSet} we obtain the following ``setification" formula: If $f \colon X \to \G\of{Z, C}$ then
\begin{equation}
\label{EqReSetify} \forall x \in X \colon f\of{x} =
    \bigcap_{z^* \in C^+\bs\{0\}} \cb{z \in Z \mid \vp_{f, z^*}\of{x} \leq z^*\of{z}}.
\end{equation}

Several important properties of $\G\of{Z, C}$-valued functions can equivalently be expressed using the family of its scalarizations $\cb{\vp_{f, z^*}}_{z^* \in C^+\bs\{0\}}$. One may say that, according to formula \eqref{EqReScalarize}, a $\G\of{Z, C}$-valued function is, as a mathematical object, equivalent to this family of extended real-valued functions. 

Topological properties like closedness pose difficulties in this context since scalarizations of a closed $\G\of{Z, C}$-valued function are not necessarily closed. A simple example is as follows: The function $f \colon \R \to \G(\R^2, \R^2_+)$ defined by $f(x) = \cb{\of{\frac{1}{x}, 0}^T} + \R^2_+$ for $x>0$ and $f\of{x} = \emptyset$ for $x \leq 0$ is closed and convex, but $\vp_{f, z^*}$ for $z^*=\of{0,1}^T$ is convex, but not closed. Below, we will deal with this issue. 

\begin{lemma}
\label{LemScalarizationConvex} Let $f \colon X \to \G\of{Z, C}$ be a function. Then:

(a) $f$ is convex if, and only if, $\vp_{f, z^*} \colon X \to \OLR$ is convex for all $z^* \in C^+\bs\{0\}$.

(b) $f$ is positively homogeneous if, and only if, $\vp_{f, z^*} \colon X \to \OLR$ is positively homogeneous for all $z^* \in C^+\bs\{0\}$.

(c) $f$ is (sub)additive if, and only if, $\vp_{f, z^*} \colon X \to \OLR$ is (sub)additive for all $z^* \in C^+\bs\{0\}$.

(d) $f$ is proper if, and only if, there is $z^* \in C^+\bs\{0\}$ such that $\vp_{f, z^*} \colon X \to \OLR$ is proper.

(e) $\dom f = \dom \vp_{f, z^*}$ for all $z^* \in C^+\bs\{0\}$.
\end{lemma}

{\sc Proof.} (a) "$\Rightarrow$" Take $t \in (0,1)$, $x, y \in X$ and $z^* \in C^+\bs\{0\}$. Then
\begin{align*}
\vp_{f, z^*}(tx + \of{1-t}y) 
	& = \inf_{z \in f(tx + \of{1-t}y)}z^*\of{z} 
	 \leq \inf_{z \in tf(tx) + \of{1-t}f(y)}z^*\of{z} \\
	& = \inf_{u \in tf(x)}z^*\of{u} +  \inf_{v \in \of{1-t}f(y)}z^*\of{v} \\
	& = t\inf_{\frac{u}{t} \in f(x)}\frac{z^*\of{u}}{t} + \of{1-t}\inf_{\frac{v}{\of{1-t}} \in f(y)}\frac{z^*\of{v}}{\of{1-t}}\\
	& = t\vp_{f, z^*}(x) + \of{1-t}\vp_{f, z^*}(y)
\end{align*}
where the inequality is a consequence of the convexity of $f$.

"$\Leftarrow$" By the way of contradiction, assume that $f$ is not convex. Then there are $t \in (0,1)$, $x, y \in X$, $z \in Z$ satisfying
\[
z \in tf(x) + \of{1-t}f(y), \quad z \not\in f(tx + \of{1-t}y).
\] 
Since the values of $f$ are closed convex sets we can apply a separation theorem and obtain $z^* \in C^+\bs\{0\}$ such that 
\[
z^*\of{z} < \vp_{f, z^*}(tx + \of{1-t}y) \leq t\vp_{f, z^*}(x) + \of{1-t}\vp_{f, z^*}(y)
\]
where the second inequality is a consequence of the convexity of the scalarizations. Since $f$ maps into $\G\of{Z, C}$, $z^* \in C^+\cb{0}$. Since $z \in tf(x) + \of{1-t}f(y)$ there are $u \in f\of{x}$ and $v \in f(y)$ such that $z = tu + (1-t)v$. Hence
\[
z^*\of{z} = tz^*\of{u} + (1-t)z^*\of{v} \geq t\vp_{f, z^*}(x) + (1-t)\vp_{f, z^*}(y)
\]
by definition of the scalarization. This contradicts the strict inequality above.

(c) If $f$ is (sub)additive, then (sub)additivity of the scalarizations $\vp_{f, z^*}$ follows from \eqref{EqAdditiveSupport}. The converse can be proven using the same separation idea as in the proof of (a).

The remaining claims are straightforward.
\pend

Finally, we link closedness and semicontinuity of $\mathcal G(Z,C)$-valued functions to corresponding properties of their scalarizations. The main result is Theorem \ref{ThmScalarFamily} below which shows that a proper closed and convex set-valued function and the family of its proper closed and convex scalarizations are equivalent as mathematical objects. We start with a characterization of the lattice-limit inferior in terms of scalarizations.

\begin{corollary}
\label{CorLiminfScalar}
Let $f \colon  X \to \G(Z,C)$ and $\bar x \in \dom f$ such that $f$ is lattice-l.s.c. at $\bar x$. Then
\[
\liminf_{x \to \bar x}f(x) = \cb{z \in Z \mid \forall z^* \in \C \colon \liminf_{x \to \bar x}\vp_{f, z^*}(x) \leq z^*(z)}.
\]
\end{corollary}

{\sc Proof.} Observing  $\vp_{f,z^*}(x) = \sigma^\triup_{f(x)}(z^*)$ for each $x \in X$ and  applying Lemma \ref{LemInfSupScalarized} we obtain
\begin{align*}
\sup\limits_{U \in \mathcal \mathcal N_X}\inf\limits_{x \in \bar x + U}f(x) 
	& =  \bigcap\limits_{z^*\in \C}\cb{z \in Z \mid \sup_{U \in \mathcal U} \sigma^\triup_{\inf\limits_{x \in \bar x + U}f(x)}(z^*) \leq z^*(z)} \\
	& =  \bigcap\limits_{z^*\in \C}\cb{z \in Z \mid \sup_{U \in \mathcal U} \inf_{x \in \bar x + U}\sigma^\triup_{f(x)}(z^*) \leq z^*(z)} \\
	& = \bigcap\limits_{z^*\in \C}\cb{z \in Z \mid \liminf\limits_{x \to \bar x}\vp_{f,z^*}(x) \leq z^*(z)}.
\end{align*}
Indeed, the first equality follows from the last equation in Lemma \ref{LemInfSupScalarized} applied to $\mathcal A = \cb{\inf_{x \in \bar x + U}f(x) \mid U \in \mathcal U}$ whereas the second follows from the first equation in Lemma \ref{LemInfSupScalarized} applied to $\mathcal A = \cb{f(x) \mid x \in \bar x + U}$ for $U \in \mathcal U$.
 \pend 

\begin{theorem}
\label{ThmScalarFamily}
Let $f \colon X \to \F(Z,C)$ be a function and $\dom f \neq \emptyset$. Then $f$ is closed, convex and either constant $Z$ or proper, if and only if,
\begin{align}\label{EqScalarizationProperF}
\forall x \in X \colon f(x) =\bigcap\limits_{\substack{z^*\in C^+\bs\{0\}\\ \cl\co\vp_{f,z^*} \colon X \to \OLR \; \text{is proper}}}
	\cb{z\in Z \mid \cl\co\vp_{f,z^*}(x)\leq z^*(z)},
  \end{align}
where $\cl\co\vp_{f,z^*}$ denotes the lower semi-continuous convex hull of $\vp_{f,z^*}$ defined by
\[
\epi\of{\cl\co\vp_{f,z^*}} = \cl\co\of{\epi \vp_{f,z^*}}.
\]
\end{theorem}

{\sc Proof.}
If the set $\cb{z^*\in \C \mid \cl\co \vp_{f,z^*} \colon X \to \OLR \; \text{is proper}}$ is empty, then \eqref{EqScalarizationProperF} produces $f(x)=Z$ for all $x \in X$ since $\dom f \neq \emptyset$. On the other hand, $f(x)=Z$ for all $x \in X$ implies the emptyness of the same set, hence  \eqref{EqScalarizationProperF} is satisfied in this case.

The graphs of $x \mapsto \cb{z\in Z \mid \cl\co\vp_{f,z^*}(x)\leq z^*(z)}$ are closed convex sets in $X \times Z$, and $\cb{z \in Z \mid \cl\co\vp_{f,z^*}(x)\leq z^*(z)}\neq Z$ for all $x \in X$ is true whenever $\cl\co\vp_{f,z^*}$ is proper. Thus, \eqref{EqScalarizationProperF} implies $f$ is closed, convex and either proper or constantly equal to $Z$. 

On the other hand, assume $f$ is closed, convex and proper. Then
\[
\forall x \in X \colon f(x) = \liminf\limits_{y \to x}f(y)  \neq Z,
\]
and all scalarizations are convex. Corollary \ref{CorLiminfScalar} yields
\[
f(x) = \sup\limits_{U \in \mathcal U}\inf\limits_{x \in \bar x + U}f(x)  = \bigcap\limits_{z^*\in \C}\cb{z \in Z \mid \cl\co\vp_{f,z^*}(x) \leq z^*(z)}.
\]
If $\cl\co\vp_{f,z^*}$ is improper, then $\cb{z \in Z \mid \cl\co\vp_{f,z^*}(x) \leq z^*(z)} = Z$ for all $x \in \dom f = \dom \vp_{f,z^*}$ (see \cite[Proposition 2.2.5]{Zalinescu02Book}), hence these scalarizations can be omitted from the intersection. This completes the proof.
\pend

We state a few more facts about relationships between semicontinuity properties of set-valued functions and their scalarizations.

\begin{proposition}
\label{PropSetToScalarSemiCont}
(a) If $f \colon X \to \F(Z,C)$ is lattice-u.s.c. at $\bar x \in X$, then $\vp_{f,z^*} \colon X \to \OLR$ is u.s.c. at $\bar x$ for all $z^*\in \C$.

(b) If $f \colon X \to \G(Z,C)$ is such that $\vp_{f,z^*} \colon X \to \OLR$ is l.s.c. at $\bar x \in X$ for all $z^*\in \C$, then $f$ is lattice-l.s.c. at $\bar x$.
\end{proposition}

{\sc Proof.}
(a) Define $A(\bar x) = \limsup\limits_{x \to \bar x}f(x) = \cl\bigcup\limits_{U\in\mathcal N_X}\bigcap\limits_{x \in \bar x + U}f(x)$ and take $z^* \in C^+\bs\{0\}$. By assumption,
\[
\vp_{f,z^*}(\bar x) \geq \sigma^\triup_{A(\bar x)}\of{z^*}.
\]
By a successive application of the first and the second relation of Lemma \ref{LemInfSupScalarized},
\[
\sigma^\triup_{A(\bar x)}\of{z^*} \geq \inf\limits_{U\in\mathcal N_X}\sup\limits_{x \in \bar x + U} \vp_{f, z^*}(x).
\]
This verifies the upper semicontinuity of the scalarizations.

(b) From Corollary \ref{CorLiminfScalar}, the lower semicontinuity of the $\vp_{f,z^*}$'s and \eqref{EqReSetify} we obtain
\begin{align*}
\liminf\limits_{x \to \bar x}f(x) & = \bigcap\limits_{z^*\in \C}\cb{z \in Z \mid \liminf\limits_{x \to \bar x}\vp_{f,z^*}(x) \leq z^*(z)} \\
	& \subseteq \bigcap\limits_{z^*\in \C}\cb{z \in Z \mid \vp_{f,z^*}(\bar x) \leq z^*(z)}  = f(\bar x)
\end{align*}
which means that $f$ is lattice-l.s.c. at $\bar x$.
\pend

\begin{corollary} 
\label{CorSetToScalarSemiCont} 
If $f \colon X \to \F(Z,C)$ is convex and lattice-u.s.c. at $\bar x \in \dom f$, then each scalarization $\vp_{f,z^*}$ is continuous at $\bar x$ and $f$ is lattice-l.s.c. at $\bar x$. Moreover, in this case $f$ also is lattice-u.s.c. and -l.s.c. at $\bar x$ as a function into $\G(Z,C)$.
\end{corollary}

{\sc Proof.} By Proposition \ref{PropSetToScalarSemiCont} (a), $\vp_{f,z^*}$ is u.s.c. at $\bar x$ for each $z^* \in \C$ (and also convex by Lemma \ref{LemScalarizationConvex} (a)). Well-known results about extended real-valued convex functions \cite[Theorem 2.2.9]{Zalinescu02Book} imply that  $\vp_{f,z^*}$ for each $z^* \in \C$ is continuous which in turn yields that $f$ is lattice-l.s.c. at $\bar x$ by \ref{PropSetToScalarSemiCont} (b).
The last claim follows from Proposition \ref{PropFGSemiCont}.
\pend

\begin{corollary}
\label{CorIntGraphUscScalar}
Let $f \colon X \to \G(Z,C)$ be a convex function and $\bar x \in X$ such that there exists a $\bar z \in Z$ with $(\bar x, \bar z) \in \Int(\gr f)$. Then $\vp_{f,z^*} \colon X \to \OLR$ is continuous on $\emptyset \neq \Int(\dom f)$ for all $z^*\in \C$.
\end{corollary}

{\sc Proof.}
If $(\bar x, \bar z) \in \Int(\gr f)$ then $\vp_{f,z^*}$ is bounded from above by $z^*(\bar z)$ on a neighborhood of $\bar x$, thus continuous on $\emptyset \neq \Int(\dom f)$ for all $z^*\in \C$ again by \cite[Theorem 2.2.9]{Zalinescu02Book}.
\pend

%%%New subsection
\subsection{Comments on convexity, semicontinuity and scalarization}

The properties which are called lattice-lower and lattice-upper semicontinuity can already be found in the 1978 paper \cite{Lechicki78PJO}. Note that in this survey, for obvious reasons, `upper' and `lower' are swapped compared to \cite{Lechicki78PJO}. Therein, the result of Proposition \ref{PropLatticeLscClosedness} is even referenced to a paper by Choquet from 1947.

Level-closedness features in \cite{Ferro96JOTA} and \cite{Ferro97SVAN} as `$D$-lower semi-continuity' and `$C$-lower semi-continuity', respectively: Proposition 2.3 in \cite{Ferro96JOTA} states the equivalence of (epi)closedness and level-closedness whenever the cone has a non-empty interior. The assumption ``pointedness of the cone" and a compactness assumption used in \cite{Ferro96JOTA} are not necessary, the latter already removed in \cite[Proposition 3.1]{Ferro97SVAN}. Compare also \cite{Muselli00JOTA}.

Of course, the lattice semicontinuity concepts of this survey differ from the definitions of lower and upper semicontinuity as used, for example, in \cite[Definition 1.4.1 and 1.4.2]{AubinFrankowska90Book}. This is one reason why lower and upper continuity replace lower and upper semicontinuity, respectively, in \cite{GoeRiaTamZal03Book}. We refer to Section 2.5 of \cite{GoeRiaTamZal03Book} for a survey about continuity concepts of set-valued functions and also a few bibliographical remarks at the end of the section.

For a more detailed discussion of (semi)continuity concepts for set-valued functions, compare \cite{HeydeLoehne11Opt}, \cite{Loehne11Book}, \cite{HeydeSchrage12JMAA}: Whereas Corollary \ref{CorLiminfScalar} seems to be new in this form, Proposition \ref{PropSetToScalarSemiCont} appears in \cite{HeydeSchrage12JMAA} with a (slightly) different proof.

The scalarization approach via \eqref{EqReScalarize} (and Lemma \ref{LemScalarizationConvex}) has many contributors. Motivated by economical applications, R. W. Shephard used it in \cite{Shephard70Book}, compare, for example, the definition of the `factor minimal cost function' \cite[p. 226, (97)]{Shephard70Book} and Proposition 72 on the following page where essentially Lemma \ref{LemScalarizationConvex} (a) is stated. Moreover, the first part of Proposition \ref{PropConvex1} corresponds to \cite[Appendix 2, Proposition 3]{Shephard70Book}. B. N. Pshenichnyi \cite[Lemma 1]{Pshenichnyi72CSA} also used the functions $\vp_{f, z^*}$ and proved Lemma \ref{LemScalarizationConvex} (a), see also \cite{Beresnev73CSA}. Another reference is \cite[Proposition 1.6]{PhamHuy83AMV}. In \cite{NguyenNguyen02VJM} as well as in \cite{NguyenNguyen02AMV} continuity concepts for set-valued functions are discussed using the $\vp_{f,z^*}$-functions as essential tool. See also \cite[Proposition 2.1]{BenoistPopovici03MMOR} and the more recent \cite[p. 188]{Gorokhovik08SVA} (see also the references therein).

Theorem \ref{ThmScalarFamily} has been established in \cite{Schrage09Diss}, \cite{Schrage12Opt} and is the basis for the scalarization approach to convex duality results for set-valued functions. Together with the ``setification" formula \eqref{EqReSetify} it basically tells us that one can either deal with the $\G(Z,C)$-valued function or a whole family of scalar functions, and both approaches are equivalent in the sense that major (convex duality) results can be expressed and proven either way: Using the ``set calculus" or ``scalarizations." The reader may compare the two different proofs for Lagrange duality in \cite{HamelLoehne13JOTA-OF}.

Finally, we mention that an alternative scalarization approach to (convex as well as non-convex) problems is based on directional translative 
extended real-valued functions which are used in many areas of mathematics and prominently in vector optimization, see \cite[Section 2.3]{GoeRiaTamZal03Book}. To the best of our knowledge, \cite{HamelLoehne02R} (eventually published as \cite{HamelLoehne06JNCA}) was the first generalization to set-valued problems. See also \cite{NishizawaOnodsukaTanaka05, NishizawaShimizuTanaka07NACA}, \cite{HernandezRodriguezMarin07JMAA}, \cite{LiTeoZhang09NA} and \cite{Araya12NA}, \cite{Maeda12JOTA}.

%New section
\section{Set-valued convex analysis}
\label{SecConvAnal}

What is convex analysis? A core content of this theory could be described as follows: Define affine minorants, directional derivatives, (Fenchel) conjugates and subdifferentials for convex functions and relate them by means of a Fenchel-Moreau type theorem, a max-formula, Young-Fenchel inequality as an equation. How can one establish such a theory for set-valued convex functions? In this section, we will define appropriate ``dual variables" for the set-valued framework, define ``affine minorants" of set-valued functions and introduce corresponding Fenchel conjugates, directional derivatives and subdifferentials. The difference in expressions involved in these constructions for scalar functions will be replaced by a residuation.

In the following, we assume that $X$ and $Z$ are non-trivial, locally convex, topological linear spaces with topological duals $X^*$ and $Z^*$, respectively. As before, $C \subseteq Z$ is a convex cone with $0 \in C$, and $C^+ = \cb{z^* \in Z^* \mid \forall z \in C \colon z^*\of{z} \geq 0}$
is its positive (topological) dual. %Moreover, $C^- = -C^+ =  \cb{z^* \in Z^* \mid \forall z \in C \colon z^*\of{z} \leq 0}$ is the negative dual cone.

%New subsection
\subsection{Conlinear functions}

What is an appropriate replacement for the dual variables $x^* \colon X \to \R$ in scalar convex analysis? A good guess might be to use linear operators $T \colon X \to Z$ instead of linear functionals in expressions like
\[
f^*\of{x^*} = \sup_{x \in X}\cb{x^*\of{x} - f\of{x}}.
\]
This has been done in most references about duality for vector/set optimization problems. A notable exception is the definition of the coderivative of set-valued functions due to B. S. Mordukhovich which goes back to  \cite{Mordukhovich80SMD} and can be found in \cite[Section 2]{Mordukhovich97NA}. Coderivatives at points of the graph are defined as sets of $x^*$'s depending on an element $z^* \in Z^*$. Another exception is the use of ``rank one" operators of the form $\hat zx^*$ whose existence can be proven using classical separation results, compare \cite[Proof of Theorem 4.1]{Corley87JOTA} and \cite[Theorem 4.1]{HernandezRodriguezMarin11JOTA} for an older and a more recent example. The constructions in \cite{Zowe74MathScan} are also based on this idea.

Another attempt to find set-valued analogues of linear functions is the theory of convex processes. See \cite[p. 55]{AubinFrankowska90Book} in which the authors state that `it is quite natural to regard set-valued maps, with closed convex cones as their graphs, as these set-valued analogues.'

In our approach, a class of set-valued functions will be utilized the members of which almost behave like linear functions. In some sense (see Proposition 8 in \cite{Hamel09SVVAN}), they are more general than linear operators and also than linear processes as defined in \cite[p. 55]{AubinFrankowska90Book}, and on the other hand, they form a particular class of convex processes. In fact, these functions are characterized by the fact that their graphs are homogeneous closed half spaces in $X \times Z$.

Let $x^* \in X^*$ and $z^* \in Z^*$ be  given. Define a function $S_{\of{x^*, z^*}} \colon X \to \mathcal P\of{Z}$ by
\[
S_{\of{x^*, z^*}}\of{x} = \cb{z \in Z \mid x^*\of{x} \leq z^*\of{z}}.
\]
The next result shows that these functions are indeed as ``linear" as one can hope for.

\begin{proposition}
\label{PropConlinearFunction} 
Let $\of{x^*, z^*} \in X^* \times Z^*\bs\{0\}$. Then
 
(a) for all $x \in X$ and for all $t > 0$
\[
    S_{\of{x^*, z^*}}\of{tx} = tS_{\of{x^*, z^*}}\of{x};
\]

(b) for all $x_1, x_2 \in X$
\[
S_{\of{x^*, z^*}}\of{x_1 + x_2}
    = S_{\of{x^*, z^*}}\of{x_1} + S_{\of{x^*, z^*}}\of{x_2},
\]
in particular
\[
S_{\of{x^*, z^*}}\of{x} + S_{\of{x^*, z^*}}\of{-x} = S_{\of{x^*, z^*}}\of{0} = H^+(z^*);
\]

(c) $S_{\of{x^*, z^*}}$ maps into $\G\of{Z, C}$, hence in particular into $\P\of{Z,C}$, if, and only if, $z^* \in C^+$;

(d) $S_{\of{x^*, z^*}}\of{x}$ is a closed half space with normal $z^*$ if, and only if, $z^* \neq 0$; and $S_{\of{x^*, 0}}\of{x} \in \cb{Z, \emptyset}$;

(e) if $\widehat{z} \in Z$ such that $z^*\of{\widehat{z}} = -1$ then
\begin{equation}
\label{SRepresentation} \forall x \in X \colon
    S_{\of{x^*, z^*}}\of{x} = x^*\of{x}\widehat{z} + S_{\of{x^*, z^*}}\of{0} = x^*\of{x}\widehat{z} + H^+(z^*).
\end{equation}
\end{proposition}

{\sc Proof.} Elementary, see, for instance, \cite{Hamel09SVVAN}. \pend

\medskip A function of the type $S_{\of{x^*, z^*}}$ is called {\em conlinear}. It will turn out that convex analysis is a ``conlinear" theory--not because convex functions are not linear, but because the image space of a convex function is a conlinear space and all properties of linear functions necessary for the theory are only the ``conlinear" ones from  the previous proposition. The following result gives a characterization of the class of conlinear functions in the class of all positively homogeneous and additive set-valued functions.

\begin{theorem}
\label{ThmConlinearFunction}
Let $f \colon X \rightarrow \G(Z, C)$ be a  function. Then, the following are equivalent:

(a) $\exists \of{x^*, z^*} \in X^* \times C^+\bs\{0\}$, $\forall x \in X$: $f(x) = S_{(x^*, z^*)}(x)$.

(b) $\gr f$ is a closed homogeneous half-space of $X \times Z$ and $f(0) \neq Z$.

(c) $f$ is positively homogeneous, additive, lattice-l.s.c. at $0 \in X$ and $f(0) \subseteq Z$ is a non-trivial, closed homogeneous half-space.
\end{theorem}

{\sc Proof. }
(a) $\Rightarrow$ (b), (c): Straightforward.

(b) $\Rightarrow$ (a): $\gr f$ is a closed homogenous half-space if, and only if,
\[
\exists \of{x^*, z^*} \in X^* \times Z^*\bs\{(0,0)\} \colon \gr F = \cb{\of{x,z} \in X \times Z \mid x^*(x) - z^*(z) \leq 0}.
\]
This implies
\[
\forall x \in X \colon f(x) = \cb{z \in Z \mid x^*(x) \leq z^*(z)} = S_{(x^*, z^*)}(x).
\]
Since $f(0) \neq Z$ and $f$ maps into $\G(Z, C)$, $z^* \in C^+\bs\{0\}$. By Proposition \ref{PropConlinearFunction} (b), $f$ is additive.

(c) $\Rightarrow$ (a): By assumption, $f\of{0} = H^+(z^*_0) = \cb{z \in Z \mid z^*_0\of{z} \geq 0}$ for some $z^*_0 \in C^+\bs\{0\}$. By additivity, $f\of{0} = H^+(z^*_0) = f\of{x} \oplus f\of{-x}$ for all $x \in X$, hence $f\of{x}$ is never $\emptyset$ nor $Z$. Moreover, additivity implies $f\of{x} = f\of{x + 0} = f\of{x} \oplus f\of{0} = f\of{x} \oplus H^+(z^*_0)$ for each $x \in X$. This means that every value $f\of{x}$ is a closed half space with normal $z^*_0$.

Next, we use \eqref{EqReSetify} which reads
\[
\forall x \in X \colon f\of{x} = \bigcap_{z^* \in C^+\bs\{0\}}\cb{z \in Z \mid \vp_{f, z^*}\of{x} \leq z^*(z)}.
\]
Since every value $f\of{x}$ is a half space with normal $z^*_0$ the intersection in the above formula can be replaced just by  $\cb{z \in Z \mid \vp_{f, z^*_0}\of{x} \leq z^*_0(z)}$.

We shall show that $\vp_{f, z^*_0}$ is linear. By Proposition \ref{LemScalarizationConvex} (b) and (c) it is additive because $f$ is additive, and $\vp_{f, z^*_0}(tx)=t\vp_{f, z^*}$ for $t\geq 0$, so it remains to show this for $t<0$ in order to prove homogeneity. Indeed,
\[
0 = \vp_{f, z^*_0}\of{0} = \inf_{z \in f\of{x} \oplus f\of{-x}}z^*_0(z) = \inf_{z_1 \in f\of{x}} z^*_0(z_1) + \inf_{z_2 \in f\of{-x}}z^*_0(z_2) =  \vp_{f, z^*_0}\of{x} + \vp_{f, z^*_0}\of{-x},
\]
which gives us
\[
\forall t<0 \colon \vp_{f, z^*_0}(tx)=\vp_{f, z^*_0}(-|t|x)=\abs{t}\vp_{f, z^*_0}(-x)=-\abs{t}\vp_{f, z^*_0}(x)=t\vp_{f, z^*_0}(x).
\]
Therefore, $\vp_{f, z^*_0}$ is a linear function and can be identified with some $x' \in X'$, the algebraic dual of $X$. Since $f$ is lower semicontinuous at $0 \in X$, Corollary \ref{CorLiminfScalar} with $\bar x = 0$ yields
\[
\liminf_{x \to 0}f(x) = \cb{z \in Z \mid \forall z^* \in \C \colon \liminf_{x \to 0}x'(x) \leq z^*(z)}.
\]
If $x'$ is not continuous then it is not bounded (from below) on every neighborhood $U \in \mathcal N_X$. Thus,
\[
\forall U \in \mathcal N_X \colon \inf_{x \in U} x'(x) = -\infty,
\]
hence
\[
\liminf_{x \to 0}x'(x) = \sup_{U \in \mathcal U} \inf_{x \in U}x'(x) = -\infty
\]
and consequently $Z = \liminf_{x \to 0}f(x)$ which contradicts $f(0) = H^+(z^*_0) \supseteq \liminf_{x \to 0}f(x)$. Hence, there is $x^* \in X^*$ such that $x^*(x) =  \vp_{f, z^*_0}\of{x}$ for all $x \in X$.
\pend

\medskip The basic idea for the development of a set-valued convex analysis simply is as follows: Replace the extended reals by $\G(Z, C)$, $\leq$ by $\supseteq$, use the inf/sup-formulas from Proposition \ref{EqInfSup}, replace continuous linear functionals by conlinear functions and the difference by inf-residuation. We start the program with conjugates.

%New subsection
\subsection{Fenchel conjugates of  set-valued functions}
\label{SubSecConjugates}

A crucial observation concerning Fenchel conjugates for extended real-valued functions $\vp \colon X \to \R\cup\cb{\pm\infty}$ is as follows: 
\[
r \geq \vp^*\of{x^*} \quad \Leftrightarrow \quad \forall x \in X \colon x^*\of{x} - r \leq \vp\of{x}.
\]
This means, $x^*$ belongs to the domain of $\vp^*$ precisely if there is an affine minorant of $\vp$ with ``slope" $x^*$. Replacing $x^*$ by $S_{\of{x^*, z^*}}$, $\leq$ by $\supseteq$  and recalling \eqref{EqInfRes} we obtain
\begin{align*}
\forall x \in X \colon S_{\of{x^*, z^*}}\of{x} - z \supseteq f\of{x} \quad
	& \Leftrightarrow \quad \forall x \in X \colon f\of{x} + z  \subseteq S_{\of{x^*, z^*}}\of{x} \\
	& \Leftrightarrow \quad \forall x \in X \colon z \in S_{\of{x^*, z^*}}\of{x} \idif f\of{x} \\
	&\Leftrightarrow \quad z \in \bigcap_{x \in X}\cb{S_{\of{x^*, z^*}}\of{x} \idif f\of{x}}.
\end{align*}
The function $x \mapsto S_{\of{x^*, z^*}}\of{x} - z$ is called an affine minorant of $f$ precisely if the above (equivalent) conditions are satisfied.
This discussion may justify the following definition.

\begin{definition}
\label{DefConjugate} The Fenchel conjugate of the function $f \colon X \to \mathcal P\of{Z, C}$ is $f^* \colon X^* \times C^+\bs\{0\} \to  \mathcal P\of{Z, C}$ defined by
\[
f^*\of{x^*, z^*} = \sup_{x \in X}\cb{S_{\of{x^*, z^*}}\of{x} \idif f\of{x}} = \bigcap_{x \in X}\cb{S_{\of{x^*, z^*}}\of{x} \idif f\of{x}}.
\]
The biconjugate of $f$ is $f^{**} \colon X \to \mathcal P\of{Z, C}$ defined by
\begin{align*}
f^{**}\of{x} & = \sup_{x^* \in X^*, \, z^* \in C^+\bs\{0\}}\cb{S_{\of{x^*, z^*}}\of{x} \idif f^*\of{x^*, z^*}} \\
	& = \bigcap_{x^* \in X^*, \, z^* \in C^+\bs\{0\}}\of{S_{\of{x^*, z^*}}\of{x} \idif f^*\of{x^*, z^*}}.
\end{align*}
\end{definition}

The Fenchel conjugate defined above shares most properties with her scalar little sister. 

\begin{proposition}
\label{PropPropertiesConjugate} Let $f, g \colon X \to \mathcal P\of{Z, C}$ be two functions. Then

(a) $f \supseteq g$ $\Rightarrow$ $g^* \supseteq f^*$.

(b) $f^*$ maps into $\G\of{Z, C}$, and each value of $f^*$ is a closed half space with normal $z^*$, or $\emptyset$, or $Z$.

(c) $f^{**} \supseteq f$ and $f^{**}$ is a proper closed convex function into $\G\of{Z, C}$, or $\equiv Z$, or $\equiv \emptyset$.

(d) $\of{f^{**}}^* = f^*$.

(e) For all $x \in X$, $x^* \in X^*$, $z^* \in C^+\bs\{0\}$, 
\[
f^*\of{x^*, z^*}  \subseteq S_{\of{x^*, z^*}}\of{x} \idif f\of{x} \; \Leftrightarrow \;
	f^*\of{x^*, z^*}  + f\of{x} \subseteq S_{\of{x^*, z^*}}\of{x}. 
\]
\end{proposition}

{\sc Proof.} The equivalence in (e) follows from the definition of $\idif$. The other relationships can be found in \cite{Hamel09SVVAN}, \cite{Schrage09Diss} and \cite{Schrage12Opt}. \pend

\begin{remark}
\label{RemNegConjugate}
In \cite{Hamel09SVVAN}, the ``negative conjugate" 
\[
(-f^*)\of{x^*, z^*}  = \inf_{x \in X}\cb{f\of{x} \oplus S_{\of{x^*, z^*}}\of{-x}} = \cl\bigcup_{x \in X}\cb{f\of{x} \oplus S_{\of{x^*, z^*}}\of{-x}}
\]
has been introduced which avoids the residuation. The transition from $f^*$ to $-f^*$ and vice versa can be done via 
\[
(-f^*)\of{x^*, z^*}  = H^+\of{z^*} \idif f^*\of{x^*, z^*}, \quad f^*\of{x^*, z^*} = H^+\of{z^*} \idif (-f^*)\of{x^*, z^*}
\]
using Proposition \ref{PropInfResCalc}. Sometimes, it even seems to be more natural to work with $-f^*$, for example, when it comes to Fenchel-Rockafellar duality results as presented in \cite{Hamel11Opt}.
\end{remark}

Set-valued conjugates can be expressed using the (scalar) conjugates of the scalarizing functions.

\begin{lemma}
\label{LemConjugateScalarized} If $f \colon X \to \mathcal P\of{Z, C}$, then
\begin{align}
\label{EqConScalarized}
\forall x^* \in X^*, \; \forall z^* \in \C & \colon f^*\of{x^*, z^*} = \cb{z \in Z \mid \of{\vp_{f, z^*}}^*\of{x^*} \leq z^*(z)}, \\[.3cm]
\label{EqBiconScalarized}
\forall x \in X & \colon f^{**}\of{x}  = \bigcap_{z^* \in \C}\negthinspace\cb{z \in Z \mid \of{\vp_{f, z^*}}^{**}\of{x} \leq z^*(z)}. 
\end{align}
\end{lemma}

{\sc Proof.} The first formula is a consequence of the definitions, the second follows from $\of{\vp_{f, z^*}}^{**} = \of{\vp_{f^{**}, z^*}}^{**}$ and Theorem \ref{ThmScalarFamily}. \pend

\begin{remark}
\label{RemSetificationConjugate} Conversely, $\vp_{f^*\of{\cdot, z^*}, z^*} = \of{\vp_{f, z^*}}^*$ is true (see \cite[Proposition 4.2]{Schrage12Opt} and \cite[Lemma 5.1]{HamelLoehne13JOTA-OF}. On the other hand, $\vp_{f^{**}, z^*}$ does not always coincide with $\of{\vp_{f, z^*}}^{**}$ since the latter is  a closed function which is not true for the former even if $f$ is proper closed convex (see the example before Lemma \ref{LemScalarizationConvex}).
\end{remark}

\medskip The following result is a set-valued version of the famous Fenchel-Moreau theorem. Note that the additional dual variable $z^*$ disappears via the definition of the biconjugate.

\begin{theorem}
\label{ThmFenchelMoreau}
Let $f \colon X \to \mathcal P\of{Z, C}$ be a function. Then $f = f^{**}$ if, and only if, $f$ is proper closed and convex, or identically $Z$, or identically $\emptyset$.
\end{theorem}

{\sc Proof.} This follows from Theorem \ref{ThmScalarFamily}, Lemma \ref{LemConjugateScalarized} and the classical Fenchel-Moreau theorem for scalar functions, see, for example, \cite[Theorem 2.3.3]{Zalinescu02Book}. \pend

\begin{remark}
\label{RemProofMR}
Another, more direct way to prove Theorem \ref{ThmFenchelMoreau} consists in applying the basic convex duality relationship `every closed convex set is the intersection of closed half spaces containing it' to the graph of $f$ (such half spaces are generated by pairs $(x^*, z^*) \in X^* \times C^+$), making sure that one can do without $z^* = 0$ and converting the result into formulas involving the $S_{\of{x^*, z^*}}$-functions. In this way, the scalar Fenchel-Moreau theorem is obtained as a special case. See \cite{Hamel09SVVAN} for details.
\end{remark}

\medskip To conclude this section, we point out that the Fenchel conjugate does not distinct between a function $f \colon X \to \mathcal P\of{Z, C}$ and the function
\[
\tilde f\of{x} = \cl\co f\of{x};
\]
we have $\tilde f^* = f^*$ since (compare Proposition \ref{PropBasicIdif})
\begin{multline*}
\forall x \in X \colon S_{\of{x^*, z^*}}\of{x} \idif f\of{x} = \cb{z \in Z \mid f\of{x} + z \subseteq S_{\of{x^*, z^*}}\of{x}} \\
	= \cb{z \in Z \mid \cl\co f\of{x} + z \subseteq S_{\of{x^*, z^*}}\of{x}}
	= S_{\of{x^*, z^*}}\of{x} \idif \tilde f\of{x}.
\end{multline*}
The function $\tilde f$ maps into $\G\of{Z, C}$. The above relationship means that when it comes to Fenchel conjugates it does not make a difference to start with a $\G\of{Z, C}$-valued function.

\medskip Under additional assumptions, the formulas for (bi)conjugates can be simplified. One such assumption is as follows: There is an element $\hat z \in C\bs\cb{0}$ such that
\[
\forall z^* \in C^+\bs\{0\} \colon z^*\of{\hat z} > 0.
\]
In this case, the set $B^+(\hat z) = \cb{z^* \in C^+ \mid z^*\of{\hat z} = 1}$ is a base of
$C^+$ with $0 \not\in \cl B^+(\hat z)$. That is, for each $z^* \in C^+\bs\{0\}$ there is a
unique representation $z^* = tz^*_0$ with $t>0$ and $z^*_0 \in B^+(\hat z) $. Compare
\cite{GoeRiaTamZal03Book}, Definition 2.1.14, Theorem 2.1.15 and 2.2.12 applied to $C^+$
instead of $C$. Clearly, a pointed closed convex cone with non-empty interior has a base, but, for example, the cone $L^2_+$ has an empty interior, but a base is generated by the constant 1 function.

The very definition of the functions $S_{\of{x^*, z^*}}$ gives
\[
\cb{S_{\of{x^*, z^*}} \mid x^* \in X^*, \; z^* \in C^+\bs\{0\}} =
    \cb{S_{\of{x^*, z^*}} \mid x^* \in X^*, \; z^* \in B^+(\hat z)}.
\]
Therefore, it is sufficient to run an intersection like in the definition of $f^{**}$ over $x^* \in X^*$ and $z^* \in B^+(\hat z)$. Moreover, one easily checks (see also Proposition \ref{PropConlinearFunction} (e)) for $z^* \in B^+(\hat z)$
\[
\forall x \in X \colon S_{\of{x^*, z^*}}\of{x} = \cb{x^*\of{x}\hat z} + H^+(z^*).
\]
Thus, the negative conjugate of a function $f \colon X \to \P\of{Z, C}$ can be written as
\[
(-f^*)\of{x^*, z^*} =
    \cl\bigcup_{x \in X} \sqb{f\of{x} - x^*\of{x}\hat z + H^+(z^*)} =
    \cl\bigcup_{x \in X} \sqb{f\of{x} - x^*\of{x}\hat z} \oplus H^+(z^*).
\]
The part which does not depend on $z^*$ (remember $\hat z$ defines a base of $C^+$ and is
the same for all $z^* \in C^+\bs\{0\}$) has been used in \cite{Loehne05Opt},
\cite{LoehneTammer07Opt} for the definition of another set-valued conjugate, namely
\[
(-f^*_{\hat z})\of{x^*} = \cl \bigcup_{x \in X} \sqb{f\of{x} - x^*\of{x}\hat z}.
\] 
In particular, if $Z = \R$, $C = \R_+$, then $C^+ = \R_+$, and $\cb{1}$ is a base of $C^+$, thus the intersection over the $z^*$'s disappears from the definition of $f^{**}$ and formulas like \eqref{EqBiconScalarized}.

%New subsection
\subsection{Directional derivatives}
\label{SubSecDirDer}

Usually, derivatives for set-valued functions are defined at points of their graphs as for example in \cite[Chapter 5]{AubinFrankowska90Book} and \cite[Chapter 5]{Jahn04Book}. Here, we use the inf-residuation in order to define a ``difference quotient" (which could be called ``residuation quotient") and take ``lattice limits." This leads to the concept of a lower Dini directional derivative for $\G(Z, C)$-valued functions as introduced in \cite{CrespiHamelSchrage14ArX}.

%%%CrespiHamelSchrage

\begin{definition}
\label{DefDirDer} The lower Dini directional derivative of a function $f \colon X \to \G\of{Z, C}$ with respect to $z^*\in C^+\bs\{0\}$ at $\bar x \in X$ in direction $x \in X$ is defined to be
\begin{align*}
\label{EqDirDerDef} 
f_{z^*}'\of{\bar x, x} & = \liminf_{t \downarrow 0} \frac{1}{t}\sqb{\of{f\of{\bar x + tx} \oplus H^+(z^*)} \idif f\of{\bar x}} \\
	& = \bigcap_{s>0} \cl\bigcup_{0 < t < s} \frac{1}{t}\sqb{\of{f\of{\bar x + tx} \oplus H^+(z^*)} \idif f\of{\bar x}}.
\end{align*}
\end{definition}

Obviously, $f_{z^*}' = f_{tz^*}'$ for $t>0$. Hence, if $C^+$ has a basis one only gets ``as many" directional derivatives as there are elements in the basis. 

One may ask why the set $H^+(z^*)$ appears in the definition of the difference quotient. The reason is that frequently the sets $f(\bar x +tx) \idif f(\bar x)$ and also corresponding ``lattice limits" are empty.

\begin{example}
Let $X = \R$, $Z = \R^2$, $C = \cb{\of{0, 1}^Ts \mid s \geq 0}$ and the function $f \colon X \to \G(Z, C)$ be defined by
\[
f(x) = \left\{
        \begin{array}{ccc}
        	[-x, x] \times \R_+  & : & x \in [0, 1] \\[.15cm]
        \emptyset & : & \text{otherwise}
        \end{array}
    \right. .
\]
Then, $f$ is convex and $f(1) = \inf_{x \in X} f(x) \neq Z$. However, $f(1 + tx) \idif f(1) = \emptyset$ whenever $x<0$ and $t < -\frac{1}{x}$, or $x>0$ and $t>0$. This means that the directional derivative of $f$ at $\bar x = 1$ (defined without $H^+(z^*)$) would be identically $\emptyset$. On the other hand,  $f_{z^*}'\of{1, x}$ is never empty for $z^* \in \C$ and provides much better information about the local behavior of $f$ at $\bar x = 1$.
\end{example}

For scalar functions, the standard definition of the lower Dini directional derivative can be adapted.

\begin{definition}\label{DefScalarDirDer}
The lower Dini directional derivative of a function $\vp \colon X \to \OLR$ at $\bar x$ in direction $x$ is
\begin{align*}
\vp^\downarrow(\bar x, x)&=\liminf\limits_{t\downarrow 0}\frac{1}{t}\sqb{\vp(\bar x + t x) \idif \vp(\bar x)}.
\end{align*}
\end{definition}

In Definition \ref{DefScalarDirDer}, it is neither assumed $\bar x \in \dom \vp$, nor $\vp$ be a proper function. This is possible since the difference operator is replaced by the residual operator. For $\G\of{Z, C}$-valued functions, the lower Dini directional derivative can be expressed by corresponding derivatives of scalarizations.

\begin{proposition}
\label{PropDiniSca}
(a) For all $\bar x \in X$, for all $x \in X$,
\begin{align}
\label{EqDirDerScalarized}
f^\downarrow_{z^*}(\bar x, x) & = \cb{z \in Z \mid \vp_{f, z^*}^\downarrow(\bar x, x) \leq -z^*\of{z}} \\
\vp_{f,z^*}^\downarrow(\bar x, x) & = \vp_{f^\downarrow_{z^*}\of{\bar x, \cdot}, z^*}\of{x}.
\end{align}
\end{proposition}

{\sc Proof. } See \cite[Proposition 3.4]{CrespiHamelSchrage14ArX}. \pend 

The next result is familiar in the scalar case for proper functions, see \cite[Theorem 2.1.14]{Zalinescu02Book}.

\begin{lemma}
\label{LemDirDerSublinear} Let $f \colon X \to \G\of{Z,C}$ be convex, $\bar x \in X$ and
$z^* \in C^+\bs\{0\}$. Then
\begin{equation}\label{EqInfDirDer}
\forall x \in X \colon f'_{z^*}\of{\bar x, x} =  \inf_{t>0}\frac{1}{t}\sqb{\of{f\of{\bar x + tx} \oplus H^+(z^*)} \idif f\of{\bar x}},
\end{equation}
and the function
\[
x \mapsto f'_{z^*}\of{x_0, x}
\]
is sublinear as a function from $X$ into $\G\of{Z,C}$. If $\bar x \in \dom f$, then $\dom f'_{z^*}\of{\bar x,
\cdot} = \cone\of{\dom f - \bar x}$. Moreover,
\[
f'_{z^*}\of{\bar x, 0} =
    \left\{
        \begin{array}{ccc}
        H^+(z^*) & : & f\of{\bar x} \oplus H^+(z^*) \not\in \cb{Z, \emptyset} \\
        Z & : & f\of{\bar x} \oplus H^+(z^*) \in \cb{Z, \emptyset}
        \end{array}
    \right. .
\]
\end{lemma}

{\sc Proof.} It relies on the monotonicity of the ``residuation quotient" 
\[
\frac{1}{t}\sqb{\of{f\of{\bar x + tx} \oplus H^+(z^*)} \idif f\of{\bar x}} 
\]
which in turn is proven using a calculus for the inf-residuation and the convexity of $f$. For details, compare \cite{HamelSchrage13PJO}. \pend

\medskip The following result tells us when the directional derivative has only ``finite" values. As usual, we denote by $\core M$ the algebraic interior of a set $M \subseteq X$.

\begin{theorem}
\label{ThmDirDerFinite} Let $f \colon X \to \G\of{Z, C}$ be convex and $\bar x \in \core\of{\dom f}$. If $f$ is proper, then there exists $z^* \in C^+\bs\{0\}$ such that $f'_{z^*}\of{\bar x, x} \not\in \cb{Z, \emptyset}$ for all $x \in X$.
\end{theorem}

{\sc Proof.} See \cite{HamelSchrage13PJO}. \pend

%New subsection
\subsection{The subdifferential}
\label{SubSecSubdifferential}

For convex functions, we define elements of the subdifferential using conlinear minorants of the sublinear directional derivative.

\begin{definition}
\label{DefSubdiff} Let $f \colon X \to \G\of{Z,C}$ be convex, $\bar x \in X$ and $z^* \in C^+\bs\{0\}$. The set
\[
\partial f_{z^*}\of{\bar x} = \cb{x^* \in X^* \mid \forall x \in X \colon
    S_{\of{x^*, z^*}}\of{x} \supseteq f_{z^*}'\of{\bar x, x}}
\]
is called the $z^*$-subdifferential of $f$ at $\bar x$.
\end{definition}

Again, the basic idea is to replace a continuous linear functional $x^*$ by $S_{\of{x^*, z^*}}$. An alternative characterization of the subdifferential is provided in the following result.

\begin{proposition}
\label{PropSubdiffInequality} Let $f \colon X \to \G\of{Z,C}$ be convex and $\bar x \in X$. The following statements are equivalent for $x^* \in X^*$, $z^* \in C^+\bs\{0\}$:

(a) $\forall x \in X$: $S_{\of{x^*, z^*}}\of{x} \supseteq f'_{z^*}\of{\bar x, x}$,

(b) $\forall x \in X$: $S_{\of{x^*, z^*}}\of{x - \bar x} \supseteq \of{f\of{x} \oplus H^+(z^*)} \idif f\of{\bar x}$.
\end{proposition}

{\sc Proof.} See \cite{HamelSchrage13PJO}. \pend

Vice versa, under some ``regularity" one can reconstruct the directional derivative from the subdifferential. This result is known as the max-formula. Here is a set-valued version.

\begin{theorem}
\label{ThmMaxFormula} Let $f \colon X \to  \G\of{Z,C}$ be a convex function, $\bar x \in \dom
f$ and $z^* \in C^+\bs\{0\}$ such that the function $x \mapsto  f\of{x} \oplus H^+(z^*)$
is proper and the function $\vp_{f, z^*} \colon X \to \R\cup\cb{+\infty}$ is upper semi-continuous
at $\bar x$. Then $\partial f_{z^*}\of{\bar x} \neq \emptyset$ and it holds
\begin{equation}
\label{EqMaxFormula} \forall x \in X \colon f'_{z^*}\of{\bar x, x} =
    \bigcap_{x^* \in \partial f_{z^*}\of{\bar x}} S_{\of{x^*, z^*}}\of{x}.
\end{equation}
Moreover, for each $x \in X$ there exists $\bar x^* \in \partial f_{z^*}\of{\bar x}$ such that
\begin{equation}
\label{EqMaxFormulaMax} f'_{z^*}\of{\bar x, x} = S_{\of{\bar x^*, z^*}}\of{x}.
\end{equation}
\end{theorem}

{\sc Proof.} See \cite{HamelSchrage13PJO}. \pend

\medskip Next, we link the subdifferential and the Fenchel conjugate.

\begin{proposition}
\label{PropSubdiffConjugate} Let $f \colon X \to \G\of{Z, C}$ be convex, $\bar x \in X$, $\dom f \neq \emptyset$ and $f\of{\bar x} \oplus H^+(z^*) \neq Z$. Then, the following statements are equivalent for $x^* \in X^*$, $z^* \in C^+\bs\{0\}$:

(a) $x^* \in \partial f_{z^*}\of{\bar x}$,

(b) $\forall x \in X$: $S_{\of{x^*, z^*}}\of{x} \idif f\of{x} \supseteq S_{\of{x^*, z^*}}\of{\bar x} \idif f\of{\bar x}$.
\end{proposition}

{\sc Proof.} See \cite{HamelSchrage13PJO}. \pend

\medskip This results basically says that $x^* \in \partial f_{z^*}\of{\bar x}$ if the supremum in the definition of the conjugate is attained at $\bar x$ since from the Young-Fenchel inequality we have
\[
S_{\of{x^*, z^*}}\of{\bar x} \idif f\of{\bar x}  \supseteq  f^*\of{x^*, z^*}
\]
whereas (b) above produces
\[
f^*\of{x^*, z^*} =  \bigcap_{x \in X}\cb{S_{\of{x^*, z^*}}\of{x} \idif f\of{x}} \supseteq S_{\of{x^*, z^*}}\of{\bar x} \idif f\of{\bar x}.
\]
This means: In the sense of Definition \ref{DefSetSolution} adapted to maximization, the set $\cb{\bar x}$ is a solution of the problem
\[
\text{maximize} \quad S_{\of{x^*, z^*}}\of{x} \idif f\of{x} \quad \text{over} \quad x \in X.
\]

Finally, we want to describe the set of points satisfying the condition $0 \in \partial_{z^*}f\of{\bar x}$.

\begin{proposition}
\label{PropZeroInSubdiff} Let $f \colon X \to \G\of{Z, C}$ be convex, $z^* \in \C$ and $\bar x \in \dom f$ such that $f(\bar x) \oplus H^+(z^*) \neq Z$. Then, the following statements are equivalent:

(a) $H^+(z^*) \supseteq f'_{z^*}\of{\bar x, x}$ for all $x \in X$,

(b) $0 \in \partial f_{z^*}\of{\bar x}$,

(c) $f(\bar x) \oplus H^+(z^*) = \sqb{\inf_{x \in X}f(x)} \oplus H^+(z^*)$,

(d) $\vp_{f,z^*}(\bar x) \leq \vp_{f,z^*}(x)$ for all $x \in X$.
\end{proposition}

{\sc Proof.} This is immediate from the previous results. \pend

\medskip We will call an $\bar x \in X$ for which there is $z^* \in \C$ satisfying (c) in Proposition \ref{PropZeroInSubdiff}  a $C^+$-minimizer of problem (P) in Definition \ref{DefSetSolution}. The question arises if there is a (full) solution of (P) consisting of $C^+$-minimizers and how such a solution can be characterized.

We conclude this section by noting that a calculus for the $z^*$-subdifferential can be derived from corresponding calculus rules for extended real-valued convex functions. The additional feature in the set-valued case is the dependence of $\partial f_{z^*}\of{\bar x}$ on $z^*$, i.e. properties of the mapping $z^* \mapsto \partial f_{z^*}\of{\bar x}$. It turns out that this is an adjoint process type relationship as pointed out in \cite{HamelSchrage13PJO}.

%%%New subsection
\subsection{A case study: Set-valued translative functions}
\label{SubSecTransFunc}

Let $X$, $Z$ be two topological linear spaces and $T \colon Z \to X$ an injective continuous linear operator. A function $f \colon X \to \mathcal P\of{Z, C}$ is called {\bf translative} with respect to $T$ (or just {\bf $T$-translative)} if
\[
\forall x \in X, \; \forall z \in Z \colon f\of{x + Tz} = f\of{x} + \cb{z}. 
\]

A special case of interest will be $Z = \R^m$, $\cb{h^1, \ldots, h^m} \subseteq X$ a set of $m$ linearly independent elements and $T \colon \R^m \to X$ defined by $Tz = \sum_{k =1}^m z_kh^k$. This construction is very close to (and motivated by) set-valued risk measures as shown below.

It is an easy exercise to show that a $T$-translative function $f$ can be represented as follows:
\begin{equation}
\label{EqFAF}
\forall x \in X \colon f\of{x} = \cb{z \in \R^m \mid x - Tz \in A_f}
\end{equation}
where $A_f = \cb{x \in X \mid 0 \in f\of{x}}$ is the zero sublevel set of $f$. This set satisfies 
\[
\forall z \in C \colon A_f - Tz \subseteq A_f
\]
since $f$ maps into $\mathcal P\of{Z, C}$. The latter property is called $\of{T, C}$-translativity of $A_f$. If $f$ maps into $\F(Z, C)$, i.e. it has closed values, then $A_f$ is also $\of{T, C}$-directionally closed, i.e.
\[
\of{\cb{z^k}_{k \in \N} \subseteq \R^m, \; \cb{z^k}_{k \in \N} \to 0, \; \forall k \in \N \colon x + Tz^k \in A_f} \quad \Rightarrow \quad x \in A_f.
\]

The representation \eqref{EqFAF} can be written as
\[
\forall x \in X \colon f\of{x} = \of{I_{A_f} \negthinspace\boxempty\negthinspace \alpha_T}\of{x} = \inf\cb{I_{A_f}\of{x_1} + \alpha_T\of{x_2} \mid x_1 + x_2 = x}
\]
where $\alpha_T \colon X \to \mathcal P\of{Z, C}$ is given by
\[
\alpha_T\of{x} = 
	\left\{
	\begin{array}{ccc}
	\cb{z} + C & : & x = Tz \\
	\emptyset & : & \text{otherwise}
	\end{array}
	\right.
\]
and $I_A$ is the set-valued indicator function of $A$: $I_A\of{x} = C$ if $x \in A$ and $I_A\of{x} = \emptyset$ if $x \not\in A$. Note that the function $\alpha_T$ is well-defined since $T$ is assumed to be injective.

We start the investigation of set-valued translative functions with their conjugates and make use of the fact that the conjugate of the infimal convolution of two functions is the sum of the two conjugates. For set-valued functions, this has been established in \cite[Lemma 2]{Hamel09SVVAN}. The conjugate of the indicator function is indeed the set-valued support function as shown in \cite{Hamel09SVVAN}:

\begin{align*}
I_{A_f}^*(x^*,z^*) = \bigcap_{x\in A_f}S_{(x^*,z^*)}(x).
\end{align*}
Moreover,
\begin{align*}
\alpha_T^*\of{x^*,z^*} & =\bigcap_{x\in X}\of{S_{(x^*,z^*)}(x)\idif\alpha_T(x)}=\bigcap_{u\in Z}\of{S_{(x^*,z^*)}(Tu)\idif\of{\cb{u}+C}}\\
	& = \bigcap_{u \in Z}\cb{z \in Z \mid z+u+C \subseteq S_{(x^*,z^*)}(Tu)}\\
	& = \cb{z \in Z \mid \forall u \in Z \colon z^*(z+u) \geq x^*(Tu)}\\
	& = \cb{z \in Z \mid z^*(z) \geq \sup_{u\in Z}(T^*x^*-z^*)(u)} 
	= 
	\left\{
	\begin{array}{ccc}
	H^+(z^*) & : & z^* = T^*x^* \\ 
	\emptyset & : &  z^* \neq T^*x^*
	\end{array}
	\right.
\end{align*}
Hence, for a $T$-translative function $f$ we get
\begin{equation}
\label{EqTransConjugate}
f^*(x^*,z^*)=I_{A_f}^*(x^*,z^*)+\alpha_T^*\of{x^*,z^*}=
	\left\{
	\begin{array}{ccc}
	\bigcap\limits_{x\in A_f}S_{(x^*,z^*)}(x) & : & z^* = T^*x^* \\ 
	\emptyset & : & z^* \neq T^*x^*
	\end{array}
	\right.
\end{equation}
and (see Remark \ref{RemNegConjugate})
\begin{align*}
(-f^*)(x^*,z^*)=H^+(z^*)\idif f^*(x^*,z^*)=
	\left\{
	\begin{array}{ccc}
	\cl\bigcup\limits_{x\in A_f}S_{(x^*,z^*)}(-x) & : & z^* = T^*x^* \\ 
	Z & : & z^* \neq T^*x^*
	\end{array}
	\right.
\end{align*}
since $H^+(z^*)\idif \emptyset = Z$ and $H^+(z^*)\idif \bigcap_{x \in A_f}S_{(x^*,z^*)}(x) = \cl\bigcup_{x \in A_f}\sqb{H^+(z^*)\idif S_{(x^*,z^*)}(x)} = S_{(x^*,z^*)}(-x)$ according to Proposition \ref{PropInfResCalc}.

If the function $f$ additionally maps into $\G(Z, C)$ and is proper, closed and convex, then the biconjugation theorem applies, and the following dual representation is obtained:
\begin{equation}
\label{EqTranslativeDual}
\forall x \in X \colon f\of{x} = \bigcap_{\substack{x^* \in X^* \\ T^*x^* \in \C}}\sqb{S_{(x^*,T^*x^*)}(x) \idif I_{A_f}^*(x^*, T^*x^*)}.
\end{equation}
If $f$ is additionally sublinear, then $A_f$ is a closed convex cone and \eqref{EqTranslativeDual} simplifies to 
\begin{equation}
\label{EqTranslativeDualSublin}
\forall x \in X \colon f\of{x} = \bigcap_{\substack{x^* \in A^-_f \\ T^*x^* \in \C}}S_{(x^*, T^*x^*)}(x)
\end{equation}
since in this case
\[
I_{A_f}^*(x^*,z^*) =
	\left\{
	\begin{array}{ccc}
	H^+(z^*) & : & x^* \in A^-_f \\
	\emptyset & : & \text{otherwise}
	\end{array}
	\right. .
\]
Of course, $A^-_f  = -\of{A_f}^+$.

The value of these formulas depends on how the dual data $x^*$, $T^*$ and $I_{A_f}^*$ can be interpreted in terms of the application at hand. We will show in Section \ref{SubSecRiskMeasures} below that this can be done very nicely.

\begin{example}
$Z = \R^m$, $T \colon \R^m \to X$ defined by $Tz = \sum_{k =1}^m z_kh^k$. Then
\[
\forall z \in \R^m \colon \of{T^*x^*}\of{z} = \sum_{k =1}^m x^*(h^k)z_k,
\]
thus $T^*x^*$ can be identified with $\of{x^*\of{h^1}, \ldots, x^*\of{h^m}}^T \in \R^m$.
\end{example}

We turn to the subdifferential of $T$-translative functions. The result reads as follows.

\begin{corollary}
\label{CorSubdifTrans}
Let $f \colon X \to \mathcal G\of{Z, C}$ be convex, $T$-translative and $z^* \in \C$. If $\partial f_{z^*}\of{\bar x} \neq \emptyset$ then
\begin{equation}
\label{EqSubdifTrans}
\partial f_{z^*}\of{\bar x} = \cb{x^* \in X^* \mid  z^* = T^*x^* \; \text{and} \; \forall x \in A_f \colon  S_{(x^*, T^*x^*)}(x) \supseteq S_{(x^*, T^*x^*)}(\bar x) \idif f(\bar x)}.
\end{equation}
\end{corollary}

{\sc Proof.} First, we show "$\subseteq$". The assumption $\partial f_{z^*}\of{\bar x} \neq \emptyset$ in conjunction with Proposition \ref{PropSubdiffInequality} implies $f(\bar x) \oplus H^+(z^*) \not\in \cb{Z, \emptyset}$.  Hence $S_{(x^*,z^*)}(\bar x) \idif f(\bar x) \not\in \cb{Z, \emptyset}$, and Proposition \ref{PropSubdiffConjugate} produces $f^*(x^*, z^*) \not\in \cb{Z, \emptyset}$. Take $x^* \in \partial f_{z^*}\of{\bar x}$. From \eqref{EqTransConjugate} we now obtain
\[
z^* = T^*x^* \quad \text{and} \quad f^*(x^*, z^*) = I_{A_f}^*(x^*,z^*).
\]
The definition of the set-valued support function yields that $x^*$ belongs to the right hand side of \eqref{EqSubdifTrans}.

Conversely, assume that $x^* \in X^*$ satisfies $z^* = T^*x^*$ as well as
\[
\forall x \in A_f \colon  S_{(x^*,z^*)}(x) \supseteq S_{(x^*,z^*)}(\bar x) \idif f(\bar x).
\]
Take $x \in \dom f$. Then
\[
\forall z \in f(x) \colon x - Tz \in A_f
\]
by $T$-translativity and hence
\[
\forall z \in f(x) \colon S_{(x^*,z^*)}(x - Tz) \supseteq S_{(x^*,z^*)}(\bar x) \idif f(\bar x) .
\]
Since $z^* = T^*x^*$ we have
\[
S_{(x^*,z^*)}(x - Tz) = S_{(x^*,z^*)}(x) +\cb{-z}
\]
and therefore
\[
\forall z \in f(x) \colon S_{(x^*,z^*)}(x) +\cb{-z} \supseteq S_{(x^*,z^*)}(\bar x) \idif f(\bar x).
\]
This means that any $\eta \in S_{(x^*,z^*)}(\bar x) \idif f(\bar x)$ satisfies
\[
\forall z \in f(x) \colon z + \eta \in S_{(x^*,z^*)}(x),
\]
thus $\eta \in S_{(x^*,z^*)}(x) \idif f(x)$. Hence
\[
\forall x \in \dom f \colon S_{(x^*,z^*)}(x) \idif f(x) \supseteq S_{(x^*,z^*)}(\bar x) \idif f(\bar x)
\]
which is, according to Proposition \ref{PropSubdiffConjugate}, equivalent to $x^* \in \partial f_{z^*}\of{\bar x}$.
\pend 

\medskip The above corollary tells us that the knowledge of $\partial f_{z^*}$ can be obtained by knowledge about $A_f$ and $T^*$. This becomes even more clear in the sublinear case.

\begin{corollary}
\label{CorSubdifTransSublin}
Let $f \colon X \to \mathcal G\of{Z, C}$ be sublinear, $T$-translative and $z^* \in \C$. If $\partial f_{z^*}\of{\bar x} \neq \emptyset$ then
\begin{equation}
\label{EqSubdifTrans}
\partial f_{z^*}\of{\bar x} = \cb{x^* \in X^* \mid z^* = T^*x^*, \; x^* \in A^-_f, \;  S_{(x^*,z^*)}(\bar x) = f(\bar x) \oplus H^+(z^*)}.
\end{equation}
\end{corollary}

{\sc Proof.} As observed above, in this case $A_f$ is a convex cone and $I^*$ can only attain the two values $H^+(z^*)$ for $x^* \in A^-_f$ and $\emptyset$ otherwise. Finally,
\[
S_{(x^*,z^*)}(\bar x) \idif f(\bar x) = H^+(z^*) \quad \Leftrightarrow \quad S_{(x^*,z^*)}(\bar x) = f(\bar x) \oplus H^+(z^*).
\]
The result now follows from Corollary \ref{CorSubdifTrans}.
\pend

%%%New subsection
\subsection{Comments on vector- and set-valued convex analysis}

The history of convex analysis for scalar functions is a continuing success story, and this area of mathematics is the theoretical basis for linear and nonlinear, in particular non-smooth, optimization and optimal control theory: compare \cite[p. 3]{Rockafellar93SIAMRev}\footnote{`In fact the great watershed in optimization isn't between linearity and nonlinearity, but convexity and nonconvexity.'} or the preface of \cite[p. xii]{BenTalNemirovski01Book}\footnote{`Theoretically, what modern optimization can solve well are {\em convex optimization problems}.'}.

Surprisingly, the gap between theory and applications (in optimization and multi-criteria decision making) is much wider for vector- or even set-valued functions. For example, there is no canonical (Fenchel) conjugate of a vector-valued function, but rather a whole bunch of different definitions which work under different assumptions (see below for references).

If one ignores for a moment scalarization approaches, then there are basically two different paths to a ``vector-valued" convex analysis.

The first one simply consists in an extended interpretation of the infimum and the supremum in formulas like the definition of the Fenchel conjugate: Under the assumption that the function maps into a conditional complete vector lattice (this means that every set which is bounded from below with respect to the vector order has an infimum in the space) one considers infima/suprema with respect to the vector order. This approach has been followed by J. Zowe \cite{Zowe75JMAA, Zowe78JMAA},  K.-H. Elster  and R. Nehse \cite{ElsterNehse75Opt}, \cite{Brumelle78MOR}, J. M. Borwein \cite{Borwein81MS}, C. Zalinescu \cite{Zalinescu83ASUC}, S. S. Kutateladze \cite{Kutateladze79RMS} and others. One may compare \cite{BorweinPenotThera84JMAA} for the state of the art in the mid 1980ies and more references. This approach has the advantage that a corresponding version of the Hahn-Banach theorem is available which is due to L. V. Kantorovich, see for example M. M. Day's book \cite{Day58Book}. Disadvantages are, of course, the strong assumptions to the image space and, even worth for applications, the fact that a vector infimum/supremum is hardly an appropriate concept when it comes to  vector optimization and multi-criteria decision making.

In the second approach, infima and suprema are therefore replaced by sets of minimal and maximal points, respectively, with respect to the vector order. This worked for (and was motivated by) applications of vector optimization, but made the task of developing a corresponding vector-valued convex analysis incredibly harder: It turns out that ``dual constructions" like conjugates or dual optimization problems become set-valued: `for a vector problem, its dual constructed by several means, is a problem whose objective function is set-valued, whatever the objective of the primal problem be' (\cite[p. 57]{Luc89Book}).  Set-valued Legendre--Fenchel conjugates with maximal points replacing the supremum appear in \cite{Luc89Book, NakayamaSawaragiTanino85Book, SawaragiTanino80JOTA}, with weakly maximal points in \cite{Malivert92LNEMS, NakayamaSawaragiTanino85Book}, with (weakly) supremal points in \cite{Kawasaki81MOR, Postolica86-1, Postolica86-2, Song97JMAA, Song98DM, Tanino92JMAA} and an even more general construction involving ``non-submitted" points is used in \cite{DoleckiMalivert93Opt}, for example.

A major difficulty for this approach is the lack of an appropriate Hahn-Banach theorem which is at the heart of convex analysis: One has to turn to scalarizations in order to apply the ``usual" Hahn-Banach argument. J. Zowe's paper \cite{Zowe74MathScan} shows how difficult it is to get back to vector-valued concepts after a scalarization.

In both approaches, continuous linear operators were used as dual variables. One way to avoid this again is a scalarization approach: An early attempt is J. Jahn's work \cite{Jahn83MP} (compare also \cite[Chapter 8]{Jahn04Book}). This approach leads to peculiar difficulties even if the problem at hand is linear: Compare \cite[Conclusions]{Jahn83MP} and the discussion at the ends of \cite[Section 8.2 and 8.3]{Jahn04Book}. A modern account is given in \cite{BotGradWanka09Book} which leads to dual problems with a, in general, non-convex feasibility set even if the original problem is convex (or linear). 

Let us mention that there are at least two quite different attempts to answer the duality question for vector problems: In \cite{Azimov88MathSb, Azimov08JOTA} as well as in \cite{BrecknerKolumban69MathScan} Fenchel conjugates of vector- or set-valued functions are defined in terms of scalar functions depending on an additional dual variable. Although in both attempts quite strong assumptions are imposed, they seem to be only a few steps short of the constructions in this section.

The approach summarized in \cite{Loehne11Book} is also based on scalarization via support functions, but it involves a set infimum/supremum which admits to obtain stronger results.

The concepts presented in this survey goes without the usual assumptions to the ordering cone $C$ (non-empty interior, pointedness, generating a lattice order etc.), and it basically produces set-valued versions of all the known (duality) formulas for scalar convex functions, and this includes the case of vector-valued functions. A crucial observation is the theoretical equivalence of a convex $\G(Z, C)$-valued function $f$ and the family $\cb{\vp_{f, z^*}}_{z^* \in C^+\bs\{0\}}$. Formula \eqref{EqTranslativeDual} is an example for how the set-valued theory tells us what kind of scalarizations should be taken into consideration. New insights can be obtained by investigating relationships between the two components of the dual variable $(x^*, z^*)$ which is essentially of adjoint process duality type (see \cite[Section 4]{HamelSchrage13PJO}). Set-valued functions satisfying (a) and (b) of Proposition \ref{PropConlinearFunction} are sometimes called linear, e.g. in \cite{Olko98AM}.

Directional derivatives for set-valued functions are usually defined at points of its graph, thus fixing (only) one element in the image set along with a point in the pre-image set. The standard reference is \cite{AubinFrankowska90Book}, and Mordukhovich's coderivative \cite{Mordukhovich97NA} is of the same type. Compare also \cite{Yang98MMOR}. Quite a different path is the attempt to embed certain subsets of $\P(Z)$ into a linear space and then use the usual ``linear" constructions, see \cite{Kuroiwa09JNCA} for an example. This, of course, only works under strong assumptions since, in general, $\G(Z, C)$ cannot be embedded into a linear space even if one drops $\emptyset$ and $Z$.

Concerning subgradients for set-valued functions, the paper \cite{HernandezRodriguezMarin11JOTA} presents an overview over the existing concepts each of which is afflicted with a peculiar difficulty: for example, the `weak subgradient' of \cite{ChenJahn98MMOR} leaves the realm of convexity, the `strong subgradient' introduced in \cite{HernandezRodriguezMarin11JOTA} needs an artificial exclusion condition in its definition and requires rather strong assumptions for its existence (see \cite[Definition 3.2 and Theorem 4.1]{HernandezRodriguezMarin11JOTA}). Other concepts define subgradients at points of the graph rather than at points of the domain, see \cite{Borwein81MS}, \cite{NakayamaSawaragiTanino85Book}, \cite[Definition 2.1]{BaierJahn99JOTA} and also the `positive subgradients' defined in \cite[Definition 2.5]{ChenLiWu09EJOR}, \cite[Definition 3.1]{HernandezLoehneRodriguezMarinTammer13Opt} and the `$k$-subgradients' of \cite[Definition 7.1.9]{BotGradWanka09Book} among many others.

Most of those concepts use linear operators as dual variables, but when it comes to existence very often operators of rank 1 show up, see, for example, \cite[Theorem 4.1]{HernandezRodriguezMarin11JOTA}. The (straightforward) relationships are discussed in \cite[p. 331]{BotGradWanka09Book} and \cite[Section 4]{HernandezLoehneRodriguezMarinTammer13Opt}.
 
We interpret this as evidence that, unless the image space is a (conditional) complete vector lattice and the Hahn-Banach-Kantorovitch theorem is available, the dual variables should involve linear functionals rather than linear operators. Using ``conlinear functions" generated by pairs of linear functionals, the constructions in Sections \ref{SubSecDirDer} and \ref{SubSecSubdifferential} offer a way to obtain results which are very close in shape to the scalar case and avoid strong assumptions to the ordering cone in $Z$. Moreover, in contrast to most of the ``vectorial" constructions in the literature (for example, see the discussion in \cite[p. 313]{BotGradWanka09Book}), our set-valued results reproduce the ones for scalar extended real-valued functions as special cases; this includes e.g. existence of subgradients and strong duality with attainment of the solution for the dual problem.

The subdifferential as given in Definition \ref{DefSubdiff} is exactly the same set which is called the `conjugate to' $f$ in \cite[Definition 2 and the remark thereafter]{Pshenichnyi72CSA} provided one assumes that every expression in B. N. Pshenichnyi's definition is finite. Section \ref{SubSecSubdifferential} should make it clear why we call it a subdifferential; the relationship to convex process duality can be found in \cite{HamelSchrage13PJO}. It should be pointed out that the complete lattice approach of this survey also adds new insights to scalar convex analysis: the improper case, in particular the function value $-\infty$ can be dealt with using the residuation. We refer to \cite{HamelSchrage12JCA}.

Scalar translative functions appear in many areas of applied mathematics, for example probability (quantile functions and lower previsions \cite{Walley91Book}), insurance and finance (constant additive insurance premiums \cite{PanjerWangYoung97IME} and cash additive risk measures, introduced in \cite{ArtDelEbeHea99MF} and reviewed in \cite{FoellmerSchied11Book}), mathematical economics (benefit and shortage functions \cite{Luenberger92JOTA}, \cite{Luenberger92JME}), vector optimization (nonlinear scalarization functions, compare \cite{GerthWeidner90JOTA} also for earlier references and \cite{GoeRiaTamZal03Book} for an overview) and idempotent analysis (compare the survey \cite{Kolokoltsov92}) as well as in max-plus algebra (see e.g. \cite{CohenGaubertQuadratSinger05}). A relationships between vector optimization and risk measures in finance is pointed out in \cite{Heyde06}.

Following an idea of \cite{JouiniMeddebTouzi04FS}, in \cite{HamelHeyde10SIFIN}, \cite{HamelHeydeRudloff11MAFE} cash additive risk measures have been generalized to set-valued risk measures for multivariate positions which turned out to be $T$-translative for some special $T$. Thus, such functions are important in applications, and they provide examples for the set optimization theory of this survey.

%%%
\section{Set-valued optimization}

%%%
\subsection{Unconstraint problems}

Within the set-up of the previous section, the basic problem again is
\[
\tag{P} \text{minimize} \quad f(x) \quad \text{subject to} \quad x \in X.
\]
The difficulty with the solution concept given in Definition \ref{DefSetSolution} is that solutions are, in general, sets rather than single points. Thus, optimality conditions such as ``zero belongs to the subdifferential of some function" should actually be taken ``at sets" rather than ``at points." Of course, this does not sound very attractive. The following construction provides a remedy.

\begin{definition}
\label{DefInfTranslation} Let $f \colon X \to \G\of{Z, C}$ be a function and $M \subseteq X$ a non-empty set. The function $\hat f\of{\cdot; M} \colon X \to \G\of{Z, C}$ defined by
\begin{equation}
\label{EqInfTranslation}
\hat f\of{x; M} = \inf_{u \in M}f\of{x+u} = \cl\co\bigcup_{u \in M}f\of{x+u}
\end{equation}
is called the inf-translation of $f$ by $M$.
\end{definition}

The function $\hat f\of{\cdot; M}$ coincides with the canonical extension of $f$ at
$M + \cb{x}$ as defined in \cite{HeydeLoehne11Opt}. A few elementary properties of the inf-translation are
collected in the following lemma.

\begin{lemma}\label{LemInfShiftProps}
Let $M \subseteq X$ be non-empty and $f \colon X \to \G\of{Z, C}$ a function.

(a) If $M \subseteq N \subseteq X$ then $\hat f\of{x; M} \subseteq \hat f\of{x; N}$ for
all $x \in X$.

(b) $\inf_{x \in X} f\of{x} = \inf_{x \in X} \hat f\of{x; M}$.

(c) If $f$  and $M$ are convex, so is $\hat f\of{\cdot; M} \colon X \to \G\of{Z, C}$, and in this case $\hat f\of{x; M} = \cl\bigcup_{u \in M}f\of{u+x}$.
\end{lemma}

{\sc Proof.} The proof can be found in \cite{HamelSchrage13PJO}. \pend

\begin{proposition}
\label{PropInfShiftChar} Let $f \colon X \to \G\of{Z, C}$ be a convex function and
$\emptyset \neq M \subseteq \dom f$. The following statements are equivalent:

(a) $M$ is an infimizer for $f$;

(b) $\cb{0} \subseteq X$ is an infimizer for $\hat f\of{\cdot; M}$;

(c) $\cb{0}$ is an infimizer for $\hat f\of{\cdot; \co M}$ and $\hat f\of{0; M} = \hat
f\of{0; \co M}$.
\end{proposition}

{\sc Proof.} The equivalence of (a) and (b) is immediate from $\hat f\of{0; M} = \inf_{u
\in M}f\of{u}$ and Lemma \ref{LemInfShiftProps}, (b). The equivalence of (a) and (c)
follows from $\hat f\of{0; \co M} = \inf_{u \in \co M}f\of{u}$ and Lemma \ref{LemInfShiftProps}
(b). \pend

The previous proposition makes clear that by an inf-translation an infimizer (set) can be reduced to a single point, namely just $0 \in X$.
Moreover, it should be apparent that we need to consider $\hat f\of{\cdot; \co M}$: Since
we want to characterize infimizers via directional derivatives and subdifferentials, a
convex function is needed, and $\hat f\of{\cdot; M}$ is not convex in general even if $f$
is convex (find a counterexample!). Obviously, an infimizer is not necessarily a convex set; on the contrary, sometimes
one prefers a nonconvex one, for example a collection of vertices of a polyhedral set
instead of whole faces.

\begin{theorem}
\label{ThmOptimality} Let $f \colon X \to \G(Z, C)$ be a convex function satisfying
\[
I\of{f} = \inf_{x \in X}f\of{x} \not\in \cb{Z, \emptyset}.
\]
Then $f$ is proper, and the set $\Gamma^+\of{f} = \cb{z^* \in C^+\bs\{0\} \mid I\of{f}\oplus H^+(z^*) \neq Z}$ is non-empty. Moreover, a set $M \subseteq X$ is an infimizer for $f$ if, and only if, $\hat f\of{0; M} = \hat f\of{0; \co M}$ and
\[
0 \in \bigcap_{z^* \in \Gamma^+\of{f}}\partial\hat f_{z^*}\of{\cdot; \co M}\of{0}.
\]
\end{theorem}

{\sc Proof.} Since $\cb{0}$ is a singleton infimizer of the function $x \mapsto \hat f\of{x; M}$, $\bar x = 0 \in X$ satisfies (c) of Proposition \ref{PropZeroInSubdiff} with $f$ replaced by $\hat f\of{\cdot; M}$ for each $z^* \in \Gamma^+\of{f}$. Now, the result follows from Proposition \ref{PropZeroInSubdiff} and Proposition \ref{PropInfShiftChar}.  \pend

\medskip Theorem \ref{ThmOptimality} highlights the use of the ``$z^*$-wise" defined directional derivatives and subdifferentials. One needs to take into consideration all reasonable (= proper) scalarizations at the same time in order to characterize infimizers.

%%%New subsection
\subsection{Constrained problems and Lagrange duality}

\label{SubSecLagrange}

Let $Y$ be another locally convex spaces with topological dual $Y^*$, and $D \subseteq Y$ a convex cone. The set $\G\of{Y, D}$ is defined in the same way as $\G\of{Z, C}$. Finally, let $f \colon X \to \G\of{Z, C}$ and $g \colon X \to \G(Y ,D)$ be two functions. We are interested in the problem
\[
\tag{PC} \mbox{minimize} \quad f(x) \quad \mbox{subject to} \quad  0 \in g\of{x}.
\]
The set
\[
\mathcal X = \cb{x \in X \mid 0 \in g\of{x}} 
\]
is called the feasible set for (PC) and $\mathcal I(f, g) = \inf\cb{f\of{x} \mid x \in \mathcal X}$ is the optimal value of the problem. With Definition \ref{DefMinSet} in view we define a solution of (PC) as follows.

\begin{definition}
\label{DefMinProblem}
A set $M \subseteq \mathcal X$ is called a solution of (PC) if

(a) $\inf\cb{f\of{x} \mid x \in M} = \mathcal I(f, g)$,

(b) $\bar x \in M$, $x \in \mathcal X$, $f\of{x} \supseteq f\of{\bar x}$ imply $f\of{x} = f\of{\bar x}$.
\end{definition}

Clearly, $M \subseteq X$ is a solution of  (PC) if, and only if $f\sqb{M}$ generates the infimum of $f
\sqb{\mathcal X} = \cb{f(x) \mid x \in \mathcal X}$ and each $f\of{\bar x}$ for $\bar x \in M$ is minimal in $f
\sqb{\mathcal X}$ with respect to $\supseteq$.

We define the Lagrangian $L \colon X \times Y^* \times C^+\bs\{0\} \to \G\of{Z, C}$ of
problem (PC) by
\begin{equation}
\label{EqDefLagrangian}
L\of{x, y^*, z^*} = f\of{x} \oplus \bigcup_{y \in g\of{x}}S_{\of{y^*, z^*}}\of{y}
    = f\of{x} \oplus \inf\cb{S_{\of{y^*, z^*}}\of{y} \mid y \in g\of{x}}.
\end{equation}

Under a mild condition, the primal problem can be reconstructed from the Lagrangian.

\begin{proposition}
\label{PropReconPrimal} If $f\of{x} \neq Z$ for each $x \in \mathcal X$, then
\[
\sup_{\of{y^*, z^*} \in Y^* \times C^+\bs\{0\}} \negthickspace\negthickspace L\of{x, y^*, z^*}
    = \negthickspace\negthickspace \bigcap_{\of{y^*, z^*} \in D^+ \times C^+\bs\{0\}} \negthickspace\negthickspace L\of{x, y^*, z^*}
    = \left\{
    \begin{array}{ccc}
    f\of{x} & : & 0 \in g\of{x} \\
    \emptyset & : & 0 \not\in g\of{x}.
    \end{array}
    \right. 
\]
\end{proposition}

{\sc Proof.} The proof is based on the assumption that the values of $f$ and $g$ are closed convex sets. See \cite{HamelLoehne13JOTA-OF} for details. \pend

\medskip The next proposition is a Lagrange sufficient condition which is a simple, but important result with an algorithmic character since it admits to test if a given set is an infimizer of (PC).

\begin{proposition}
\label{PropLagrangeSuff}
Let $M \subseteq \mathcal X$ be a non-empty set of feasible points for (PC). Assume that for each $z^* \in \C$ there is $y^* \in D^+$ satisfying
\begin{equation}
\label{EqLagrangeInfimizer}
\hat f\of{0; M} \oplus  \inf_{y \in \hat g\of{0; M}}S_{\of{y^*, z^*}}\of{y} = \inf_{x \in X} L\of{x, y^*, z^*}
\end{equation}
and
\begin{equation}
\label{EqCompSlack}
\inf_{y \in \hat g\of{0; M}} S_{\of{y^*, z^*}}\of{y} = H^+(z^*).
\end{equation}

Then, $M$ is an infimizer for (PC).
\end{proposition}

{\sc Proof.} Using \eqref{EqCompSlack} and \eqref{EqLagrangeInfimizer} we obtain
\begin{align*}
\hat f\of{0; M} \oplus H^+(z^*) & = \hat f\of{0; M} \oplus \inf_{y \in \hat g\of{0; M}} S_{\of{y^*, z^*}}\of{y} 
	 = \inf_{x' \in X} L\of{x', y^*, z^*} \\
	 & \supseteq f(x) \oplus \inf_{y \in g\of{x}}S_{\of{y^*, z^*}}\of{y}  \supseteq f(x) \oplus H^+(z^*)
\end{align*}
for all $x \in \mathcal X$ since $S_{\of{y^*, z^*}}\of{0} = H^+(z^*)$. Taking the infimum over the feasible $x$ on the right hand side and then the intersection over $z^* \in \C$ on both sides while observing $\hat f\of{0; M} = \inf_{u \in M}f\of{u}$ we obtain that $M$ indeed is an infimizer for (PC). 
\pend

\medskip Condition \eqref{EqCompSlack} serves as set-valued complementary slackness condition. If one considers the Lagrange function $(x, y^*, z^*) \mapsto  \hat L\of{x, y^*, z^*; M}$ for the ``inf-translated" problem
\[
\text{minimize} \quad \hat f(x; M) \quad \text{subject to} \quad 0 \in \hat g(x; M)
\]
then condition \eqref{EqLagrangeInfimizer} means that the infimum of the Lagrange function for the original problem coincides with $\hat L\of{0, y^*, z^*; M}$. Finally, if $z^* \in \C$ and $y^* \in D^+$ satisfy \eqref{EqCompSlack} and \eqref{EqLagrangeInfimizer} then $y^*$ is nothing else than a Lagrange multiplier for the by $z^*$ scalarized problem. One may therefore expect that strong duality is something like ``strong duality for all reasonable scalarized problems." This idea works as shown in the following.

Define the function $h \colon Y^* \times C^+\bs\{0\} \to \G\of{Z, C}$ by
\[
h\of{y^*, z^*} = \inf_{x \in X}L\of{x, y^*, z^*} = \cl\bigcup_{x \in X} L\of{x, y^*, z^*}.
\]
Since the values of $L(\cdot, y^*, z^*)$ are closed half spaces with the same normal $z^*$, the convex hull can be dropped in the
infimum. The dual problem,
\[
{\f \tag{DC}} \mbox{maximize} \quad h\of{y^*, z^*} \quad \mbox{subject to} \quad  y^* \in Y^*, \;z^* \in
 C^+\bs\{0\},
\]
thus consists in finding
\[
d = \sup_{y^* \in Y^*, \, z^* \in C^+\bs\{0\}} h\of{y^*, z^*} =
    \bigcap _{y^* \in Y^*, \, z^* \in C^+\bs\{0\}} h\of{y^*, z^*}
\]
and corresponding (full) solutions. The following weak duality result is immediate.

\begin{proposition}
\label{PropWeakDuality} Let $f\colon X \to \F\of{Z, C}$
and $g\colon X \to \F\of{Y, D}$. Then 
\[
\sup\cb{h\of{y^*, z^*} \mid y^* \in Y^*, \; z^* \in C^+\bs\{0\}} \supseteq \inf\cb{f\of{x} \mid x \in X, \; 0 \in g\of{x}}.
\]
\end{proposition}

{\sc Proof.} This is true since for $\of{y^*, z^*} \in Y^* \times C^+\bs\{0\}$ and $x \in
X$ satisfying $0 \in g\of{x}$ we have
\[
h\of{y^*, z^*}  \supseteq f\of{x} \oplus
        \cl\bigcup_{y \in g\of{x}}S_{\of{y^*, z^*}}\of{y} 
     \supseteq f\of{x} \oplus S_{\of{y^*, z^*}}\of{0} = f\of{x} \oplus H^+\of{z^*}.
\]
\pend

As usual, a constraint qualification condition is needed as part of  sufficient conditions for
strong duality. The following condition is called the {\em Slater condition} for problem (PC):
\begin{equation}
\label{EqSlater}
\exists \bar{x} \in \dom f \colon g\of{\bar{x}} \cap \Int\of{-D} \neq \emptyset.
\end{equation}
The implicit assumption is $\Int D \neq \emptyset$. 

\begin{theorem}\label{ThmStrongDuality}
Assume $p=\inf\cb{f(x) \mid x \in \mathcal X} \neq Z$. If $f \colon X \to \mathcal \G\of{Z, C}$ and $g \colon X \to \mathcal \G\of{Y, D}$ are convex and the Slater condition for problem (PC) is satisfied then strong duality holds for (PC), that is 
\begin{align}\label{EqSDSet1}
& \inf\cb{f\of{x} \mid 0 \in g\of{x}}
    = \sup \cb{h\of{y^*, z^*} \mid y^* \in Y^*, \; z^* \in C^+\bs\{0\}}\\[.2cm]
\label{EqSDSet2}
& z^* \in C^+\bs\{0\},\; p \oplus H^+\of{z^*}\neq Z \quad \Rightarrow \quad \exists y^* \in Y^* \colon p \oplus H^+\of{z^*} = h\of{y^*, z^*}
\end{align}
\end{theorem}

{\sc Proof.} \cite{HamelLoehne13JOTA-OF}. \pend

\medskip Note that the assumption $p \neq Z$ implies the existence of $z^* \in \C$ with $p \oplus H^+\of{z^*}\neq Z$.  Thus, \eqref{EqSDSet2} is attainment of the supremum for the dual problem ``$z^*$-wise."

\begin{corollary}
\label{CorExistenceDualSol} Under the assumptions of the strong duality theorem, the set
\[
\Delta = \cb{\of{y^*, z^*} \in Y^* \times C^+\bs\{0\} \mid Z \neq p \oplus H^+(z^*) = h\of{y^*, z^*}}    
\]
is non-empty and a full solution of the dual problem (DC).
\end{corollary}

{\sc Proof.} See \cite{HamelLoehne13JOTA-OF}. \pend

%%%New subsection
\subsection{Comments on set optimization duality}

\label{SecCommDuality}

Among the first papers in which optimization problems with a set-valued constraint have been systematically studied are \cite{Borwein77SICON}, \cite{Borwein81MS} and \cite{Oettli82}. It is, for example, instructive to realize that the Lagrange function in \eqref{EqDefLagrangian} is nothing else, but a set-valued version of the one in \cite[p. 197]{Oettli82}. Compare also \cite[Theorem 3.28]{Loehne11Book}.

Whereas in \cite[Problem (P) in (3.1)]{Borwein81MS} the vector infimum serves as the building block for optimality, in \cite[Theorem 3]{Borwein77SICON} a Lagrange duality result is established for properly efficient points of vector optimization problems. The dual variables are rank one linear operators. Similarly, in \cite[Theorem 3.3]{Song96JMAA} and also \cite[Theorem 3.3]{Song97JOTA}, rank one linear operators and a set-valued Lagrange function (see equation \eqref{EqLagrangeOperator} below) are used under strong assumptions (cones with weakly compact base). A similar idea can be found in the proof the the Lagrangian duality theorem, \cite[Theorem 1.6 on p. 113]{Luc89Book} under the assumption that the ordering cone in $Z$ has non-empty interior. These examples may suffice with respect to vector optimization problems in view although the literature is huge.

In \cite{Kuroiwa98-3}, \cite{Kuroiwa99RIMS} the same type of set-valued Lagrangian has been used (without giving proofs) in connection with set relations, i.e., basically the solution concept IIa of Section \ref{SubsecMinConcepts}. The more recent \cite{HernandezRodriguezMarin07JOTA} and \cite{HernandezRodriguezMarin07PJO} proceed similarly: Theorem 3.3 in \cite{HernandezRodriguezMarin07JOTA} (basically the same as Theorem 4.2 in \cite{HernandezRodriguezMarin07PJO}) is a Lagrange duality result for weakly $\lel_C$-minimal solutions with Lagrange function
\begin{equation}
\label{EqLagrangeOperator}
f(x) + (T \circ g)(x) = f(x) + \cb{Ty \mid y \in g(x)}
\end{equation}
where $T \in \mathcal L(Y, Z)$, the set of continuous linear operators from $Y$ to $Z$. It is again based on rank one operators, an idea which at least dates back to \cite[Theorem 4.1]{Corley87JOTA}. The same set of dual variables is used in \cite{HamelEtAl04JCA} for a Lagrangian approach to linear vector optimization. However, the Lagrange function, even for a vector-valued problem, already is a set-valued one.

A thorough discussion of optimality conditions of Fermat and Lagrange type for (non-convex) set-valued optimization problems based on the minimality concept III can be found in \cite{DureaStrugariu11Opt} (compare also the references therein). These conditions are formulated in terms of the Mordukhovich subdifferential. It might be worth noting that the use of $\F(Z, C)$-valued functions `gives better conclusions' \cite[Remark 3.10]{DureaStrugariu11Opt}.

A complete lattice approach based on infimal and supremal sets was developed in \cite{Loehne11Book} and \cite{HernandezLoehneRodriguezMarinTammer13Opt}. The Lagrange function for a vector-valued function $f$ and a set-valued $G$ has the form
\[
f(x) + \Inf\cb{y^*(y)c \mid y \in G(x)}
\]
where $\Inf$ stands for the infimal set and $c \in \Int C$ is a (fixed) element. Assumptions, of course, include $\Int C \neq \emptyset$. The same assumption also is crucial in \cite{LiLiuSunTeoYao11NFAO}; Theorem 3.2 and 3.3 therein is probably as far as one get in terms of conjugate duality based on ``suprema" of a set, i.e. the elements which belong to the closure of the set, but are not dominated with respect to the relation which is generated by the interior of the ordering cone.

Other approaches rely on other set-valued derivatives, for example on contingent epiderivatives \cite{GoetzJahn00SIOpt} or coderivatives \cite{Truong05JMAA}, \cite{MordukhovichTruong07AM}.

In virtually all approaches for set/vector optimization problems known to the authors, the strong duality assertion is based on the assumption of the existence of a (weakly, properly etc.) minimal element of the primal problem either with respect to the vector order (see \cite{Corley87JOTA}, \cite[Theorem 1.6 on p. 113, Theorem 2.7 on p. 119]{Luc89Book}, \cite[Theorem 3.4, 3.5]{Song98DM}, \cite[Theorem 5.2.4, 5.2.6]{BotGradWanka09Book}) or with respect to a set relation (see \cite[Theorem 3.3]{HernandezRodriguezMarin07JOTA},  \cite[Theorem 4.2]{HernandezRodriguezMarin07PJO}). The two exceptions are the approaches in \cite{Loehne11Book} and \cite{HamelLoehne13JOTA-OF} where the primal problems only have finite values in some sense and still existence for the dual problems is obtained--which is standard in the scalar case. In \cite[p. 98]{Loehne11Book} (see also Open Problem 3.6 therein with respect to Fenchel duality) and \cite{HernandezLoehneRodriguezMarinTammer13Opt} it is discussed that the approach based on infimal/supremal sets indeed yields strong duality, but it is not clear wether the existence of the dual solution can be guaranteed without the $z^*$-component of the dual variable.

By means of the ``complete lattice approach" surveyed here, the type of results which is known from the scalar case can be transferred to a ``set level." Strong duality then indeed means ``inf equals sup" and includes the existence of dual solutions: compare \cite{Loehne11Book}, \cite{HamelLoehne13JOTA-OF} for Lagrange duality and \cite{Hamel11Opt} for Fenchel-Rockafellar duality.

The reduction of a ``set solution" in the sense of Definition \ref{DefMinProblem} to a ``point solution" via an inf-translation (see Definition \ref{DefInfTranslation}) is due to \cite{HamelSchrage13PJO}. The exploitation of this construction seems to be very promising for obtaining optimality conditions and algorithms.

The complementary slackness condition given in Proposition \ref{PropLagrangeSuff} seems to be new although it clearly is in the spirit of \cite[formulae (10), (12)]{Borwein77MP}.

%%%New section
\section{Applications}

%%%
\subsection{Vector optimization}

In this section, let $X$ and $Z$ be as in Section \ref{SecConvAnal} and $C \subseteq Z$ a closed, convex, pointed (i.e. $C \cap -C = \{0\}$) and non-trivial cone. Then, $\leq_C$ is a partial order (i.e. also antisymmetric). Moreover, let a function $F \colon X \to Z\cup\cb{-\infty,+\infty}$ be given. Defining a function $f \colon X \to \G(Z, C)$ by
\[
f(x) = 
	\left\{
	        \begin{array}{ccc}
				F(x) + C & : &  F(x) \in Z \\
				   Z  & : & F(x) = -\infty \\
				\emptyset & : & F(x) = +\infty
		    \end{array}
		  \right.
\]
we observe
\[
f(x_1) \supseteq f(x_2) \quad \Leftrightarrow \quad F(x_1) \leq_C F(x_2),
\]
where it is understood that $-\infty \leq_C z \leq_C +\infty$ for all $z \in Z\cup\cb{-\infty,+\infty}$. Hence the two problems
\begin{align*}
\tag{VOP} \text{find minimizers w.r.t.} \; \leq_C \; \text{of} \; F(x) \; \text{subject to} \; 0 \in g(x), \\
\tag{SOP} \text{find minimizers w.r.t.} \; \supseteq \; \text{of} \; f(x) \; \text{subject to} \; 0 \in g(x)
\end{align*}
have the same feasible elements and the same minimizers. The minimizers of (VOP) are called 'minimal solutions' \cite[Definition 7.1]{Jahn04Book} or  `efficient solutions' \cite[Definition 2.5.1]{BotGradWanka09Book}. In most cases, it does not make sense to look  for the infimum in (VOP) with respect to $\leq_C$: It may not exist (not even for simple polyhedral cones $C$, see e.g. \cite[Example 1.9]{Loehne11Book}), and even if it does, it is not useful in practice at it refers to so-called utopia points which are typically not realizable by feasible points (i.e. ``decisions").

The (PC) version of (SOP) considered as an $\F(Z, C)$- or $\G(Z, C)$-valued problem is called the lattice extension of (VOP), and a solution of (VOP) is defined to be a solution of its lattice extension (see \cite{HeydeLoehne11Opt}, compare Definition \ref{DefMinProblem}). In this way, the notion of an ``infimum" makes a strong comeback, and the infimum attainment becomes a new feature in vector optimization, which is useful for theory and applications: It ensures that the decision maker possesses a sufficient amount of information about the problem if (s)he knows a solution. For a detailed discussion see \cite[Chapter 2]{Loehne11Book}. Note that one possibly obtains different solutions depending on the choice of $\F(Z, C)$ or $\G(Z, C)$ as image space. Since the infimum in $\G(Z, C)$ involves the convex hull, solutions of $\G(Z, C)$-valued problems may include ``fewer" elements, and this is in particular preferable for convex problems.

If $f$ is the ``lattice extension" of a vector-valued function $F$ as given above, the Lagrange function for (PC) takes the form
\begin{align*}
L\of{x, y^*, z^*} & = f(x) \oplus \inf_{y \in g(x)}S_{\of{y^*, z^*}}\of{y} = F(x) + \inf_{y \in g(x)}S_{\of{y^*, z^*}}\of{y} \\
	& = \inf_{y \in g(x)}\cb{z + F(x) \in Z \mid y^*(y) \leq z^*(z)} \\
	& = \cb{z \in Z \mid \inf_{y \in g(x)}y^*(y) + z^*\of{F(x)} \leq z^*(z)}
\end{align*}
whenever $F(x) \in Z$, $L\of{x, y^*, z^*} = \emptyset$ whenever $F(x) = +\infty$ or $g(x) = \emptyset$, and $L\of{x, y^*, z^*} = Z$ whenever $F(x) = -\infty$ and $g(x) \neq \emptyset$. The function $\Lambda_{z^*}\of{x, y^*} := z^*\of{F(x)} + \inf_{y \in g(x)}y^*(y)$ (with the convention $z^*(\pm\infty) = \pm\infty$) is the (classical) Lagrange function of the (scalar) problem
\[
\inf\cb{z^*\of{F(x)} \mid 0 \in g(x)}
\]
(see, for example, already \cite[p. 197]{Oettli82}). Moreover, if $g$ is generated by a vector-valued function $G \colon X \to Y\cup\cb{-\infty, +\infty}$ in the same way as $f$ by $F$, then 
\[
\inf_{y \in g(x)}y^*(y) =
	\left\{
	        \begin{array}{ccc}
				y^*\of{G(x)} & : &  G(y) \in Z, \; y^* \in D^+ \\
				   -\infty  & : & G(y) = -\infty, \; \text{or} \; G(y) \in Z \; \text{and} \; y^* \not\in D^+  \\
				+\infty & : & G(y) = +\infty.
		    \end{array}
		  \right.
\]
Thus, $\Lambda_{z^*}\of{x, y^*} = z^*\of{F(x)} + y^*\of{G(y)}$ whenever $F(x) \in Z$, $G(x) \in Y$ and $y^* \in D^+$. The dual objective becomes
\[
h(y^*, z^*) = \inf_{x \in X} L\of{x, y^*, z^*} = \cb{z \in Z \mid \inf_{x \in X}\Lambda_{z^*}\of{x, y^*} \leq z^*(z)}.
\]

\begin{corollary}
\label{CorStrongDualityVector}
Let $F$ be $C$-convex, $f$ its lattice extension and $g \colon X \to \G(Y, D)$ convex such that the Slater condition \eqref{EqSlater} is satisfied. If $I(f, g) = \inf\cb{f(x) \mid 0 \in g(x)} \not\in \cb{Z, \emptyset}$, then $\Gamma^+(f,g) = \cb{z^* \in \C \mid I(f, g) \oplus H^+(z^*) \neq Z}$ is non-empty and
\begin{align}\label{EqSDVec1}
I(f, g) & = \cl\bigcup\cb{F\of{x} \mid 0 \in g\of{x}}
    = \bigcap_{y^* \in D^*, \, z^* \in \Gamma^+(f,g)} \cb{z \in Z \mid \Lambda_{z^*}\of{x, y^*} \leq z^*(z)},\\[.2cm]
\label{EqSDVec2}
& \forall z^* \in \Gamma^+(f,g) \; \exists y^* \in Y^* \colon I(f, g) \oplus H^+\of{z^*} = h\of{y^*, z^*}.
\end{align}
\end{corollary}
{\sc Proof.} Of course, $f$ is convex if, and only if, $F$ is $C$-convex (see \cite[Definition 1.6 on p. 29]{Luc89Book} for a definition). Theorem \ref{ThmStrongDuality} and the above discussion produce the result. \pend

\medskip It might be worth it to compare Corollary \ref{CorStrongDualityVector} with standard duality results in vector optimization. First, there is no assumption about the existence of (weakly, properly) minimal solutions: This is in contrast to most results in vector optimization such as \cite[Theorem 3.7.4, 3.7.7]{GoeRiaTamZal03Book}, \cite[Theorem 8.7]{Jahn04Book}, \cite[Theorem 4.1.2, 4.1.4]{BotGradWanka09Book}. Secondly, there are no interior point assumptions to the cone $C$. Thirdly, with Corollary \ref{CorExistenceDualSol} in view, the existence of a dual solution in a set-valued sense is provided in the sense of the ``maximization" version of Definition \ref{DefSetSolution}. Finally, classical duality results in vector optimization can be obtained from Corollary \ref{CorStrongDualityVector} as it is described in \cite[Section 3.5]{Loehne11Book}.

%%%New subsection
\subsection{A case study: Linear vector optimization} 

\label{SubSecLVO}

We proceed with an exemplary application of the set-valued theory to linear vector optimization problems and show that we obtain what we expect in view of scalar linear programming duality: a dual program of the same type. In this section, we will write $\leq$ and $\geq$ for $\leq_{\R^m_+}$ and $\geq_{\R^m_+}$, respectively, for any $m \in \cb{1,2, \ldots}$.

Consider the linear vector optimization problem
\begin{equation} \tag{$P_L$}\label{PL}
{\textstyle\min_C} \quad P x \quad \text{subject to} \quad Ax \geq b,
\end{equation}
where $P \in \R^{q \times n}$, $A\in \R^{m \times n}$, $b \in \R^m$, and the cone $C$ is polyhedral convex with nonempty interior. A representation $C = \cb{z \in \R^q \mid W^T z \geq 0}$ by a matrix $W \in \R^{q \times k}$ is given. The feasible set is denoted by $S := \cb{x \in \R^n \mid Ax \geq b}$.

With (PC) in view we define $f(x) = Px + C$ and $g(x) = b - Ax + \R^m_+$. Then, the set $\cb{(x,z) \in \R^n \times \R^q \mid z \in f(x),\; 0 \in g(x)}$ is polyhedral convex. We modify the solution concept in Definition \ref{DefMinProblem} by adding the requirement that a solution is a finite set of vectors and directions, see \cite{Loehne11Book} and also \cite{HamelLoehneRudloff13JOGO}, the latter also including ``$\eps$-variants." The reason is that every polyhedral set can be expressed as the generalized convex hull of finitely many vectors and directions. Such a solution is called {\em finitely generated solution}, but we call it just {\em solution} if the context of polyhedral convex set-valued problems or the subclass of linear vector optimization problems is clear. To keep the notation simple, we only consider {\em bounded} problems here, that is, we assume
\begin{equation}
	\exists \bar z  \in \R^q \colon \forall x \in S \colon \bar z \leq_C Px. 
\end{equation}
Under this assumption, a solution consists of finitely many vectors only. For the general case, see \cite[Chapter 4]{Loehne11Book}. A solution to ($P_L$) is a nonempty finite set $\bar S \subseteq S$ of minimizers ('efficient solutions' in the most textbooks) such that $P[S] \subseteq P[\bar S]+C$, where the latter condition refers to infimum attainment in $\bar S$ with respect to the lattice extension (compare Definition \ref{DefSetSolution}).

Considering the lattice extension of ($P_L$) we show that the Lagrange technique from Section \ref{SubSecLagrange} leads to a dual problem, which enjoys nice properties and is useful for applications and algorithms. Re-labeling the dual variables by $u = y^*$, $w = z^*$ we obtain the Lagrangian
\begin{align*}
	L(x,u,w) &= Px + C + \cl\bigcup_{z \geq b - Ax} \cb{z \in \R^q \mid u^Ty \leq w^T z} \\
	         &=	Px + C + \cl\bigcup_{r \in \R^m_+} \cb{z \in \R^q \mid u^T (r-Ax+b) \leq w^T z}.
\end{align*}
The dual objective is
\begin{align*}
	h(u,w) &= \cl\bigcup_{x \in \R^n} L(x,u,w) \\
	       &= \cl\bigcup_{r \in \R^m_+,\, x \in \R^n,\, v \in C}  \cb{z \in \R^q \mid u^T (r-Ax+b) \leq w^T (z-Px-v)} \\
	       &= \cl\bigcup_{r \in \R^m_+,\, x \in \R^n,\, v \in C}  \cb{z \in \R^q \mid (w^T P-u^T A) x \leq w^T (z - v) -u^T (b+r)}	 \\
	       &= \left\{
	             \begin{array}{ccc} 
		           \cb{z \in \R^q \mid 0 \leq w^T z - u^T b} & : & A^T u = P^T w,\, u \geq 0,\, w \in C^+\bs\cb{0} \\
				   \R^q                                             & : & \text{otherwise.}
			     \end{array}
			   \right.
\end{align*}
Let $C^+=\cb{w \in \R^q \mid V^T w \geq 0}$ be a representation of $C^+$ by a matrix $V \in \R^{q \times l}$.
Note that a basis of $C^+$ is already sufficient to cover all values of the dual objective $h$ (see the end of Section \ref{SubSecConjugates}). If we fix some $c \in \Int C$, we obtain the (set-valued) dual problem
\begin{equation}\tag{$D_L$}\label{DL}
  \text{ maximize } D(u,w) \quad \text{subject to} \quad (u,w) \in T
\end{equation}
with objective function 
\begin{equation*}
	D \colon \R^m \times C^+ \to \G(\R^q,C),\quad D(u,w):= \cb{ z \in \R^q \mid u^T b \leq w^T z }
\end{equation*}
and feasible set
\begin{equation*}
	T:=\cb{(u,w) \in \R^m \times \R^q \mid A^T u = P^T w,\; u \geq 0,\; V^T w \geq 0,\; c^T w = 1}.
\end{equation*}
This dual problem has a very simple structure: linear constraints, a halfspace-valued objective function and maximization means to take the intersection over these halfspaces. The objective function is conlinear in $b$ and in $u$, i.e., $D(u,w)=S_{(u,w)}(b)=S_{(b,w)}(u)$, and therefore a natural replacement of the dual objective ``$b^T u$'' in (scalar) linear programming. A (finitely generated) solution of \eqref{DL} is a nonempty set $\bar T \subseteq T$ of maximizers with respect to the ordering $\supseteq$ satisfying $\bigcap_{(u,w)\in \bar T} D(u,w) = \bigcap_{(u,w) \in T} D(u,w)$, where the latter conditions means supremum attainment in $\bar T$.

\begin{remark} Using the construction of Example \ref{ExSelfInfimalSets}, we obtain an equivalent problem with a hyperplane-valued objective. This shows that we indeed have a very natural generalization of scalar linear programs to the vectorial case because in $\R$, a real number and a hyperplane are the same object. In more general linear spaces, vectors and half-spaces are dual in some sense. Compare the footnote on p. 2.
\end{remark}

Weak duality (see Proposition \ref{PropWeakDuality}) means that $x \in S$ and $(u,w) \in T$ imply
$D(u,w) \supseteq Px + C$. As a consequence, for every subset $\tilde T \subseteq T$ of feasible points, the set 
$\bigcap_{(u,w) \in \tilde T} D(u,w)$ is a superset (``outer approximation") of the set $\P:=\cb{Px \mid Ax \geq b}+C$, which is just the optimal value (the infimum) of the lattice extension. Likewise, for every subset $\tilde S \subseteq S$ of feasible points of \eqref{PL}, the set $\cl\co\bigcup_{x \in \tilde S} Px+C$ is a subset (``inner approximation") of $\P$.

Strong duality means that $\bigcap_{(u,w) \in T} D(u,w) = \P$. A constraint qualification is not needed as in the case of linear constraints in (scalar) convex programming. Note further that, if $\emptyset \neq \bar S \subseteq S$ such that $P[\bar S]$ is the set of vertices of $\P$, then $\bar S$ is a solution to \eqref{PL}. Likewise, a set $\emptyset \neq \bar T \subseteq T$ such that $\cb{D(u,w) \mid (u,w) \in \bar T}$ is the family of half-spaces supporting $\P$ in facets, then $\bar T$ is a solution of \eqref{DL}.

\begin{remark} In the vector optimization literature one can observe the longstanding paradigm that the dual of a vector optimization problem should be a vector optimization problem with the same ordering cone. To fulfill this requirement, problems of the type 
\begin{equation}
	{\textstyle\max_C} \quad z \quad \text{subject to} \quad z \in D(u,w),\; (u,w)\in T
\end{equation} 
have been considered, see e.g. \cite[Section 4.5.1]{BotGradWanka09Book} and in the linear case \cite{BotGradWanka12OL}. The price is high. In general, important properties like linearity of the constraints and convexity of the feasible set get lost by such a transformation.
\end{remark}

To emphasize the ``linear" character of problem \eqref{DL}, we transform it into an equivalent linear vector optimization problem:

\begin{equation}\tag{D$^*_L$}\label{DLstar}
	\textstyle{\max_K} \quad D^*(u,w) \quad \text{subject to} \quad (u,w) \in T,
\end{equation}
where the objective function $D^* \colon \R^q\times\R^m \to\R^q$, given by
\begin{equation*}
	  D^*(u,w) := (w_1,\dots,w_{q-1},b^T u)^T,
\end{equation*}
is linear and vector-valued, and the ordering cone is $K := \cb{z \in \R^q \mid z_1 = \dots = z_{q-1} = 0, z_q \geq 0}$. A (finitely generated) solution of \eqref{DLstar} is a nonempty set $\bar T \subseteq T$ of maximizers with respect to $\le_K$ in $\R^q$ induced by $K$ satisfying $D^*[T] \subseteq \co D^*[\bar T] - K$, where the latter condition refers to supremum attainment in $\bar T$ (with respect to the lattice extension with image space $\G(\R^q, K)$). 

\begin{proposition}
The problems \eqref{DL} and \eqref{DLstar} have the same solutions.
\end{proposition}

{\sc Proof.} See \cite[Theorems 4.38 and 4.57]{Loehne11Book}. \pend

\medskip In the sense of the previous proposition, \eqref{DL} and \eqref{DLstar} are equivalent. This means that the set-valued dual problem \eqref{DL} can be expressed as a linear vector optimization problem, however, with a different ordering cone $K$ and an interpretation of the duality relation which differs from the one in standard references.

Of course, we can derive a set-valued dual problem to \eqref{DLstar} by an analogous procedure. This leads to outer and inner approximations and different representations of $\D:=\cb{D^*(u,w) -K \mid (u,w)\in T}$, i.e., the optimal value of the lattice extension of \eqref{DLstar}.

Problem \eqref{DLstar} is called the {\em geometric dual problem}, and there is a further duality relation called {\em geometric duality} \cite{HeydeLoehne08SIAMOpt} between \eqref{PL} and \eqref{DLstar}: There is an inclusion-reversing one-to-one map between the proper faces of $\P$ and the proper $K$-maximal faces of $\D$. This means, for instance, that a vertex of one set can be used to describe a facet of the other set and vice versa. For a detailed explanation of geometric duality see \cite{HeydeLoehne08SIAMOpt, Loehne11Book}. Geometric duality has been extended to convex vector optimization problems, see \cite{Heyde13JCA}. The paper \cite{Luc11EJOR} is in the same spirit.

%%%New subsection
\subsection{Approximate solutions and algorithms}

In this section, we assume that $C$ is a closed convex cone. Let $f \colon X \to \G(Z, C)$ be a function. The starting point for constructing algorithms for solving the problem (P) (see Section \ref{SubsecMinConcepts}), i.e.
\[
\tag{P} \text{minimize} \quad f(x) \quad \text{subject to} \quad x \in X
\]
should be Definition \ref{DefSetSolution}: It involves minimal values of $f$ as well as the infimum taken in $\G(Z,C)$. In order to make algorithms reasonable, both notions should be replaced by appropriate approximate versions. 

Recall $I(f) = \inf_{x \in X} f(x)$. Two sets $A, B \in \G(Z, C)$ are called an outer approximation and an inner approximation of $I(f)$, respectively, if $A \supseteq I(f) \supseteq B$.  Outer and inner approximations of $I(f)$ could be generated by sets $M \subseteq \dom f$ or by dual admissible elements.

\begin{definition}
Let $D \colon \R_+ \to \G(Z, C)$ be a function satisfying

(i) $D(\eps_2) \supseteq D(\eps_1)$ for all $\eps_1, \eps_2 \in \R_+$ with $0 < \eps_1 \leq \eps_2$, and

(ii) $C = D(0) = \bigcap_{\eps > 0} D(\eps)$.

A set $M \subseteq \dom f$ is called a $(D, \eps)$-solution of (P) if
\[
\inf f[M] \oplus D(\eps) \supseteq I(f),
\]
and each $x \in M$ is a minimizer of $f$. 
\end{definition}

A similar concept applies to supremum problems which can be useful in connection with duality.  If $M$ is a $(D, \eps)$-solution of (P), then
\[
\inf f[M] \oplus D(\eps) \supseteq I(f) \supseteq \inf f[M],
\]
i.e., $\inf f[M]$ trivially is an inner approximation of $I(f)$. 

The condition that elements of $M$ be minimizers for $f$ might be relaxed to any type of approximate minimizers, thus producing sets of $(D, \eps)$-solutions consisting of approximate minimizers. Similarly, the intersection in (ii) might be replaced by any type of set convergence which is sometimes useful if $C \subseteq D(\eps)$ is not satisfied for some (or all) $\eps > 0$.

It turned out that effective algorithms for vector and set optimization problems generate $(D, \eps)$-solutions, for example with
\[
D(\eps) = C - \eps c
\]
with some $c \in C\bs(-C)$, even $c \in \Int C$ under the assumption that the latter set is non-empty. This idea has been exploited with Benson's outer approximation algorithm as the building block, see \cite[Remark 4.10]{HamelLoehneRudloff13JOGO} and \cite[Proposition 4.8]{LoehneRudloffUlus13JOGO}. The obtained algorithms indeed produce approximations of the set-valued infimum for (linear, convex) vector optimization problems. In \cite{LoehneSchrage13Opt}, it is shown that the same idea can be used for minimizing a polyhedral set-valued function (i.e., a $\G(\R^q, C)$-valued function whose graph is a polyhedral set): The corresponding algorithm produces solutions in the sense of Definition \ref{DefSetSolution} and might be considered as the first ``true set-valued" algorithm. Its extension to non-polyhedral problems is highly desirable and another challenge for the future. 

%Algorithms of this type for linear vector optimization problems can also be found in \cite{EhrgottShao} .   

We note that a different algorithmic attempt for producing minimizers with respect to a set relation can be found in \cite{Jahn13R}. In particular, it provides a numerical test if two (compact) sets $A, B \subseteq Z$ are in relation with respect to $\lel_C \cap \leu_C$ (compare the closely related Section \ref{SubSecScalarization} of this survey and \cite{Jahn13JOTA-OF}). In the polyhedral case, this test can be implemented on a computer. An algorithm is given which produces minimizers of a set-valued function if the set of feasible point is finite, and a descent method \cite[Algorithm 4.1]{Jahn13R} for problem (P) generates feasible points which are minimal with respect to some subset of the set of feasible points. 

%%%New subsection
\subsection{Risk measures}
\label{SubSecRiskMeasures}

Set-valued risk measures shall serve as a prominent example of set-valued translative functions as discussed in Section \ref{SubSecTransFunc}. The framework will be the following. By a slight abuse of notation, $X$ is this section does not denote a linear space, but rather a random variable etc.

Let $\of{\Omega, \mathcal{F}_T, P}$ be a probability space. A multivariate random variable is a $P$-measurable function $X \colon \Omega \to \R^d$ for some positive integer $d \geq 2$. If $d=1$, the random variable is called univariate. Let us denote by $L^0_d =L^0_d\of{\Omega, \mathcal{F}_T, P}$ the linear space of the equivalence classes of all $\R^d$-valued random variables which coincide up to sets of $P$-measure zero ($P$-almost surely). As usual, we write
\[
\of{L^0_d}_+ =  \cb{X \in L^0_d \mid
    P\of{\cb{\omega \in \Omega \mid X\of{\omega} \in \R^d_+}} = 1}
\]
for the closed convex cone of $\R^d$-valued random vectors with $P$-almost surely
nonnegative components. An element $X \in L^0_d$ has components $X_1, \ldots, X_d$ in
$L^0 = L^0_1$. In a similar way, we use $L^p_d$ for the spaces of equivalence classes of $d$-dimensional random variables whose components are to the $p$-th power integrable (if $0 < p < \infty$) and essentially bounded (if $p = \infty$). The symbol $\One$ denotes the random variable in $L^0_1$ which has $P$-almost surely the value $1$.

Let $M \subseteq \R^d$ be a linear subspace. We set $M_+ = M \cap \R^d_+$ and assume $M_+ \neq \cb{0}$ in the following. 

\begin{definition}[\cite{HamelHeydeRudloff11MAFE}]
\label{DefRegRiskMea} A function $R \colon L^p_d \to
\mathcal{P}\of{M, M_+}$ is called a risk measure if it is

(R0) finite at $0 \in L^p_d$: $R\of{0} \neq \emptyset$, $R\of{0} \neq M$;

(R1) $M$-translative:
\begin{equation}
\label{EqTrans}
 \forall X \in L^p_d, \; \forall u \in M \colon R\of{X + u\One} = R\of{X} - u;
\end{equation}

(R2) $\of{L^p_d}_+$-monotone: $X^2 - X^1 \in \of{L^p_d}_+$ $\Rightarrow$ $R\of{X^2} \supseteq R\of{X^1}$.
\end{definition}

Set-valued risk measures are indeed recognized as $T$-translative if, within the notation of Section \ref{SubSecTransFunc}, $X=L^p_d$, $Z=M$, $C = M_+ \subseteq M$ and the linear operator $T \colon M \to L^p_d$ is defined by $Tu = -u\One$. This means that $T$ assigns to each $u\in M$ the random vector being constantly equal to $-u$. 

A financial interpretation is as follows. A multivariate random variable is understood as a model for an unknown future portfolio or payoff of $d$ assets where each component indicates the number of units of the corresponding asset in the portfolio. The elements of $R(X)$ are understood as deposits, to be given at initial time, which compensate for the risk of $X$. The collection of all such risk compensating initial portfolios is understood as a measure of the risk associated to $X$. Such deposits usually involve fewer assets than the original portfolio, for example cash in a few currencies. This motivates the introduction of the space $M$ which is called the space of eligible portfolios. A typical example is $M = \R^m \times \cb{0}^{d-m}$ for $1 \leq m \leq d$ with $m \ll d$.

The axiom (R1) roughly means that the risk of $X + u\One$ is the risk of $X$ reduced by $u$ whenever $u \in M$. Axiom (R2) also has a clear interpretation: if a random vector $Y\in L^p_d$ dominates another random vector $X\in L^p_d$, then there should be more possibilities to compensate for the risk of $Y$ (in particular cheaper ones) than for $X$.  Finiteness at zero means that there is an eligible portfolio which covers the risk of the zero payoff, but not all eligible portfolios {\f do}. Convexity is an important property as it allows to invoke diversification effects.
 
From $M$-translativity and $\of{L^p_d}_+$-monotonicity it follows that $R$ maps into $\mathcal{P}\of{M,M_+}$. Clearly, the image space of a closed convex risk measure is $\mathcal{G}\of{M,M_+}$.

If trading is allowed a market model has to be incorporated. Here, a one-period market with proportional transaction costs as in \cite{Kabanov99FS, Schachermayer04MF} is considered. It is given by closed convex cones $K_0$ and $K_T = K_T\of{\omega}$ with $\R_+^d \subseteq K_t(\omega) \subsetneq \R^d$ for all $\omega \in \Omega$ and $t \in \cb{0, T}$ such that $\omega \mapsto K_T\of{\omega}$ is $\mathcal F_T$-measurable. These cones, called solvency cones, include precisely the set of positions which can be exchanged into a nonnegative portfolio at time $0$ and $T$, respectively, by trading according to the prevailing exchange rates. We set $K_0^M := M \cap K_0 \subseteq M$ which is the cone containing the ``solvent" eligible portfolios.  The set
\begin{equation*}
L^p_d\of{K_T} = \cb{X \in L^p_d \mid
    P\of{\cb{\omega \in \Omega \mid X\of{\omega} \in K_T\of{\omega}}} = 1}
\end{equation*}
is a closed convex cone in $L^p_d$.

\begin{definition}[\cite{HamelHeydeRudloff11MAFE}]
A risk measure $R \colon L^p_d \to \mathcal{P}\of{M,M_+}$ is called market-compatible if it maps into $\mathcal{P}\of{M,K_0^M}$ and is
 $L^p_d\of{K_T}$-monotone, that is $X^2 - X^1 \in L^p_d\of{K_T}$  implies  $R\of{X^2} \supseteq R\of{X^1}$.
\end{definition}

Let $1 \leq p \leq \infty$. We consider the dual pairs $(L^p_d,L^q_d)$ with $\frac{1}{p}+\frac{1}{q} = 1$ and endow them with the norm topology if $ p<\infty$ and the $\sigma\of{L^\infty_d, L^1_d}$-topology on $L^\infty_d$ in the case $p = +\infty$, respectively.  The duality pairing is given by $(X, Y) \mapsto E[Y^TX]$ for $X \in L^p_d$, $Y \in L^q_d$. The adjoint operator $T^* \colon L^q_d \to M$ is given by $T^*Y=\Pr_M E\sqb{-Y}$ where $\Pr_M$ denotes the projection operator onto the linear subspace $M$.

The biconjugation theorem, Theorem \ref{ThmFenchelMoreau}, can be used to obtain a dual description of a closed convex market-compatible set-valued risk measure of the form
\begin{equation}
\label{EqDualRep0}
R(X)= R^{**}(X) = \bigcap_{Y\in L^q_d, \, v \in \of{K_0^M}^+\bs\cb{0}} \of{S_{(Y,v)}(X) + (-R^*)(Y,v)}
\end{equation}
with 
\[
S_{(Y,v)}(X)=\cb{u\in M \mid  v^Tu\ge E\sqb{Y^T X}}
\]
and $\of{K_0^M}^+=\cb{v\in M \mid \forall u\in K_0^M \colon v^Tu\ge 0}$.

Using the considerations of Section \ref{SubSecTransFunc} and taking into account that $L^p_d(K_T)$-monotonicity implies $(-R^*)\of{Y, v}=M$ if $-Y\not\in L^q_d\of{K^+_T}$ we get
\begin{equation}
\label{EqRiskMeasureConjugate1} (-R^*)\of{Y, v} =
    \left\{\begin{array}{ccc}
   \displaystyle \cl \bigcup_{X \in A_R} S_{\of{-Y, v}}\of{X} & : & -Y \in L^q_d\of{K^+_T}, \; v=\Pr_M E\sqb{-Y}  \\[1ex]
    M & : & \mbox{else.}
    \end{array}\right.
\end{equation}
Recall $A_R = \cb{X \in L^p_d \mid 0 \in R(X)}$ from Section \ref{SubSecTransFunc}.

The next lemma admits a change of variables from vector densities $Y$ to vector
probability measures $Q$. This allows a formulation of the dual representation result in
terms of probability measures as it is common in the scalar case. 

In the following, $\diag\of{w}$ with $w \in \R^d$ denotes the diagonal matrix with the
components of $w$ as entries in its main diagonal and zero elsewhere. Moreover,
$\mathcal{M}^P_d = \mathcal{M}^P_d\of{\Omega, \mathcal{F}_T}$ denotes the set of
all vector probability measures with components being absolutely continuous with respect
to $P$, i.e. $Q_i \colon \mathcal{F}_T \to \sqb{0,1}$ is a probability measure on
$\of{\Omega,\mathcal{F}_T}$ such that $\frac{dQ_i}{dP} \in L^1$ for $i = 1, \ldots, d$. 

\begin{lemma}
\label{LemDualTransform} (a) Let $Y \in L^q_d\of{K^+_T}$, $v =\Pr_M E\sqb{Y} \in \of{K_0^M}^+\bs\{0\}$. Then there are $Q \in \mathcal{M}^P_d$, $w \in K_0^+\bs M^\perp + M^\perp$ such that $\diag\of{w}\frac{dQ}{dP} \in L^q_d\of{K^+_T}$ and
$S_{\of{Y, v}} = {F}^M_{\of{Q, w}}$ with
\begin{equation}
\label{EqTransform}
 {F}^M_{\of{Q, w}}\sqb{X} = \cb{{\f z} \in M \mid w^TE^Q\sqb{X}
 \leq w^T{\f z}} = \of{E^Q\sqb{X} + H^+(w)} \cap M.
\end{equation}
(b) Vice versa, if $Q \in \mathcal{M}^P_d$, $w \in K_0^+\bs M^\perp + M^\perp$ such
that $\diag\of{w}\frac{dQ}{dP} \in L^q_d\of{K^+_T}$ then there is $Y \in
L^q_d\of{K^+_T}$ such that $v :=\Pr_M E\sqb{Y} \in \of{K_0^M}^+\bs\{0\}$ and
${F}^M_{\of{Q, w}} = S_{\of{Y, v}} $.
\end{lemma}

{\sc Proof.} See \cite{HamelHeydeRudloff11MAFE}. \pend

\medskip Let us denote the set of dual variables by
\[
\mathcal{W}^q =
    \cb{\of{Q, w} \in \mathcal{M}^P_d \times \R^d \mid
        w \in K_0^+\bs M^\perp + M^\perp, \;
        \diag\of{w}\frac{dQ}{dP} \in L^q_d\of{K^+_T}}.
\]
The preceding considerations lead to the following dual representation result.

\begin{theorem}
\label{ThmDualRep} A function $R \colon L^p_d \to \mathcal{G}\of{M,K_0^M}$ is a market-compatible
closed ($\sigma\of{L^\infty_d, L^1_d}$-closed if $p = \infty$) convex risk measure
if, and only if, there is a set $\mathcal W^q_R \subseteq \mathcal W^q$ such that
\begin{equation}
\label{ThmDualRep1}
\forall X \in L^p_d \colon R\of{X} =
 \bigcap_{\of{Q, w} \in \mathcal{W}^q}
 \sqb{(-\alpha_R)\of{Q, w} + \of{E^Q\sqb{-X} + H^+(w)}\cap M},
\end{equation}
where the function $-\alpha_R \colon \mathcal W^q \to \mathcal G(M, M_+)$ is defined by
\[
\forall \of{Q, w} \in \mathcal W^q_R \colon (-\alpha_R)\of{Q, w} = \cl\bigcup_{X' \in A_R} \of{E^Q\sqb{X'} + H^+(w)}\cap M
\]
and $(-\alpha_R)\of{Q, w} = M$ whenever $\of{Q, w} \in \mathcal W^q \bs\mathcal W^q_R$.
\end{theorem}

{\sc Proof.} See \cite{HamelHeydeRudloff11MAFE}. \pend

\medskip Lemma \ref{LemDualTransform} shows that the set $\mathcal{W}^\infty$ for $M = \R^d$ coincides with the set of so-called consistent price systems (or processes). Strictly consistent price systems  are crucial for market models with proportional transaction cost: In finite discrete time, the existence of such a price system is equivalent to the fundamental robust-no-arbitrage condition (see \cite{Schachermayer04MF} for conical and \cite{PennanenPenner10SIAMFM} for convex market models). Therefore, results like Theorem \ref{ThmDualRep}, derived with set-valued duality tools, fit nicely into the mathematical finance background: They produce the correct dual variables, and they yield formulas which look like the corresponding scalar ones.

%%%New subsection
\subsection{Comments on applications}

Duality for vector optimization problems is already discussed in Section \ref{SecCommDuality}. We add a few remarks about the linear case. It is an astounding fact that there still is no consensus on what to consider as the ``canonical" dual of a linear vector optimization problem. After early contributions of J. S. H. Kornbluth \cite{Kornbluth74ORQ}, H. Isermann \cite{Isermann78ZOR, Isermann78MCPS} and W. R\"odder \cite{Roedder77}, E. H. Ivanov and R. Nehse \cite{IvanovNehse85Opt} discuss five different duals for a given linear vector optimization problem which illustrates the ambiguity even in the ``simplest", i.e. linear, case. The difficulty is further illustrated by means of the examples in \cite{Brumelle81MOR} and \cite[Discussion after Theorem 8.13]{Jahn04Book}. A set-valued approach has been presented in \cite{HamelEtAl04JCA} and later compared to several ``vector-valued" duals in \cite{BotGradWanka12OL}. Compare also \cite{HeydeLoehneTammer09MMOR} and Dinh The Luc \cite{Luc11EJOR}.  We believe that this ambiguity and the mathematical difficulties that come with it are rooted in the non-totalness of the order: A two-player matrix game with vector payoffs is hardly in equilibrium since the decisions of the players also depend on their ``vertical preferences" (as well as on their guesses about the vertical preference of the opponent), i.e. the weight they put on the components of the payoff vectors. This topic, essentially the link between set-valued convex duality and games with vector payoffs (more general, payoffs which are not totally ordered), seems to be one of the most interesting open questions that can be derived from the material presented in this survey.

One advantage of the complete lattice approach presented in this survey is that the set-valued calculus deals with all ``vertical preferences", i.e. all reasonable scalarizations at the same time. This admits to re-discover the scalar duality results on a ``set level."

In 1998, Harold P. Benson \cite{Benson98JOGO}, \cite{Benson98JOTA} proposed an ``outer approximation algorithm'' to solve linear vector optimization problems ``in the outcome space'' (see also \cite{EhrgottLoehneShao12JOGO}, report already from 2007). Benson motivated this kind of solving a linear vector optimization problem by three practical reasons: First, the set of minimal elements in the outcome space $\R^q$ has a simpler structure than the set of minimizers in the decision space $\R^n$, because one usually has $q \ll n$. The second reason is that a decision maker prefers to base her decision on objectives rather than directly on a set of efficient decisions. The third argument is that many feasible points are mapped on a single image point which may lead to redundant information. 

Later it turned out that Benson's algorithm just computes solutions to $\eqref{PL}$ and $\eqref{DL}$ as defined above, see \cite{Loehne11Book, HamelLoehneRudloff13JOGO}. Therefore Benson's arguments motivate the solution concepts introduced in Section \ref{SubsecMinConcepts} from an application oriented viewpoint, compare also \cite{LoehneSchrage13Opt}. The geometric duality theory \cite{HeydeLoehne08SIAMOpt}, \cite{Heyde13JCA} briefly discussed in Section \ref{SubSecLVO} is a fundamental tool to develop dual algorithms for solving linear and convex vector optimization problems, see \cite{EhrgottLoehneShao12JOGO, HamelLoehneRudloff13JOGO, LoehneRudloffUlus13JOGO}. 

Set-valued risk measures have been introduced in \cite{JouiniMeddebTouzi04FS}. It contains a dual representation result for the sublinear case, basically a combination of the formulae \eqref{EqTranslativeDualSublin} and \eqref{EqDualRep0}. A more systematic development including the extension to the general convex case has been presented in \cite{HamelHeyde10SIFIN} while market compatibility is due to \cite{HamelHeydeRudloff11MAFE}. A link to depth-trimmed regions, yet another set-valued object from statistics, can be found in \cite{CascosMolchanov07FS}. Currently, the set-valued approach for evaluating multivariate risks is gaining more and more attention, see for example \cite{JoeLi11MCAP, MaggisLaTorre12ISOR, FeinsteinRudloff14} and also \cite{FarkasKochMedinaMunari14FS, FarkasKochMedinaMunari13ArX}. Applications of Benson's algorithm and its variants to financial problems can be found in \cite{LoehneRudloff11ArX, HamelRudloffYankova13MAFE, HamelLoehneRudloff13JOGO} and related approaches in \cite{RouxTokarzZastawniak08AAM, RouxZastawniak11ArX} as well as in \cite{Csirmaz13ArX}.

\bibliography{SetOptBib}

\end{document}